\documentclass{article}
\pdfoutput=1
\usepackage{geometry}
\usepackage[utf8]{inputenc}
\usepackage{amsfonts, amsmath, amssymb, amsthm}
\usepackage{physics}
\usepackage{graphicx}
\usepackage[hidelinks]{hyperref}
\usepackage{makecell}
\usepackage{float}
\usepackage[title]{appendix}

\usepackage{tikz}
\usetikzlibrary{arrows.meta}

\usepackage{xcolor}
\definecolor{myblue}{rgb}{0, 0, 0}

\usepackage[backend=bibtex,style=numeric,maxbibnames=99,sorting=none]{biblatex}
\renewbibmacro{in:}{}
\usepackage{mathtools}
\mathtoolsset{showonlyrefs}

\renewcommand{\norm}[1]{\left\lVert#1\right\rVert}

\DeclareMathOperator{\Div}{div}

\renewcommand{\dv}[2]{\frac{d#1}{d#2}}

\numberwithin{equation}{section}

\newtheorem{theorem}{Theorem}[section]
\newtheorem{corollary}[theorem]{Corollary}
\newtheorem{definition}[theorem]{Definition}
\newtheorem{lemma}[theorem]{Lemma}
\newtheorem{proposition}[theorem]{Proposition}
\newtheorem{remark}[theorem]{Remark}
\author{Anatoly Neishtadt and Alexey Okunev}
\date{}
\title{
On the phase change for perturbations of Hamiltonian systems with separatrix crossing\thanks{
The work was supported by the Leverhulme Trust (Grant No. RPG-2018-143).}}
\begin{document}
\maketitle

\abstract{
  { \color{myblue}
  We study the evolution of angular variable (phase) for general (not necessarily Hamiltonian) perturbations of Hamiltonian systems with one degree of freedom near separatrices of the unperturbed system. To this end, we use averaged system of order 2. We obtain estimates for the accuracy of order 2 averaged system near separatrices and use these estimates to prove a formula for the phase change when solutions of the perturbed system approach separatrices of the unperturbed system (such formula is known when the perturbation is Hamiltonian). As an application of this formula, we show that two natural definitions of probability of capture into different domains after separatrix crossing proposed by V.I. Arnold and D.V. Anosov lead to the same formula for this probability.
  }
}

{ \color{myblue}

\section{Introduction}
A standard object of perturbation theory are systems described by differential equations of the form
\begin{equation} \label{for_start}
  \dot x=A_0(x)+\varepsilon A_1(x,\varepsilon), \quad x \in \mathbb R^m.
\end{equation}
Here $\varepsilon$ is small parameter, when $\varepsilon=0$ we have \emph{unperturbed system}, and system~\eqref{for_start} with $0 < \varepsilon \ll 1$ is called \emph{perturbed system}.
An important particular case with numerous applications (e.g., in study of oscillations of mechanical systems, celestial mechanics, dynamics of charged particles) is when in the phase space of~\eqref{for_start} there is a domain filled by periodic trajectories of the unperturbed system. In this case one says that~\eqref{for_start} is \emph{one-frequency system} (or a system with one rotating phase).

For this case, the classical averaging method allows to describe the dynamics of the perturbed system with a high accuracy.  One can introduce new coordinates
\begin{equation}
  y \in \mathbb R^{m-1}, \qquad \varphi \in \mathbb R^1  \bmod 2\pi
\end{equation}
such that $y$ enumerates periodic trajectories of the unperturbed system and $\varphi$ is an angular variable (phase) on these trajectories. In order to approximately describe the evolution of the variable $y$ one should just average the rate of change of $y$ over the phase $\varphi$.
This allows to describe the behaviour of $y$ with accuracy $\sim \varepsilon$ on time intervals $\sim 1/\varepsilon$ (\cite{fatou1928mouvement},~\cite{bogolyubov61}). Such accuracy for $y$ is not enough to approximately describe the evolution of the phase $\varphi$, one should introduce the averaged system of order~$2$ to this end. Averaged system of order~$2$ allows to describe the behaviour of $y$ with accuracy $\sim \varepsilon^2$ and the dynamics of $\varphi$ with accuracy $\sim \varepsilon$
on time intervals $\sim 1/\varepsilon$ (\cite{bogolyubov61}).

This classical approach should be modified when the foliation of the phase space by periodic trajectories of the unperturbed system has singularities. This situation is general (i.e., cannot be destroyed by small perturbations) and is frequently encountered in applications. For example, effects of multiple separatrix crossings explain the emerging of chaos in the dynamics of charged particles in the tail of the Earth’s magnetosphere and play a key role in Wisdom’s theory of the origin of the Kirkwood gap in the asteroid belt at the $3:1$ resonance (references and more information on separatrix crossing can be found in \cite[§6.4.7]{arnold2007mathematical}).
It is typical for many applications (e.g., mentioned above) that such singularities appear as follows. The unperturbed system can be considered as a Hamiltonian system with $1$ degree of freedom for canonical variables $(q, p)\in \mathbb R^2$ and the Hamiltonian $H(p, q, z)$ depending on a vector parameter $z=(z_1,\ldots, z_k), \, k=m-2$. For all values of $z$ the Hamiltonian $H$ has a saddle point $C(z)$ and two separatrix loops passing through this point (cf. Figure~\ref{f:up} in Section~\ref{s:assumptions}). These separatrices split the phase portrait of $H$ into three domains filled by closed contour lines of $H$ (periodic trajectories of the unperturbed system). One of these domains is adjacent to both separatrix loops (we call it \emph{outer} domain) and the other two are adjacent to only one separatrix loop (\emph{inner} domains).
The unperturbed system is
\begin{equation} \label{e:unperturbed}
  \dot q = \pdv{H}{p}, \qquad \dot p = - \pdv{H}{q}, \qquad \dot z = 0
\end{equation}
and the perturbed system~\eqref{for_start} takes the form
\begin{align} \label{e:perturbed}
\begin{split}
  \dot q = \pdv{H}{p} + \varepsilon f_q(p, q, z, \varepsilon), \qquad
  \dot p = - \pdv{H}{q} + \varepsilon f_p(p, q, z, \varepsilon), \qquad
  \dot z = \varepsilon f_z(p, q, z, \varepsilon).
\end{split}
\end{align}
In this case we have $y = (I, z)$ (here $I$ is the action variable of the unperturbed system) and $\varphi$ is the angle from the pair of action-angle variables.
For the perturbed system, $I$ and $z$ are slow variables (they change with speed $\sim \varepsilon$) and $\varphi$ is a fast variable. Solutions of the perturbed system may cross separatrices of the unperturbed system.

Separatrix crossing is much better studied when the perturbed system is Hamiltonian. We will discuss this case further in the introduction, let us now focus on separatrix crossing for arbitrary perturbations~\eqref{e:perturbed}. Note that if in a physical problem there is some kind of friction (e.g., tidal friction in celestial mechanics), then the perturbed system is non-Hamiltonian.

The right hand sides of the equations of motion in the variables $I, \varphi$ are singular at the separatrices of $C$. Near the separatrices the angle variable $\varphi$  behaves badly: the unperturbed frequency is small  while  the rate of change of $\varphi$ along the solutions of the perturbed system  is unbounded.  Many functions used in the averaging method are unbounded.
The averaging method can still be used, but estimates for these functions are required.
One can use the variable $h = H(p, q, z) - H(C(z), z)$ instead of the action $I$, this is more convenient near the separatrices. Suitably, $h$ measures how far a point in the phase space is from the separatrices, $h=0$ on the separatrices. Averaged system can be written using $h$ instead of $I$, then the averaged system in the outer domain can be glued with the averaged systems in the inner domains by the line $h=0$.
Accuracy of averaging method $O(\varepsilon |\ln \varepsilon|)$ over times $\sim \varepsilon^{-1}$ holds (\cite{neishtadt17} and references therein) for solutions crossing the separatrices
for most initial data, except a set of measure $O(\varepsilon^r)$. Here $r$ can be taken as large as needed, but larger $r$ give worse constant in the $O$-estimate for accuracy of averaging method.

Crossing of separatrices leads to a remarkable scattering of trajectories. Suppose a trajectory starts in the outer domain and approaches the separatrices. Then it moves into one of the inner domains, we call the choice of this inner domain the \emph{outcome} of separatrix crossing. Initial data with different outcomes alternate in the phase space with a step $\sim \varepsilon$. For fixed initial data, different outcomes alternate when $\varepsilon$ changes.
Because of this, the scattering of trajectories on the separatrices does not have a deterministic description when $\varepsilon \to 0$. One should consider different outcomes as random events and give a definition of probabilities of these events. There are two natural definitions of such probabilities. Definition by Arnold~\cite{arnold63} uses how different outcomes alternate in the phase space, while Anosov's definition uses how they alternate for fixed initial data when $\varepsilon \to 0$ (see Section~\ref{s:probability}).
Formulas for the probabilities of different outcomes were proven in~\cite{neishtadt17} for Arnold's definition, and it was suggested that for Anosov’s definition the formulas are the same.
Stochastic perturbations of some particular classes of systems~\eqref{e:perturbed} were studied in~\cite{wolansky1990limit} and~\cite{freidlin1998random}. In this case there are true probabilities of moving to different domains, and zero-noise limit of such probabilities is given by the same formulas.

In this paper we establish estimates on order 2 averaging method near the separatrices. This allows to track the evolution of the phase $\varphi$ with accuracy $O(\varepsilon h^{-1} |\ln^{-1} h|)$ (Theorem~\ref{t:avg-2} below).
Such estimates existed before only when the perturbation is Hamiltonian. For non-Hamiltonian perturbations, let us mention~\cite{bourland1991averaging}, where a formula for the evolution of phase was written using averaging method for a particular case of motion in one-dimensional slowly time dependent potential with an additional dissipative perturbation. Estimates for accuracy of this formula were not obtained.

The phase $\varphi$ is not defined on the separatrices, a special parameter is used in literature to describe the behaviour of fast variables $p, q$ immediately before the moment of separatrix crossing (it is called \emph{crossing parameter} in~\cite{cary1989phase}, and \emph{pseudo-phase} in~\cite{neishtadtVasiliev2005}).
We use our estimates for the evolution of the phase to prove a formula~\eqref{e:phase} for the pseudo-phase (or, alternatively, one may think that this formula describes the phase change when approaching the separatrices). This formula is similar to the formulas~\cite{cary1989phase} and~\cite{neishtadtVasiliev2005} for Hamiltonian perturbations.
For non-Hamiltonian perturbations, there is formula~\cite{bourland1994connection} (see also~\cite{bourland1990separatrix}) for the aforementioned particular case of motion in slowly time-dependent one-dimensional potential with an additional dissipative perturbing force.
Unlike the works for Hamiltonian perturbations and our result, accuracy of this formula is not estimated; validity of the formula is justified by a comparison with the results of numerics.

Finally, as a corollary of our formula for pseudo-phase, we prove that the formula for probability given by Anosov's definition is indeed the same as for Arnold's definition.

Let us now discuss the case when the perturbed system is Hamiltonian.
This class is much better studied, and our results for the general case do not give new corollaries.
The perturbed system may be a slow-fast Hamiltonian system, or (this could be treated as a particular case of slow-fast Hamiltonian systems) a system with the Hamiltonian slowly depending on the time.
In both cases the action $I$ is constant for the solutions of the averaged system, thus $I$ is called the \emph{adiabatic invariant}.
One can also introduce \emph{improved adiabatic invariant} $J$ (e.g., \cite[§6.4.4]{arnold2007mathematical}), it is preserved with accuracy $O(\varepsilon^2)$ far from separatrices over times of order $\varepsilon^{-1}$ for the one degree of freedom case.
Separatrix crossing leads to a jump of the improved adiabatic invariant of order $\varepsilon \ln \varepsilon$.
There are formulas for the value of this jump~(\cite{timofeev1978constancy}, \cite{cary1986adiabatic}, \cite{neishtadt86}, \cite{neishtadt1987change}), and this value depends on the pseudo-phase.
Formulas for pseudo-phase were obtained (using the averaging method)
in~\cite{cary1989phase} for Hamiltonian systems with one degree of freedom and slow time dependence and in~\cite{neishtadtVasiliev2005} for slow-fast Hamiltonian systems with one degree of freedom corresponding to fast motion.
Formulas for pseudo-phase together with the formulas for the change of the adiabatic invariant allow to study trajectories with multiple separatrix crossings. Let us mention the remarkable existence of stability islands~\cite{neishtadt1997stable}, \cite{vasiliev2007stability}. Our work is a step towards the study of multiple separatrix crossings for non-Hamiltonian perturbations. The remaining ingredient is an analogue of the formulas for the change of adiabatic invariant, hopefully, it will be obtained in future works.
For particular cases such analogues are suggested in~\cite{bourland1994connection, bourland1990separatrix}.

The structure of the paper is as follows. In Section~\ref{s:order-2-overview} we briefly discuss order 2 averaged system, this allows us to state the results in Section~\ref{s:results}. This is followed by Section~\ref{s:avg-2-details} with detailed description of order 2 averaging and formulas for the coefficients of order 2 averaged system. In the remaining sections we prove these results, overview of the proofs with a plan of the rest of the paper is presented in Section~\ref{s:proofs-overview}.

\section{Averaged system of order 2 (overview)} \label{s:order-2-overview}
The first-order averaged system is given by the formulas
\begin{align}
\begin{split} \label{e:avg-1-order}
  \hat h' &= \hat{f}_{h, 1}(\hat h, \hat w), \\
  \hat w' &= \hat{f}_{w, 1}(\hat h, \hat w). \\
\end{split}
\end{align}
Here $\hat h \in \mathbb R$ tracks the evolution of $h$, $\hat w \in \mathbb R^k$ tracks the evolution of $z$ and $\psi '$ denotes the derivative of $\psi$ with respect to the slow time $\tau = \varepsilon t$. The functions $\hat f_{h, 1}$ and $\hat f_{w, 1}$ are averages over $\varphi$ of the rates of change of $h$ and $z$, respectively, divided by $\varepsilon$ (this is repeated more formally in Section~\ref{s:avg-2-details}).

Writing the averaged system of order 2 explicitly is cumbersome, we postpone it until Section~\ref{s:avg-2-details}. Let us briefly discuss this system here so that we can state the results of this paper.
The separatrices split the phase space of~\eqref{e:perturbed} into three domains, let us discuss the averaged system in one of these domains $G$.
It is obtained by the following standard procedure:
\begin{enumerate}
  \item The perturbed system is rewritten in the energy-angle chart $(h, \varphi, w)$.
    Here
    \begin{itemize}
      \item $h = |H(p, q, z) - H(C, z)|$,
      \item $\varphi$ is the angle variable from the pair of action-angle variables in $G$,
      \item $w = z$ (this helps to distiguish the energy-angle chart and the $(p, q, z)$ chart).
    \end{itemize}
  \item One finds a coordinate change $(h, \varphi, w) \mapsto (\hat h, \hat \varphi, \hat w)$ such that after this coordinate change the $\varphi$-dependent terms on the right-hand side of the perturbed system are small (with order $\varepsilon^3$ for $\dot h$ and $\dot w$ and with order $\varepsilon^2$ for $\dot \varphi$ far from the separatrices; however, these terms might be unbounded near the separatrices). This coordinate change is described by functions that we will denote by $u$ with different lower indices. The functions $u_{h, 1}$ and $u_{w, 1}$ are particularly important because they are needed to connect given initial data of the perturbed system with the corresponding initial data of averaged system of order 2. We have $u_{h, 1}, u_{w, 1} = O(1)$ even near the separatrices, unlike other functions $u$ that may be unbounded near the separatrices.
  Formulas for $u_{h, 1}$ and $u_{w, 1}$ can be found in Section~\ref{s:avg-2-details}. They depend on $\varepsilon$, we denote by $u^0_{h, 1}$ and $u^0_{w, 1}$ their values when $\varepsilon = 0$.
  \item One drops the $\varphi$-dependent terms on the right-hand side and uses the slow time $\tau = \varepsilon t$ instead of the normal time $t$, this gives the averaged system of order 2.
\end{enumerate}
The averaged system of order 2 has the form
\begin{align}
\begin{split} \label{e:avg-form}
  \hat h' &= \hat{f}_{h, 1}(\hat h, \hat w) + \varepsilon \hat{f}_{h, 2}(\hat h, \hat w), \\
  \hat w' &= \hat{f}_{w, 1}(\hat h, \hat w) + \varepsilon \hat{f}_{w, 2}(\hat h, \hat w), \\
  \hat \varphi' &= \varepsilon^{-1} \omega(\hat h, \hat w) + \omega_1(\hat h, \hat w).
\end{split}
\end{align}
Again, $\psi '$ denotes the derivative of $\psi$ with respect to the slow time $\tau = \varepsilon t$. The functions $\hat f_{h, 1}$ and $\hat f_{w, 1}$ here are as in the first order averaged system.
Formulas for the other functions $\hat f_{*, *}$ and the function $\omega_1$ can be found in Section~\ref{s:avg-2-details}. One should note that they may be unbounded near the separatrices, estimates on these functions are gathered in Table~\ref{t:est} below.

\section{Results} \label{s:results}
\subsection{Assumptions} \label{s:assumptions}
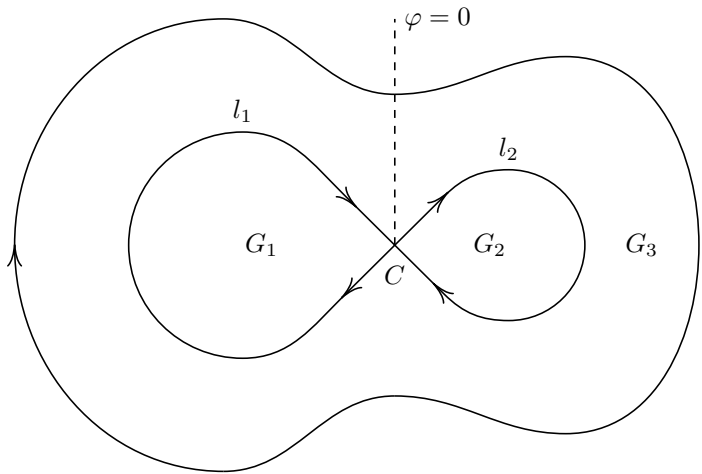
\begin{figure}[H]
  \centering
  \begin{tikzpicture}[scale=0.5]
    \draw[-{To[length=3mm,width=2mm]},semithick](1.1, -1.1)--(1,-1);
    \draw[-{To[length=3mm,width=2mm]},semithick](-1.4, -1.4)--(-1.41,-1.41);
    \draw[-{To[length=3mm,width=2mm]},semithick](-1.1, 1.1)--(-1, 1);
    \draw[-{To[length=3mm,width=2mm]},semithick](1.4, 1.4)--(1.41, 1.41);
    \draw[-{To[length=3mm,width=2mm]},semithick][semithick] (0, 0)
      to [out=45, in=45] (1.3, 1.3) to[out=45, in=180] (3, 2) node[above]{$l_2$}
      to[out=0, in=90] (5, 0)
      to[out=-90, in=0] (3, -2) to[out=180, in=-45] (1.3, -1.3) to[out=135, in=135] cycle;
    \draw[semithick] (0, 0) to[out=135, in=135] (-1.3, 1.3)
      to[out=135, in=0] (-4, 3) node[above]{$l_1$} to[out=180, in=90] (-7, 0)
      to[out=-90, in=180] (-4, -3) to[out=0, in=-135] (-1.3, -1.3) to[out=45, in=45] cycle;

    \draw[semithick] (0, 4) to[out=0, in=180] (4.5, 5) to[out=0, in=90] (8, 0)
      to[out=-90, in=0] (4.5, -5) to[out=180, in=0] (0, -4);
    \draw[semithick] (0, 4) to[out=180, in=0] (-4.5, 6) to[out=180, in=90] (-10, 0);
    \draw[-{To[length=3mm,width=2mm]},semithick] (-4.5, -6) to[out=180, in=-90] (-10, 0.03);
    \draw[semithick] (-4.5, -6) to[out=0, in=180] (0, -4);

    \draw (0, -0.3) node[below]{$C$};
    \draw (-3.5, 0) node{$G_1$};
    \draw (2.5, 0) node{$G_2$};
    \draw (6.5, 0) node{$G_3$};

    \draw[dashed, semithick] (0, 0) -- (0, 6) node[right]{$\varphi = 0$};
  \end{tikzpicture}
  \caption{The unperturbed system.}
  \label{f:up}
\end{figure}
\begin{itemize}
  \item Denote by $f_h = \pdv{h}{p} f_p + \pdv{h}{q} f_q + \pdv{h}{z} f_z$ the rate of change of $h$ for the perturbed system, divided by $\varepsilon$. Denote the separatrices by $l_1$ and $l_2$, they split the phase space of the unperturbed system into three domains that we denote $G_1, G_2, G_3$. Suppose that $G_1$ and $G_2$ are bounded by the separatrices $l_1$ and $l_2$, respectively, while $G_3$ is bounded by $l_1 \cup l_2$ (cf. Figure~\ref{f:up}).
  Set
  \begin{equation}
    \Theta_i(z) = -\oint_{l_i} f_h(p(t), q(t), z, 0) dt, \qquad i=1,2
  \end{equation}
  (here $t$ is the time for the unperturbed system). Let $\Theta_3(z) = \Theta_1(z) + \Theta_2(z)$.
  The results below describe trajectories in some domain $G_i$ approaching separatrix/separatrices.
  We assume that $\Theta_i > 0$ for all considered $z$. This is not a restrictive assumption: during one wind around separatrix/separatrices the value of $h$ changes by approximately $\varepsilon \Theta_i$, thus $\Theta_i < 0$ means moving away from the separatrices.
  \item There exists an open domain $D \subset \mathbb R^{k+2}_{p, q, z}$ such that $H(p, q, z)$ is analytic in $D$ and $f(p, q, z, \varepsilon)$ is $C^2$ in $D \times [0, \varepsilon_0]$ for some $\varepsilon_0 > 0$; here $f=(f_p, f_q, f_z)$. The separatrices lie in $D$ for all encountered $z$: we have
  $\cup_{z \in \pi_z(D)} (l_1(z) \cup l_2(z) \cup C(z)) \subset D$,
  where $\pi_z$ denotes the projection $(p, q, z) \mapsto z$. For each $z$ the domain $D_z = \{(p, q): (p, q, z) \in D\}$ is foliated by the level sets of $H$. This means that we can write $(h, z) \in D$ (with a slight abuse of notation).
  \item Initial data $p_0, q_0, z_0$ is such that the corresponding solution of first-order averaged system (i.e., with initial data $h_0 = h(p_0, q_0, z_0), w_0 = z_0$) remains in $D$ until it crosses the separatrices (i.e., reaches $h=0$).
  \item The estimates obtained in this paper are asymptotic with respect to $\varepsilon$: we prove that they hold if $\varepsilon$ is sufficiently small.
  \item
  For each $z$ the transversal $\varphi = 0$ on the plane $(p, q)$ is tangent to the bisector of an angle between the separatrices with the vertex at $C$. In $G_1$ and $G_2$ there is only one such angle, so this bisector is uniquely defined. In $G_3$ there are two such angles, and the choice of the bisector in $G_3$
  should match with the enumeration of the separatrices $l_1$ and $l_2$: the separatrix $l_2$ should correspond to $0 < \varphi < \pi$ and $l_1$ to $\pi < \varphi < 2\pi$ (cf. Figure~\ref{f:up}).
  This transversal is $C^2$.

\end{itemize}
\subsection{Order 2 averaging near separatrices}
Consider some initial data in $G_i$, $i=1,2,3$.
As we have $\Theta_i > 0$, solutions of the perturbed system approach separatrices.
Solution of order $1$ averaged system cross separatrices, i.e., reach $h=0$; we prove (Lemma~\ref{l:avg-cross} below) that same holds for solutions of order $2$ averaged system.
Our main result about order 2 averaging is the following estimate for the accuracy of order $2$ averaging method until $O(\varepsilon)$-close to the separatrices. Unlike order $1$ averaging method, order $2$ averaging method allows to track the evolution of the phase $\varphi$.

We will use the notation $X = (p, q, z)$ and will denote by $h(X), \varphi(X)$ and $w(X)$ the components of $X$ rewritten in the energy-angle chart.
\begin{theorem} \label{t:avg-2}
  There exists $C > 0$ such that the following holds.
  Consider a solution $X(t)$ of the perturbed system~\eqref{e:perturbed} with initial data $X(0)$. Recall that $\tau = \varepsilon t$.
  Let $(\hat h(\tau), \hat w(\tau), \hat \varphi(\tau))$ be the solution of averaged system of order 2 given by~\eqref{e:avg-form} with initial data
  \begin{equation}
    (\hat h(0), \hat w(0), \hat \varphi(0))
    =
    (h_0 - u_{h, 1}(h_0, w_0, \varphi_0), w_0 - u_{w, 1}(h_0, w_0, \varphi_0), \varphi_0),
  \end{equation}
  where $h_0=h(X(0))$, $w_0=w(X(0))$, $\varphi_0=\varphi(X(0))$. Then
  \begin{align}
  \begin{split}
    \varphi(t) &= \hat \varphi(\tau)
      + O(\varepsilon h^{-1} \ln^{-1} h) \qquad \mod 2 \pi, \\
      h(t) &= \hat h(\tau) + \varepsilon u^0_{h, 1}(\hat h(\tau), \hat w(\tau), \hat \varphi(\tau)) + O(\varepsilon^2 h^{-1}), \\
      w(t) &= \hat w(\tau) + \varepsilon u^0_{w, 1}(\hat h(\tau), \hat w(\tau), \hat \varphi(\tau)) + O(\varepsilon^2 h^{-1}) \\
  \end{split}
  \end{align}
  if $t$ is such that $\hat h(\tau) > C \varepsilon$ for all $\tau \in [0, \varepsilon t]$.
  In the error terms in these estimates $h$ is a shorthand for $\hat h(\tau)$.
\end{theorem}
\begin{remark} \label{r:any-transversal}
  In Theorem~\ref{t:avg-2} it is not important that the transversal $\varphi=0$ is tangent to the bisector of an angle between the separatrices. We can take any transversal $\Gamma'(z)$ to the union of the separatices passing through $C(z)$ for all $z$ that lies in the same angle between the separatices as the bisector. Alternatively, we can take any transversal $\Gamma'(z)$ to one of the separatices with $C(z) \not \in \Gamma'(z)$ for all $z$.
\end{remark}
\noindent This remark is proved in Section~\ref{s:approx-proof}.

\begin{lemma} \label{l:avg-cross}
  Solutions of averaged system of order $2$ cross the separatrices (i.e., reach $\hat h = 0$) for small $\varepsilon$.
\end{lemma}
\noindent
This lemma holds because we have $\hat f_{h, 1} < 0$ in~\eqref{e:avg-form} due to $\Theta_i > 0$. However, it needs to be proved, because the $\varepsilon \hat f_{h, 2}$ term in~\eqref{e:avg-form} might destroy the fact that $h$ decreases for the averaged system of order 2. We prove that it is not the case by providing a good estimate for this term.

\subsection{Formula for pseudo-phase} \label{s:formula}
The parameter \emph{pseudo-phase} (we use the name from~\cite{neishtadtVasiliev2005}, another name used in~\cite{cary1989phase} is \emph{crossing parameter}) describes the behaviour of the fast variables $p, q$ immediately before the moment of separatrix crossing.
This parameter is important for several reasons. First of all, one can match the values of pseudo-phase when approaching the separatrices and when moving away\footnote{Moving away from the separatrices can be studied in the same way as approaching the separatrices, one just needs to reverse the time.} from them to track the evolution of $\varphi$ after separatrix crossing. Secondly, when the perturbation is Hamiltonian, the $O(\varepsilon |\ln \varepsilon|)$ jump of adiabatic invariant after separatrix crossing is determined by the pseudo-phase, and we expect that in non-Hamiltonian case the $O(\varepsilon |\ln \varepsilon|)$ jump of the slow variables of the perturbed system caused by separatrix crossing will also depend on the pseudo phase (however, as far as we know, there are no formulas for this jump yet).

Let us state the definition of the pseudo-phase.
Suppose that we study approach to the separatrices in $G_i$. If $i=3$, it is convenient to assume that we have $\Theta_1, \Theta_2 > 0$ in addition to $\Theta_3 > 0$ (we discuss what happens if this does not hold in Remark~\ref{r:h-minus-one} below). We are given a solution $X(t)$ of the perturbed system~\eqref{e:perturbed}. Denote $h(t) = h(X(t))$, $v=(h, z)$ and $v(t) = v(X(t))$. Let $h_{-1}$ be the value of $h$ at the last crossing of the transversal $\varphi=0$
before $X(t)$ crosses the separatrices of the unperturbed system.
Let us consider a solution $(\hat h(\tau), \hat z(\tau), \hat \varphi(\tau))$ of order 2 averaged system with initial data determined by $X(0)$ as in Theorem~\ref{t:avg-2}, set $\hat v(\tau) = (\hat h(\tau), \hat w(\tau))$. Denote by $w_*$ the value of $\hat w$ when $\hat h$ reaches $0$ (this always happens by Lemma~\ref{l:avg-cross}).
Alternatively, we can pick as $w_*$ the value of $w$ when the solution of the \emph{first-order} averaged system with initial data $X(0)$ reaches $h=0$, these two quantities differ by $O(\varepsilon)$ by Remark~\ref{r:w-star} below.
Set $\Theta_{i*} = \Theta_i(w_*)$. The pseudo-phase $\xi$ is defined by the formula~\cite{neishtadtVasiliev2005}
\begin{equation}
  \xi = \frac{h_{-1}}{\varepsilon \Theta_{i*}}.
\end{equation}
As $h$ decreases by approximately $\varepsilon \Theta_{i*}$ during one turn, we have $0 \le h_{-1} < \varepsilon \Theta_{i*} + O(\varepsilon^{3/2})$, this is proved below in Remark~\ref{r:est-h-m-1}.

We will assume
$c_1 \varepsilon^{3/2} < h_{-1} < \varepsilon \Theta_{i*}$
for large enough\footnote{see Lemma~\ref{l:delta-h-near-sep} below} $c_1 > 0$. This holds for most initial conditions, what happens for other initial data is discussed in Remark~\ref{r:h-minus-one} below.
Denote $u_* = \frac 1 4 (\Theta_1(w_*) - \Theta_2(w_*))$ for $i=3$ and $u_* = 0$ for $i=1, 2$.
Then we have
\begin{equation} \label{e:phase}
  \xi
  = \bigg\{ \frac{1}{2\pi}
    \bigg(
      \varphi_0 + \frac 1 \varepsilon
      \int_{\tau=0}^{\tau_*} \big(
        \omega(\hat{v}(\tau)) + \varepsilon \omega_1(\hat{v}(\tau))
      \big) d\tau
    \bigg)
    + \frac{u_*}{\Theta_{3*}}
    + O(\varepsilon^{1/3} \ln^{1/3} \varepsilon)
  \bigg\}.
\end{equation}
Here the curly brackets $\{ \cdot \}$ denote the fractional part. This formula is proved in Section~\ref{s:formula-proof}.

\begin{remark} \label{r:h-minus-one}
  Denote the values of $h$ at the crossings of $\varphi = 0$ before $h_{-1}$ by $h_{-2}, h_{-3}, \dots$.
  We have assumed earlier that
  $c_1 \varepsilon^{3/2} < h_{-1} < \varepsilon \Theta_{i*}$.
  If $h_{-1} < c_1 \varepsilon^{3/2}$,  the right-hand side of~\eqref{e:phase} gives $\big\{ \frac{h_{-2}}{\varepsilon \Theta_{i*}} \big\}$ and if $h_{-1} > \varepsilon \Theta_{i}$, it gives $\big\{ \frac{h_{-1}}{\varepsilon \Theta_{i*}} \big\}$.

  We have also assumed that for $i=3$ we have $\Theta_1, \Theta_2 > 0$. If one of $\Theta_1$, $\Theta_2$ is less than zero with $\Theta_3 = \Theta_1 + \Theta_2 > 0$, the last transversal crossing can happen for $h > \varepsilon \Theta_3$.
  In this case we should find the first $k \in \mathbb N$ such that $h(t) > c_1 \varepsilon^{3/2}$ during all the time before the moment corresponding to $h_{-k}$.
  Then the right-hand side of~\eqref{e:phase} gives $\big\{ \frac{h_{-k}}{\varepsilon \Theta_{3*}} \big\}$.
\end{remark}

\subsection{Probabilities} \label{s:probability}
Assume $\Theta_1, \Theta_2 > 0$. Then a trajectory starting in $G_3$ may be captured into $G_1$ or $G_2$ after separatrix crossing and initial data corresponding to different outcomes are finely mixed in $G_3$. It is natural to consider captures in $G_1$ and $G_2$ as random events with some probabilities.
One natural definition of the probability of capture is stated in~\cite{arnold63}.
Denote by $U^\delta(X_0)$ the $\delta$-neighborhood of the point $X_0 = (p_0, q_0, z_0)$ in the action-angle coordinates:
\begin{equation}
  U^\delta(p_0, q_0, z_0) = \{ p, q, z:
    \norm{z-z_0} < \delta, \; |I(p, q, z) - I_0| < \delta, \;
    |\varphi(p, q, z) - \varphi_0| < \delta
  \}.
\end{equation}
Here $I_0 = I(X_0)$ and $\varphi_0 = \varphi(X_0)$. Denote by $U^\delta_1(X_0)$ and $U^\delta_2(X_0)$ the subsets of initial data in $U^\delta(X_0)$ captured in $G_1$ and $G_2$, respectively. Let $m$ denote the Lebesgue measure in $\mathbb R^{k+2}_{p, q, z}$. Fix some initial data $X_0$.
\begin{definition}[V.I. Arnold,~\cite{arnold63}] \label{d:prob-Arnold}
  The probability of capture in $G_j$, $j=1, 2$ is
  \begin{equation}
    P_j(X_0) = \lim_{\delta \to 0} \lim_{\varepsilon \to 0} \frac{m(U^\delta_j(X_0))}{m(U^\delta(X_0))}.
  \end{equation}
\end{definition}
\noindent For this definition, the following formula for the probability of capture is proved in~\cite{neishtadt17}:
\begin{equation}
  P_j(X_0) = \frac{\Theta_j(w_*)}{\Theta(w_*)}, \qquad j = 1, 2.
\end{equation}
Here $w_*$ is the value of $w$ when the solution of averaged system of order~$1$ with initial data $(h(X_0), z(X_0))$ reaches $h=0$. This formula for probability can also be obtained as a corollary of our formula for pseudo-phase~\eqref{e:phase}, see Remark~\ref{r:prob-Arnold} below.

Another way of defining the probability of capture was suggested by~D.V.~Anosov\footnote{This was a comment in a meeting of the Moscow Mathematical Society, this definition was first discussed in literature in~\cite{neishtadt17}.}.
Denote by $\mu$ the Lebesgue measure on $\mathbb R$.
Let us fix the initial data $X_0 = (p_0, q_0, z_0)$.
Denote
\begin{equation} \label{e:prob-capture}
  U_j(\varepsilon_0) = \big\{
    \varepsilon \in (0, \varepsilon_0)
    :
    \text{ the trajectory of } X_0 \text{ is captured into } G_j
  \big\}.
\end{equation}
\begin{definition}[D.V. Anosov] \label{d:prob-Anosov}
  The probability of capture into $G_j$, $j=1,2$ is
  \begin{equation}
    \lim_{\varepsilon_0 \to 0} \mu(U_j(\varepsilon_0)) / \varepsilon_0.
  \end{equation}
\end{definition}
\noindent It was suggested in~\cite{neishtadt17} that for this definition the probability of capture into $G_j$ is also $\frac{\Theta_j}{\Theta_3}$. Using~\eqref{e:phase}, we can prove this statement:
\begin{proposition} \label{p:prob-Anosov}
  We have
  \begin{equation} \label{e:est-prob-capture}
    \frac{\mu(U_j(\varepsilon_0))}{\varepsilon_0} = \frac{\Theta_j(w_*)}{\Theta_3(w_*)} + O(\varepsilon_0^{1/3} \ln^{1/3} \varepsilon_0).
  \end{equation}
\end{proposition}

}

\section{Details on order 2 averaging method} \label{s:avg-2-details}
\subsection{Energy-angle variables} \label{s:angle}
We are interested in trajectories starting in one of the domains $G_1, G_2, G_3$ (let us denote this domain $G_i$) and approaching the separatrix/separatrices of the unperturbed system.

Let us consider the action-angle variables $I, \varphi$; $\varphi \in [0, 2\pi)$ for the unperturbed system in the domain $G_i$.
We will assume that $\varphi = 0$ corresponds to a specific transversal $\Gamma(z)$ that is chosen in Section~\ref{s:Moser}. It will be tangent to the bisector of the angle between the separatrices.
Denote $h(p, q, z) = |H(p, q, z) - H(C, z)|$. Denote $w = z$. We will use the "energy-angle" variables $h, w, \varphi$. The notation $w=z$ is useful in order to distinguish $\pdv{}{z}$, which is taken for fixed $p, q$, and $\pdv{}{w}$, which is taken for fixed $h, \varphi$.
In these variables the unperturbed system~\eqref{e:unperturbed} is written as
$\dot h = 0, \dot w = 0, \dot \varphi = \omega(h, w)$.
Denote by $T(h, w) = \frac{2\pi}{\omega}$ the period of the unperturbed system.
We will sometimes use the time $t$ passed from the last crossing of the transversal $\varphi = 0$ instead of $\varphi$. We have $t = \frac{\varphi T}{2 \pi}$.

Denote by $f_h, f_w, f_\varphi$ the components of $f$ in the energy-angle variables:
$f_y = f_q \pdv{y}{q} + f_p \pdv{y}{p} + f_z \pdv{y}{z}$ for $y = h, \varphi$ and $f_w = f_z$.
Then the perturbed system~\eqref{e:perturbed} is written as
\begin{align} \label{e:init}
\begin{split}
  \dot h &= \varepsilon f_h(h, w, \varphi, \varepsilon), \\
  \dot w &= \varepsilon f_w(h, w, \varphi, \varepsilon), \\
  \dot \varphi &= \omega(h, w) + \varepsilon f_\varphi(h, w, \varphi, \varepsilon).
\end{split}
\end{align}

\subsection{Averaging chart} \label{s:varchange}
We start with the system~\eqref{e:init}.
In line with the general approach of the averaging method, let us find a change of variables\footnote{Here the functions $u_{*. *}$ depend on the small parameter $\varepsilon$, this is a convenient way to deal with the case when the perturbation $f_*$ in~\eqref{e:perturbed} depends on $\varepsilon$. If $f$ does not depend on $\varepsilon$, we can take $u_{*. *}$ independent of $\varepsilon$.}
\begin{align}
\begin{split} \label{e:varchange}
  h &= \overline h + \varepsilon u_{h, 1}(\overline h, \overline w, \overline \varphi, \varepsilon) + \varepsilon^2 u_{h, 2}(\overline h, \overline w, \overline \varphi, \varepsilon), \\
  w &= \overline{w} + \varepsilon u_{w, 1}(\overline h, \overline w, \overline \varphi, \varepsilon) + \varepsilon^2 u_{w, 2}(\overline h, \overline w, \overline \varphi, \varepsilon), \\
  \varphi &= \overline \varphi + \varepsilon u_{\varphi, 1}(\overline h, \overline w, \overline \varphi, \varepsilon)
\end{split}
\end{align}
that transforms~\eqref{e:init} to the following form:
\begin{align}
\begin{split} \label{e:init-avg}
  \dot{\overline h} &= \varepsilon \overline f_{h, 1}(\overline h, \overline w, \varepsilon) + \varepsilon^2 \overline f_{h, 2}(\overline h, \overline w, \varepsilon) + \varepsilon^3 \overline f_{h, 3}(\overline h, \overline w, \overline \varphi, \varepsilon), \\
  \dot{\overline w} &= \varepsilon \overline f_{w, 1}(\overline h, \overline w, \varepsilon) + \varepsilon^2 \overline f_{w, 2}(\overline h, \overline w, \varepsilon) + \varepsilon^3 \overline f_{w, 3}(\overline h, \overline w, \overline \varphi, \varepsilon), \\
  \dot{\overline \varphi} &= \omega(\overline h, \overline w) + \varepsilon \overline f_{\varphi, 1}(\overline h, \overline w, \varepsilon) + \varepsilon^2 \overline f_{\varphi, 2}(\overline h, \overline w, \overline \varphi, \varepsilon).
\end{split}
\end{align}
Let us call the new chart $\overline h, \overline w, \overline \varphi$ the \emph{averaging chart}.
For brevity we will often omit the dependence of the functions $f_*$, $\overline f_{*, *}$ and $u_{*, *}$ on $\varepsilon$.

It is convenient to denote by $v$ the column vector $(h, w)$ and by
$\overline v$ the column vector $(\overline h, \overline w)$. Let
$\overline f_{v, i} = (\overline f_{h, i}, \overline f_{w, i}), \;
u_{v, i} = (u_{h, i}, u_{w, i})$.
\begin{lemma} \label{l:varchange-short}
  For $k=h,w, \; i=1,2$ and for $k=\varphi, \; i=1$ we have
  \begin{equation} \label{e:overline-f}
      \overline f_{k, i}(h, w) = \langle Y_{k, i}(h, w, \varphi) \rangle_\varphi,
  \end{equation}
  \begin{equation} \label{e:homological}
      \overline f_{k, i}(h, w) + \omega(h, w) \pdv{u_{k, i}}{\varphi}\/(h, w, \varphi) = Y_{k, i}(h, w, \varphi)
  \end{equation}
  with
  \begin{align}
  \begin{split} \label{e:homological-rhs}
    Y_{h, 1} &= f_{h}, \\
    Y_{w, 1} &= f_{w}, \\
    Y_{\varphi, 1} &= f_{\varphi} + \pdv{\omega}{v} u_{v, 1}, \\
    Y_{h, 2} &=
    \pdv{f_h}{v} u_{v, 1}
    + \pdv{f_h}{\varphi} u_{\varphi, 1}
    - \pdv{u_{h, 1}}{v} \overline f_{v, 1}
    - \pdv{u_{h, 1}}{\varphi} \overline f_{\varphi, 1}, \\
    Y_{w, 2} &=
    \pdv{f_w}{v} u_{v, 1}
    + \pdv{f_w}{\varphi} u_{\varphi, 1}
    - \pdv{u_{w, 1}}{v} \overline f_{v, 1}
    - \pdv{u_{w, 1}}{\varphi} \overline f_{\varphi, 1}.
  \end{split}
  \end{align}
  The formulas for $\overline f_{h, 3}$, $\overline f_{w, 3}$ and $\overline f_{\varphi, 2}$ are stated in Lemma~\ref{l:varchange} below.
\end{lemma}
\noindent We will prove this lemma in Section~\ref{s:varchange-formulas}. The formulas above uniquely define $\overline f_{k, i}$ and $u_{k, i}$ under an additional assumption that
for $k=h,w; \; i=1,2$ and for $k=\varphi; \; i=1$ we have (in the formula below
$\langle * \rangle_\varphi$ denotes averaging with respect to $\varphi$)
\begin{equation}
  \langle u_{k, i} \rangle_\varphi = 0. \label{e:u-avg}
\end{equation}
We will always assume this to hold.

For $h \to 0$ many expressions introduced above tend to infinity. The estimates for these expressions are gathered in Table~\ref{t:est} below.

\begin{lemma} \label{l:invertible}
  There exists a constant $C_{inv} > 0$ depending on the perturbed system~\eqref{e:perturbed} such that for $h > C_{inv} \varepsilon$ the coordinate change given by~\eqref{e:varchange} is invertible.
\end{lemma}
\noindent This lemma is proved in Section~\ref{s:proof_prelim} below.

\noindent Using that
$\langle \pdv{u_{k, i}}{\varphi} \rangle_\varphi = 0,
\; \langle \pdv{u_{k, i}}{h} \rangle_\varphi
= \pdv{}{h} \langle u_{k, i} \rangle_\varphi = 0$,
we can simplify~\eqref{e:overline-f} for $\overline f_{h, 2}$ and $\overline f_{w, 2}$:
\begin{align} \label{e:f-h-2-init}
\begin{split}
  \overline f_{h, 2} &= \langle
  \pdv{f_h}{h} u_{h, 1} + \pdv{f_h}{w} u_{w, 1} + \pdv{f_h}{\varphi} u_{\varphi, 1} \rangle_\varphi, \\
  \overline f_{w, 2} &= \langle
  \pdv{f_w}{h} u_{h, 1} + \pdv{f_w}{w} u_{w, 1} + \pdv{f_w}{\varphi} u_{\varphi, 1} \rangle_\varphi.
\end{split}
\end{align}
{\color{myblue}
As $\langle u_{v, 1} \rangle_\varphi = 0$, we can also simplify the formula for $\overline f_{\varphi, 1}$:
\begin{equation} \label{e:avg-f-phi}
  \overline f_{\varphi, 1} = \langle f_\varphi \rangle_\varphi.
\end{equation}
}

\noindent The following formula is similar to Formula $2$ from~\cite{neishtadt86}.
\begin{lemma} \label{l:u-int-t}
  \begin{equation} \label{e:u-int-t}
    u_{a, 1}(h, w, t_0) = \frac{1}{T} \int_0^T \Big(t - \frac T 2\Big) f_a(h, w, t + t_0) dt
    \qquad \text{for } a=h, w_1, \dots, w_k.
  \end{equation}
  Here the third argument in $u_{a, 1}$ and $f_a$ is not $\varphi$, as usual, but the time $t = \varphi T / (2\pi)$. We use the notation $f_a(h, w, t) = f_a(h, w, \varphi(h, w, t))$ and a similar notation for $u_{a, 1}$.
\end{lemma}
\noindent This lemma is proved in Section~\ref{s:varchange-formulas} below.
Note that this formula for $u$ can also be rewritten as follows:
\begin{equation} \label{e:u-int-phi}
  u_{a, 1}(h, w, t_0) = \frac{1}{2 \pi} \int_0^T (\varphi(t) - \pi) f_a(h, w, t+t_0) dt
  \qquad \text{for } a=h, w_1, \dots, w_k.
\end{equation}

\subsection{Averaged system of order~$2$} \label{s:avg-eq}
The coefficients of the perturbed system~\eqref{e:init-avg} in the averaging chart depend on $\varepsilon$.
We would like the coefficients of the averaged system that we define in this section to be independent of $\varepsilon$. To this end, let us introduce some notation.
First, let us expand
\begin{equation} \label{e:f-series}
  f(p, q, z, \varepsilon) = f^0(p, q, z) + \varepsilon f^1(p, q, z) + \varepsilon^2 f^2(p, q, z, \varepsilon),
\end{equation}
where $f^0(p, q, z) = f(p, q, z, 0)$ and $f^1(p, q, z) = \pdv{f}{\varepsilon}\/(p, q, z, 0)$. Clearly, $f^0_q$, $f^0_p$, $f^0_z$, $f^1_q$, $f^1_p$ and $f^1_z$ are smooth functions of $p$, $q$ and $z$. The functions $f^2_p$, $f^2_q$ and $f^2_z$ are smooth functions of $p$, $q$ and $z$ that depend on $\varepsilon$ and are uniformly bounded by some constant independent of $\varepsilon$ (by Taylor's theorem with the Lagrange form of remainder).
Let us also consider the perturbed system~\eqref{e:init} with the perturbation $\varepsilon f^0(h, w, \varphi)$ instead of $\varepsilon f(h, w, \varphi, \varepsilon)$. For such system we may also consider a coordinate change of form~\eqref{e:varchange} that transforms it to the form~\eqref{e:init-avg}. Let us add an upper index $0$ to the coefficients of these equations (e.g. $u^0_{h, 1}$, $\overline f^0_{\varphi, 1}$) to show that we started with the perturbation $\varepsilon f^0$. The coefficients $u^0_{*, *}$ and $\overline f^0_{*, *}$ are determined by the same formulas as $u_{*, *}$ and $\overline f_{*, *}$, but we should plug $f^0$ instead of $f$ into those formulas.

Now let us rewrite~\eqref{e:init-avg} in such way that only the coefficients next to the largest powers of $\varepsilon$ depend on $\varepsilon$. This is done simply by expanding the coefficients similarly to~\eqref{e:f-series}. The resulting system will be
\begin{align}
\begin{split} \label{e:init-avg-trimmed}
  \dot{\overline h} &= \varepsilon \hat {f}_{h, 1}(\overline h, \overline w) + \varepsilon^2 \hat f_{h, 2}(\overline h, \overline w) + \varepsilon^3 \hat f_{h, 3}(\overline h, \overline w, \overline \varphi, \varepsilon), \\
  \dot{\overline w} &= \varepsilon \hat {f}_{w, 1}(\overline h, \overline w) + \varepsilon^2 \hat f_{w, 2}(\overline h, \overline w) + \varepsilon^3 \hat f_{w, 3}(\overline h, \overline w, \overline \varphi, \varepsilon), \\
  \dot{\overline \varphi} &= \omega(\overline h, \overline w) + \varepsilon \hat f_{\varphi, 1}(\overline h) + \varepsilon^2 \hat f_{\varphi, 2}(\overline h, \overline w, \overline \varphi, \varepsilon),
\end{split}
\end{align}
where
\begin{equation} \label{e:def-f-hat}
  \hat f_{*, 1} = \overline f^0_{*, 1} \text{ for } * = h, w, \varphi, \qquad
  \hat f_{h, 2} = \overline f^0_{h, 2} + \langle f_h^1(h, w, \varphi) \rangle_\varphi, \qquad
  \hat f_{w, 2} = \overline f^0_{w, 2} + \langle f_w^1(h, w, \varphi) \rangle_\varphi
\end{equation}
(here $f^1_h$ and $f^1_w$ are the $h$- and $w$-components of $f^1$ written in $(h, w, \varphi)$ coordinates), and $\hat f_{\varphi, 2}$, $\hat f_{h, 3}$ and $\hat f_{w, 3}$ satisfy the estimates in Table~\ref{t:est}. The estimates for $\hat f_{\star, \star}$ will be proved in Lemma~\ref{l:est-f-hat} below, one can also find formulas for $\hat f_{\varphi, 2}, \hat f_{h, 3}, \hat f_{w, 3}$ there. Also note that by~\cite[Corollary 3.1]{neishtadt17} we have $\int_0^T f^0_h dt = -\Theta_i(w) + O(h \ln h)$, so we have
\begin{equation} \label{e:hat-f-h-1}
  \hat f_{h, 1} = \frac{ -\Theta_i(w) + O(h \ln h) }{T}.
\end{equation}

The \emph{averaged system of order~$2$} is obtained from the system~\eqref{e:init-avg-trimmed} by
removing all terms on the right hand side that depend on $\overline \varphi$:
\begin{align}
\begin{split} \label{e:avg}
  \dot{\hat h} &= \varepsilon \hat{f}_{h, 1}(\hat h, \hat w) + \varepsilon^2 \hat{f}_{h, 2}(\hat h, \hat w), \\
  \dot{\hat w} &= \varepsilon \hat{f}_{w, 1}(\hat h, \hat w) + \varepsilon^2 \hat{w}_{h, 2}(\hat h, \hat w), \\
  \dot{\hat \varphi} &= \omega(\hat h, \hat w) + \varepsilon \omega_1(\hat h, \hat w).
\end{split}
\end{align}
Here we denote $\omega_1(\hat{h}) = \hat f_{\varphi, 1}(\hat h, \hat w)$ in order to match with~\cite{neishtadtVasiliev2005}. We will sometimes call this system simply the \emph{averaged system}.

Finally, let us use the slow time $\tau = \varepsilon t$ instead of $t$, then~\eqref{e:avg} rewrites as~\eqref{e:avg-form}.

\section{Proofs: overview} \label{s:proofs-overview}
{ \color{myblue}
We use the following extra assumptions in the proofs below.
\begin{itemize}
  \item All proofs are written for the domain $G_3$, averaging in the domains $G_1$ and $G_3$ can be treated in the same way. We will also assume that
  the separatrices in the phase portrait of the unperturbed system form a figure eight
  (Figure~\ref{f:up}).
  \item We will assume $H = 0$ on the separatrices, this can be achieved by replacing $H$ with $H - H(C(z), z)$. We will additionaly assume $H > 0$ in $G_3$ and $H < 0$ in $G_1 \cup G_2$, if the signs are opposite, we can swap $p$ and $q$, replacing $H$ by $-H$. Then we have $H = h$.
  \item We will use a certain transversal $\varphi = 0$ tangent to the bisector of an angle between separatrices, as stated in Section~\ref{s:assumptions}, when we obtain estimates on functions encountered in the formulas for order 2 averaging. This transversal is determined by Moser's normal form as described in Section~\ref{s:Moser} below. We will show that results of this paper (proved for one special transversal) also hold for any other transversal tangent to the bisector in Remark~\ref{r:any-transversal} and Remark~\ref{r:phase-transversal}.

\end{itemize}

Let us now state a plan of the proofs. We state this plan in the logical order here; the order in the paper is different, because we moved more technical parts closer to the end.
\begin{enumerate}
  \item A certain relation between partial derivatives of the perturbation $f$ in energy-angle variables (namely, between $\pdv{f_\varphi}{\varphi}$ and $\pdv{f_h}{h}$) will be important throughout the whole paper.
  This relation comes from the fact that divergence of the perturbation $f$ written in the action-angle variables is the same as in the coordinates $p, q$, as the coordinate change $(p, q) \mapsto (I, \varphi)$ is volume-preserving, so
  we call this relation \emph{divergence lemma}. It is stated and proved in Section~\ref{s:div}.

  \item The angle variable $\varphi$ behaves badly near the separatrices, the $\varphi$-component of the perturbation $f_\varphi$ may be unbounded.
  Partial derivative $\pdv{}{h}$ also behaves badly, as it is taken for fixed $\varphi$.
  We obtain estimates describing the perturbation in energy-angle variables near separatrices:
  we estimate $f_\varphi$ and and its partial derivatives (with respect to $h, w, \varphi$) and also partial derivatives of first and second order of $f_h$ and $f_w$.

  Estimates for derivatives of $f_h$ and $f_w$ are obtained using Moser's normal form~\cite{moser1956analytic} in the following way: we use extra coordinate chart connected with this normal form and compute partial derivatives via chain rule through this extra chart (for second order derivatives analogue~\eqref{e:d2-composition} of chain rule is used).
  Estimates for $f_\varphi$ are obtained using estimates for $\pdv{f_h}{h}$ and divergence lemma that connects $\pdv{f_h}{h}$ and $\pdv{f_\varphi}{\varphi}$.
  This is done in Section~\ref{s:est-action-angle}.

  \item We write formulas describing averaging method of order 2, namely, for the functions $u_{*, *}$ and $\overline f_{*, *}$ introduced in Section~\ref{s:avg-eq}. This is done in Section~\ref{s:est-u}.

  \item We estimate the functions $u_{*, *}$ and $\overline f_{*, *}$. Most of these estimates are obtained in a straightforward way using the formulas for these functions and estimates for $f$ in energy-angle variables (Section~\ref{s:est-u}).

  Particular care is needed when estimating $\overline f_{h, 2}$ (and also $\overline f_{w, 2}$, but let us focus on $\overline f_{h, 2}$), as good estimate for $\overline f_{h, 2}$ is crucial to prove that solutions of averaged system of order~2 cross separatrices. Straightforward estimate is not enough, we use a careful argument based on rewriting formula for $\overline f_{h, 2}$ using integration by parts and then using divergence lemma to obtain the estimate~$\overline f_{h, 2} = O(\ln^{-1} h)$ (Section~\ref{s:est-f-2}).
  \item
  Using estimates above and the standard approach for justification of order 2 averaging method gives Theorem~\ref{t:avg-2}.
  This is done in Section~\ref{s:proof_prelim} (some preliminary lemmas) and Section~\ref{s:approx-proof}.
  \item
  Formula for the pseudo-phase is proved using the same scheme as in~\cite{neishtadtVasiliev2005}. In this scheme one uses averaging until very close to separatrices ($h \approx \varepsilon^{2/3}$ up to some power of $\ln \varepsilon$) and in the immediate vicinity of separatrices one uses estimate on the change of $h$ during one wind around the separatrices. Estimates on the averaging are provided by the current paper and estimates used when very close to separatrices are taken from~\cite{neishtadt17}.

  In addition to this scheme~\cite{neishtadtVasiliev2005}, we also need a certain cancellation lemma, it is proved in Section~\ref{s:cancel} using divergence lemma.
  \item
  Formula for probabilities of capture (using Anosov's definition) easily follows from the formula for pseudo-phase (Section~\ref{s:prob-proofs}).
\end{enumerate}
}

\begin{table}
\begin{center}
\begin{tabular}{ |c|c|c| }
 \hline
 Expression & Estimates & Obtained in \\
 \hline \hline
$T$ & $T, \pdv{T}{w}, \pdv[2]{T}{w} = O(\ln(h)); \; \pdv{T}{h}, \pdv[2]{T}{h}{w} = O(h^{-1}); \; \pdv[2]{T}{h} = O(h^{-2})$
& Section~\ref{s:est-T} \\
 \hline
$\omega$ & $\omega, \pdv{\omega}{w}, \pdv[2]{\omega}{w} = O(\ln^{-1} h); \;
\pdv{\omega}{h}, \pdv[2]{\omega}{h}{w} = O(h^{-1}\ln^{-2}h); \;
\pdv[2]{\omega}{h} = O(h^{-2}\ln^{-2}h)$
& Section~\ref{s:est-T} \\
\hline
$f_{w_i}$
& \makecell{
  $f_{w_i}, \pdv{f_{w_i}}{w}, \pdv[2]{f_{w_i}}{w} = O(1); \;
  \pdv{f_{w_i}}{h}, \pdv[2]{f_{w_i}}{h}{w} = O_*(h^{-1}\ln^{-1} h);$ \\
  $\pdv{f_h}{\varphi}, \pdv[2]{f_{w_i}}{w}{\varphi} = O_*(\ln h);$ \;
  $\pdv[2]{f_{w_i}}{h} = O_*(h^{-2} \ln^{-1} h);$ \\
  $\pdv[2]{f_{w_i}}{h}{\varphi} = O_*(h^{-1}); \;
  \pdv[2]{f_{w_i}}{\varphi} = O_*(\ln^2 h)$
}
& Section~\ref{s:est-f} \\
\hline
$f_h$
& \makecell{
  $f_h, \pdv{f_h}{w}, \pdv[2]{f_h}{w} = O_*(1)$, other estimates as for $f_{w_i}$
}
& Section~\ref{s:est-f} \\
\hline
$\Div f$
& \makecell{
  As for $f_{w_i}$
}
& Section~\ref{s:est-f} \\
\hline
$f_\varphi$ &
\makecell{
  $f_\varphi, \pdv{f_\varphi}{w} = O_*(h^{-1}\ln^{-2} h); \; f_\varphi(h, w, 0) = O(h^{-1/2} \ln^{-1} h);$ \\
  $\pdv{f_\varphi}{\varphi} = O_*(h^{-1}\ln^{-1} h); \; \pdv{f_\varphi}{h} = O_*(h^{-2} \ln^{-2} h);$
}
& Section~\ref{s:est-f} \\
\hline
$u_{h, 1}$ &
\makecell{
  $u_{h, 1}, \pdv{u_{h, 1}}{w}, \pdv[2]{u_{h, 1}}{w} = O(1); \;
  \pdv{u_{h, 1}}{\varphi}, \pdv[2]{u_{h, 1}}{\varphi}{w} = O(\ln h);$ \\
  $\pdv{u_{h, 1}}{h}, \pdv[2]{u_{h, 1}}{h}{w} = O(h^{-1} \ln^{-1} h); \;
  \pdv[2]{u_{h, 1}}{h}{\varphi} = O(h^{-1}); \;
  \pdv[2]{u_{h, 1}}{h} = O(h^{-2} \ln^{-1} h)$
}
& Section~\ref{s:est-u} \\
\hline
$u_{w_i, 1}$ & As for $u_{h, 1}$ & Section~\ref{s:est-u} \\
\hline
$u_{\varphi, 1}$ &
$u_{\varphi, 1}, \pdv{u_{\varphi, 1}}{\varphi}, \pdv{u_{\varphi, 1}}{w} = O(h^{-1} \ln^{-1} h); \;
\pdv{u_{\varphi, 1}}{h} = O(h^{-2} \ln^{-1} h);$
& Section~\ref{s:est-u} \\
\hline
$\overline f_{w_i, 1}$ & $\overline f_{w_i, 1}, \pdv{\overline f_{w_i, 1}}{w} = O(1), \; \pdv{\overline f_{w_i, 1}}{h} = O(h^{-1} \ln^{-2} h)$
& Section~\ref{s:est-u} \\
\hline
$\overline f_{h, 1}$ & $\overline f_{h, 1}, \pdv{\overline f_{h, 1}}{w} = O(\ln^{-1} h)$; \;
$\pdv{\overline f_{h, 1}}{h} = O(h^{-1} \ln^{-2} h)$
& Section~\ref{s:est-u} \\
\hline
$\overline f_{\varphi, 1}$ & $\overline f_{\varphi, 1}, \pdv{\overline f_{\varphi, 1}}{w} = O(h^{-1} \ln^{-3} h); \; \pdv{\overline f_{\varphi, 1}}{h} = O(h^{-2}\ln^{-3} h)$
& Section~\ref{s:est-u} \\
\hline
$u_{h, 2}$ &
$u_{h, 2}, \pdv{u_{h, 2}}{w} = O(h^{-1}); \;
\pdv{u_{h, 2}}{\varphi} = O(h^{-1} \ln h); \;
\pdv{u_{h, 2}}{h} = O(h^{-2})$
& Section~\ref{s:est-u} \\
\hline
$u_{w_i, 2}$ & As for $u_{h, 2}$ & Section~\ref{s:est-u} \\
\hline
$\overline f_{h, 2}$ & $\overline f_{h, 2} = O(\ln^{-1} h), \; \pdv{\overline f_{h, 2}}{h} = O(h^{-2} \ln^{-1} h), \; \pdv{\overline f_{h, 2}}{w} = O(h^{-1} \ln^{-1} h)$
& Section~\ref{s:est-u} \\
\hline
$\overline f_{w_i, 2}$ & $\overline f_{w_i, 2} = O(h^{-1} \ln^{-3} h), \; \pdv{\overline f_{w_i, 2}}{h} = O(h^{-2} \ln^{-1} h), \; \pdv{\overline f_{w_i, 2}}{w} = O(h^{-1} \ln^{-1} h)$
& Section~\ref{s:est-u} \\
\hline
$\overline f_{\varphi, 2}$ & \makecell{
  $\overline f_{\varphi, 2} = O(h^{-2} \ln^{-2} h) + O_*(h^{-2} \ln^{-1} h)$ for $h > C_h \varepsilon$; \\
  $\overline f_{\varphi, 2} = O(h^{-2} \ln^{-2} h)$ for $h > C_h \varepsilon \abs{\ln \varepsilon}^{0.5}$.
}
& Section~\ref{s:est-u} \\
\hline
$\overline f_{h, 3}$ & $\overline f_{h, 3} = O(h^{-2} \ln^{-1} h) + O_*(h^{-2})$ for $h > C_h \varepsilon$.
& Section~\ref{s:est-u} \\
\hline
$\overline f_{w, 3}$ & $\overline f_{w, 3} = O(h^{-2} \ln^{-1} h) + O_*(h^{-2})$ for $h > C_h \varepsilon$.
& Section~\ref{s:est-u} \\
\hline
$\hat f_{*, *}$ & The estimates for $\hat f_{a, i}$ and its derivatives are as for $\overline f_{a, i}$.
& Section~\ref{s:est-u} \\
\hline
\end{tabular}
\end{center}
\caption{Estimates used in this paper. The notation $g = O_*(h^\alpha \ln^\beta h)$ means that $g = O(h^\alpha \ln^\beta h)$ and $g \ll h^\alpha \ln^\beta h$ near the saddle $C$, see details in Remark~\ref{r:O-star} below.
The constant $C_h > 0$ depends only on $H$ and $f$, it will be defined in Section~\ref{s:est-u}.
Let us also note that some of the expressions above are vectors, for them their norm is estimated, i.e. $\pdv{T}{w} = O(\ln^{-1} h)$ means $\norm{\pdv{T}{w}} = O(\ln^{-1} h)$. \label{t:est}}
\end{table}

\newpage
\section{Divergence lemma} \label{s:div}

\begin{lemma} \label{l:trace}
  \begin{equation} \label{e:trace}
    \pdv{f_h}{h} + \pdv{f_\varphi}{\varphi} + \sum_{w_i} \pdv{f_{w_i}}{w_i} + \frac{1}{T} \Big( \pdv{T}{h} f_h + \sum_{w_i} \pdv{T}{w_i} f_{w_i} \Big) = \Div(f),
  \end{equation}
  where $ \Div(f) = \pdv{f_q}{q} + \pdv{f_p}{p} + \sum_{z_i} \pdv{f_{z_i}}{z_i}.$
\end{lemma}
\begin{proof}
The lemma immediately follows from the Voss--Weyl formula for the divergence (cf., e.g.,~\cite[\S 9.8]{grinfeld2013introduction}). Let us now state this formula in a slightly modified form, with the product rule applied to one of the terms.
Suppose $\tilde x_i$ are curvilinear coordinates and $x_i$ are cartesian coordinates. Let $\tilde f_i$ and $f_i$ be components of a vector field $f$ in these coordinates. Let $D$ be the Jacobian of the coordinate change $T$ given by $x = T(\tilde x)$. Then
\begin{equation}
    \Div f = \sum \pdv{f_i}{x_i} = \sum \pdv{\tilde f_i}{\tilde x_i}
    + D^{-1} \sum \pdv{D}{\tilde x_i} \tilde f_i .
\end{equation}

Let us apply this formula the to coordinate systems $\tilde x = (h, \varphi, w)$ and $x = (p, q, z)$. For any fixed $z$ the map $(h, \varphi) \mapsto (p, q)$ has the Jacobian $(\frac{\partial h}{\partial I})^{-1} = (\omega(h, z))^{-1}$. Thus, $D = (\omega(h, z))^{-1} = \frac{T(h, z)}{2\pi}$. The $2\pi$ factor cancels out, and we get~\eqref{e:trace}.
\end{proof}

\section{Lemmas on order 2 averaging} \label{s:proof_prelim}
\begin{proof}[Proof of Lemma~\ref{l:avg-cross}]
  { \color{myblue}
  We assumed in Section~\ref{s:assumptions} that solutions of order 1 averaged system cross separatrices. This means that solutions of order 2 averaged system come close to separatrices for small $\varepsilon$. In the rest of the proof we consider solution of order $2$ averaged system starting near separatices.
  }

  By~\eqref{e:hat-f-h-1} and the estimate on $\hat f_{h, 2}$ from Table~\ref{t:est} we get that
  \begin{equation} \label{e:d-h-tau}
    \dv{\hat h}{\tau} = \frac{-\Theta_i(w) + O(\hat h \ln \hat h) + O(\varepsilon)}{T}
  \end{equation}
  along solutions of averaged system~\eqref{e:avg-form} of order~$2$.
  As $\Theta_i > 0$, this means that any solution $\hat h(\tau), \hat w(\tau), \hat \varphi(\tau)$ of the averaged system of order~$2$ starting close to the separatrices crosses the separatrix of the initial unperturbed Hamiltonian equation.
\end{proof}
Denote by $\tau_*$ the slow time at the moment of crossing, $\hat h(\tau_*) = 0$.
From~\eqref{e:d-h-tau} we also see that for small $\varepsilon$, $h$ and $\tau < \tau_*$ the function $\hat h(\tau)$ is decreasing.
By~\eqref{e:d-h-tau} we also have that along solutions of the averaged system of order~$2$
\begin{equation} \label{e:d-tau-h}
  \dv{\tau}{\hat h} = - \frac{T}{\Theta_i(\hat w)} (1 + O(\hat h \ln \hat h) + O(\varepsilon)).
\end{equation}

\begin{proof}[Proof of Lemma~\ref{l:invertible}]
  Set $v = (h, w)$ and $\overline v = (\overline h, \overline w)$.
  Let us denote by $F$ the map
  $(\overline v, \overline \varphi) \to (v, \varphi)$
  given by~\eqref{e:varchange}.
  Let us consider the domain $\overline h > C \varepsilon$, where the constant $C > 0$ is large. From Table~\ref{t:est} (note that as the values of $u_{*, *}$ are taken at $(\overline v, \overline \varphi)$, so we should plug $h = \overline h$ in the estimates in Table~\ref{t:est}) and $\varepsilon \overline h^{-1} < C^{-1}$ we have $\abs{h - \overline h} = O(\varepsilon)$, so for large enough $C$ we have
  $0.5 h < \overline h < 2h$.
  This means that $h > C_{inv} \varepsilon$ implies $\overline h > C \varepsilon$ for $C_{inv} = 2 C$. This also means that we can write $O(h^{-1})$ instead of $O(\overline h^{-1})$, and so on.

  We can estimate the coefficients of the Jacobian matrix of $F$
  using Table~\ref{t:est} and $\varepsilon \overline h^{-1} < C^{-1}$. For $C \to \infty$ all coefficients tend to the corresponding coefficients of the identity matrix
  except $\pdv{\varphi}{h} = O(\varepsilon h^{-2} \ln^{-1} h)$. However, as all elements of the last column of $DF$ except the diagonal one (i.e. $\pdv{*}{\overline \varphi}, \; * = h, w_1, \dots, w_k$) are $O(\varepsilon \ln^{-1} h)$ for $\varepsilon \lesssim h$, any summand in $\det DF$ containing $\pdv{\varphi}{\overline h}$ is $O(\varepsilon^2 h^{-2})$. So
  $\det DF \to 1$ for $C \to \infty$, hence for some $C$ this determinant lies in $[0.5, 2]$ for $\overline h > C \varepsilon$. By the inverse function theorem this implies that $F$ is a local diffeomorphism.
  Moreover, for $\overline h > C \varepsilon$ we have $\norm{F(x) - x} = O(\ln^{-1} \varepsilon)$. Indeed, $\varepsilon u_{\varphi, 1} = O(\varepsilon h^{-1}\ln^{-1} h) = O(C^{-1} \ln^{-1} \varepsilon)$; we can estimate $\varepsilon u_{v, 1}$ and $\varepsilon^2 u_{v, 2}$ in the same way. Therefore, $F$ is invertible as a local diffeomorphism that is $C^0$-close to the identity.
\end{proof}

{ \color{myblue}
\begin{remark} \label{r:w-star}
  Let $w_{*,2}$ denote the value of $w$ when a solution of order $2$ averaged system with some initial data $v_0 = (h_0, w_0)$ crosses separatrices. Take initial data $v_0'$ with $\norm{v_0 - v_0'} = O(\varepsilon)$. Let $w_{*, 1}$ denote the value of $w$ when the solution of order $1$ averaged system with initial data $v_0'$ crosses separatrices. Then
  \begin{equation}
    \norm{w_{*,2} - w_{*,1}} = O(\varepsilon).
  \end{equation}
\end{remark}
\begin{proof}
  It is easy to check that for order 1 averaged system $O(\varepsilon)$ change in initial data leads to $O(\varepsilon)$ change in the value of $w$ at separatrix crossing. Thus it is enough to prove the lemma for $v_0 = v_0'$, and we will assume that this holds.

  Denote by $(\hat h_1, \hat w_1)$ and $(\hat h_2, \hat w_2)$ solutions of order $1$ and $2$, respectively, averaged systems with the same initial data. Both $\hat h_1(\tau)$ and $\hat h_2(\tau)$ decrease, so we may use $h$ as an independent variable instead of the slow time $\tau$. We have
  \begin{equation}
    \dv{\hat w_1}{\hat h_1} = \frac{\hat f_{w, 1}}{\hat f_{h, 1}},
    \qquad
    \dv{\hat w_2}{\hat h_2} = \frac{\hat f_{w, 1} + \varepsilon \hat f_{w, 2}}{\hat f_{h, 1} + \varepsilon \hat f_{h, 2}}.
  \end{equation}
  Denote $g(h, w) = \frac{\hat f_{w, 1}}{\hat f_{h, 1}}$. We have $\hat f_{h, 1} \sim \ln^{-1} h$. Plugging in the second formula the estimates from Table~\ref{t:est}, we get
  \begin{equation}
    \dv{\hat w_2}{\hat h_2} = g + O(\varepsilon h^{-1} \ln^{-2} h).
  \end{equation}
  We also have (by Table~\ref{t:est} and the quotient rule) $\norm{\pdv{g}{w}} = O(\ln h)$.

  We have two solutions $\hat w_1(h)$ and $\hat w_2(h)$ with the same initial data $\hat w_i(h_0) = w_0$. Denote $\Delta w(h) = \hat w_2(h) - \hat w_1(h)$.
  For $h \in [h_*, h_0)$ we have the estimate
  \begin{equation}
    \norm{\dv{\Delta w}{h}} \le O(\ln h) \norm{\Delta w} + O(\varepsilon h^{-1} \ln^{-2} h).
  \end{equation}
  Set $u(h) = \varepsilon + \norm{\Delta w}$, then $u(h_0) = \varepsilon$ and $u' \le O(h^{-1} \ln^{-2} h) u$ (as $\ln h \lesssim h^{-1} \ln^{-2} h$).
  By Gronwall's inequality we get $u(0) \le u(h_0) O(1) = O(\varepsilon)$, as claimed.
\end{proof}
}

{\color{myblue}
The construction of averaged system of order $2$ depends on the choice of angle variable $\varphi$ (this choice is determined by the transversal $\varphi = 0$). It turns out that the first two equations (that describe the evolution of slow variables) do not depend on the choice of $\varphi$ and the last equation (it describes the evolution of $\varphi$) depends on the choice of $\varphi$ in a predictable way.
\begin{lemma} \label{l:change-phi}
  \leavevmode
  \begin{enumerate}
    \item The functions $\hat f_{a, i}$, where $a=h, w$, $i=1, 2$, and $u_{h, 1}$, $u_{w, 1}$ do not depend on the choice of $\varphi$.
    \item Suppose we have another angle variable $\psi$ connected with original angle $\varphi$ by the relation $\psi = \varphi + \Delta \varphi(h, z)$. Denote by $\omega_{1, \psi}$ the value of $\omega_1$ for angle variable $\psi$.
    Let $\hat v(\tau)$, $\hat v=(\hat h, \hat w)$ be a solution of the first two equations in~\eqref{e:avg-form} (they are the same for $\varphi$ and $\psi$ by the first part of this lemma).
    Then
    \begin{equation} \label{e:change-transversal}
      \int_{\tau_1}^{\tau_2} \omega_{1, \psi}(\hat v(\tau)) d\tau
      - \int_{\tau_1}^{\tau_2} \omega_1(\hat v(\tau)) d\tau
      = \Delta \varphi(\hat v(\tau_2)) - \Delta \varphi(\hat v(\tau_1))
      - \varepsilon \int_{\tau_1}^{\tau_2} \pdv{\Delta \varphi}{v} \hat f_{v, 2}(\hat v(\tau)) d\tau.
    \end{equation}
  \end{enumerate}
\end{lemma}
\begin{proof}[Proof of Lemma~\ref{l:change-phi}] $\qquad$

  \noindent 1.
  Clearly, $f_h$ and $f_w$ do not depend on the choice of $\varphi$. This means that $\hat f_{h, 1}$ and $\hat f_{w, 1}$, averages of these functions over $\varphi$, also do not depend on the choice of $\varphi$. The function $u_{h, 1}$ is uniquely determined by two conditions:
  $\omega \pdv{u_{h, 1}}{\varphi} = f_h - \langle f_h \rangle_\varphi$
  and
  $\langle u_{h, 1} \rangle_\varphi = 0$.
  As both $\langle \cdot \rangle_\varphi$ and $\pdv{}{\varphi}$ do not depend on the choice of $\varphi$, neither does $u_{h, 1}$. One can check in the same way that $u_{w, 1}$ does not depend on the choice of $\varphi$.

  Let us now focus on the functions $\hat f_{h, 2}$ and $\hat f_{w, 2}$.
  In order to define order 2 averaged system, we first find coordinate change~\eqref{e:varchange} that transforms~\eqref{e:init} to the form~\eqref{e:init-avg}; here in~\eqref{e:init-avg} order $1$ terms (in $\varepsilon$) for $\dot \varphi$ and order $2$ terms for $\dot {\overline h}$ and $\dot {\overline w}$ do not depend on $\varphi$; this gives the functions $\overline f_{h, 2}$ and $\overline f_{w, 2}$. We also require that averages over $\varphi$ of the functions $u$ from~\eqref{e:varchange} are zero.
  Then coordinate change~\eqref{e:varchange} and system~\eqref{e:init-avg} are uniquely defined. Let us drop the condition $\langle u_{\varphi, 1} \rangle_\varphi = 0$. Then the coordinate change is no more unique, but $\overline f_{h, 2}$ and $\overline f_{w, 2}$ are still uniquely defined due to~\eqref{e:f-h-2-init}, as $\langle \pdv{f_h}{\varphi} \rangle_\varphi = \langle \pdv{f_w}{\varphi} \rangle_\varphi = 0$.
  Thus adding some $\varphi$-independent function $\Delta u_\varphi(h, w)$ to $u_{\varphi, 1}$ does not change $\overline f_{h, 2}$ and $\overline f_{w, 2}$. But adding $\Delta u_\varphi(h, w)$ means a change of the transversal $\varphi=0$. Thus $\overline f_{h, 2}$ and $\overline f_{w, 2}$ do not depend on the choice of such transversal.
  By~\eqref{e:def-f-hat} this means that $\hat{f}_{h, 2}$ and $\hat{f}_{w, 2}$ also does not depend on the choice of $\varphi$.

  \noindent 2. We use formula~\eqref{e:avg-f-phi}:
  $\omega_1 = \overline f^0_{\varphi, 1} = \langle f^0_{\varphi} \rangle_\varphi$.
  We have $f_\psi = f_\varphi + \pdv{\Delta \varphi}{v} f_v$, $v = (h, w)$, thus
  \begin{equation}
    \omega_{1, \psi} = \omega_1 + \pdv{\Delta \varphi}{v} \hat f_{v, 1} = \omega_1 + \pdv{\Delta \varphi}{v} \dv{\hat v}{\tau} - \varepsilon \pdv{\Delta \varphi}{v} \hat f_{v, 2}.
  \end{equation}
  This gives
  \begin{align}
  \begin{split}
    \int_{\tau_1}^{\tau_2} \omega_{1, \psi}(\hat v(\tau)) d\tau - \int_{\tau_1}^{\tau_2} \omega_1(\hat v(\tau)) d\tau
    &=
    \int_{\tau_1}^{\tau_2} \pdv{\Delta \varphi}{v} \dv{\hat v}{\tau} d\tau -
    \varepsilon \int_{\tau_1}^{\tau_2} \pdv{\Delta \varphi}{v} \hat f_{v, 2} d\tau
    = \\
    &= \Delta \varphi(\hat v(\tau_2)) - \Delta \varphi(\hat v(\tau_1)) - \varepsilon \int_{\tau_1}^{\tau_2} \pdv{\Delta \varphi}{v} \hat f_{v, 2}(\hat v(\tau)) d\tau,
  \end{split}
  \end{align}
  as claimed.
\end{proof}

}

\section{Approximation lemma and proof of Theorem~\ref{t:avg-2}} \label{s:approx-proof}

The lemma below estimates how the solutions of the averaged system of order 2 approximate the solutions of~\eqref{e:init-avg} while approaching the separatrices. It will be proved in Section~\ref{s:proof-approx}.
\begin{lemma}[Approximation lemma] \label{l:close-avg-init}
  There exists $C_2 > 0$ such that the following holds.
  Consider a solution $\overline v(t), \overline \varphi(t)$ of~\eqref{e:init-avg}, where $\overline v(t) = (\overline h(t), \overline w(t))$, with initial condition
  $\overline v(0), \overline \varphi(0)$.
  Consider also a solution
  $\hat v(\tau) = (\hat h(\tau), \hat w(\tau))$ of the first two equations of the averaged system~\eqref{e:avg-form} of order~$2$ with initial condition
  $\hat v(0)$ such that $\norm{\overline v(0) - \hat v(0)} \le C_1 \varepsilon^2$ for some $C_1 > 0$.
  Then for all small enough $\varepsilon$ for any $t$ such that
  \begin{equation} \label{e:bound-h-eps}
    \hat h(\varepsilon t) > C_2 \varepsilon
  \end{equation}
  we have the following estimates (in the error terms below we write $h$ for $\hat h(\varepsilon t)$, e.g. $O(h)$ instead of $O(\hat h(\varepsilon t))$):
  \begin{align} \label{e:close-avg-init}
  \begin{split}
    &\norm{\overline{v}(t) - \hat{v}(\varepsilon t)} = O(\varepsilon^2 h^{-1}), \\
    &\overline \varphi(t) - \overline \varphi(0)
      = \varepsilon^{-1} \int_{0}^{\varepsilon t} \bigg( \omega(\hat v(\tau')) + \varepsilon \omega_1(\hat v(\tau')) \bigg) d\tau'
      + O(\varepsilon h^{-1} \ln^{-1} h).
  \end{split}
  \end{align}
\end{lemma}

{ \color{myblue}
\noindent Let us derive Theorem~\ref{t:avg-2} from Lemma~\ref{l:close-avg-init} and estimates in Table~\ref{t:est}.
\begin{proof}[Proof of Theorem~\ref{t:avg-2}]
  We are given a solution $X(t)$ of perturbed system, rewriting it in energy-angle variables and applying coordinate change~\eqref{e:varchange} (it is invertible by Lemma~\ref{l:invertible} when $h \gtrsim \varepsilon$ and this holds for all considered solutions if $C$ in the statement of Theorem~\ref{t:avg-2} is large enough) yields a solution $\overline X(t) = (\overline h(t), \overline w(t), \overline \varphi(t))$ of~\eqref{e:init-avg}.
  Theorem~\ref{t:avg-2} prescribes us to take solution of averaged system of second order~\eqref{e:avg-form} with initial data
  \begin{equation}
    (\hat h(0), \hat w(0), \hat \varphi(0))
    =
    (h_0 - u_{h, 1}(h_0, w_0, \varphi_0), w_0 - u_{w, 1}(h_0, w_0, \varphi_0), \varphi_0).
  \end{equation}
  For the initial data we have $h \sim 1$. This means $\norm{X(0) - \overline X(0)} = O(\varepsilon)$ and
  \begin{equation}
    |\hat h(0) - \overline h(0)| \le
    \varepsilon|u_{h, 1}(X(0)) - u_{h, 1}(\overline X(0))|
    + \varepsilon^2 |u_{h, 2}(\overline X(0))| =O(\varepsilon^2)
  \end{equation}
  by~\eqref{e:varchange} and as $\norm{\pdv{u_{h, 1}}{X}} = O(1)$ far from separatrices. Similarly, $\norm{\hat w(0) - \overline w(0)} = O(\varepsilon^2)$ and $|\overline \varphi(0) - \hat \varphi(0)| = O(\varepsilon)$. Thus the condition on initial data in Lemma~\ref{l:close-avg-init} is satisfied (for large enough $C$). This lemma gives~\eqref{e:close-avg-init}.

  Note that the integral in last line of~\eqref{e:close-avg-init} is equal to
  $\hat \varphi(\varepsilon t) - \hat \varphi(0)$.
  Thus last line of~\eqref{e:close-avg-init} and $\overline \varphi(0) = \hat \varphi(0) + O(\varepsilon)$ imply
  \begin{equation}
    \overline \varphi(t) = \hat \varphi(\tau) + O(\varepsilon h^{-1} \ln^{-1} h).
  \end{equation}
  We have $\overline h(t) \sim \hat h(\tau)$ due to the first line of~\eqref{e:close-avg-init} (for large enough $C$).
  From~\eqref{e:varchange} we get
  \begin{equation}
    \varphi(t) = \overline \varphi(t) + \varepsilon u_{\varphi, 1}(\overline X)
    = \overline \varphi(t) + O(\varepsilon h^{-1} \ln^{-1} h)
    = \hat \varphi(t) + O(\varepsilon h^{-1} \ln^{-1} h).
  \end{equation}
  This is the first estimate of Theorem~\ref{t:avg-2}.

  We have by~\eqref{e:varchange}
  \begin{align}
  \begin{split}
    h(t) &= \overline h(t) + \varepsilon u_{h, 1}(\overline X) + \varepsilon^2 u_{h, 2}(\overline X)
    = \overline h(t) + \varepsilon u_{h, 1}(\hat X) + \varepsilon \pdv{u_{h, 1}}{X}\/(\overline X - \hat X) + O(\varepsilon^2 h^{-1}) \\
    &= \overline h(t) + \varepsilon u^0_{h, 1}(\hat X) + O(\varepsilon^2 h^{-1})
  \end{split}
  \end{align}
  We have used the estimate
  \begin{align}
  \begin{split}
    &\pdv{u_{h, 1}}{X}\/(\overline X(t) - \hat X(\tau)) =
    \pdv{u_{h, 1}}{h}\/(\overline h(t) - \hat h(\tau))
    + \pdv{u_{h, 1}}{w}\/(\overline w(t) - \hat w(\tau))
    + \pdv{u_{h, 1}}{\varphi}\/(\overline \varphi(t) - \hat \varphi(\tau)) \\
    &\qquad= O(h^{-1} \ln^{-1} h) O(\varepsilon^2 h^{-1})
    + O(1) O(\varepsilon^2 h^{-1})
    + O(\ln h) O(\varepsilon h^{-1} \ln^{-1} h) = O(\varepsilon h^{-1}).
  \end{split}
  \end{align}
  obtained using Table~\ref{t:est} and the estimate $|u_{h, 1}(X, \varepsilon) - u^0_{h, 1}(X)| = O(\varepsilon)$ proved below in Lemma~\ref{l:est-f-hat}.
  Thus $h(t) = \hat h(\tau) + \varepsilon u_{h, 1}(\hat X) + O(\varepsilon^2 h^{-1})$.
  Similarly, $w(t) = \hat w(\tau) + \varepsilon u_{w, 1}(\hat X) + O(\varepsilon^2 h^{-1})$. This completes the proof of Theorem~\ref{t:avg-2}.
\end{proof}

\begin{proof}[Proof of Remark~\ref{r:any-transversal}]
  So far Theorem~\ref{t:avg-2} is proved for a special transversal $\varphi=0$ described in Section~\ref{s:Moser}. Suppose we have other transversal, denote by $\psi$ the angle variable such that $\psi=0$ on this other transversal. Denote by $\Delta \varphi(h, w)$ the phase shift between $\varphi$ and $\psi$: $\psi = \varphi + \Delta \varphi(h, w)$.
  Consider two cases.
  \begin{enumerate}
    \item Transversal $\psi=0$ passes through saddle $C(z)$ for all $z$.
    By Lemma~\ref{l:change-phi} the evolution of slow variables $\hat h(\tau)$, $\hat w(\tau)$ is the same for averaged systems of order 2 written using angle variable $\varphi$ and angle variable $\psi$. The function $\omega_1$ is different for $\varphi$ and $\psi$. We will use integral form of the last equation in order 2 averaged system~\eqref{e:avg-form}:
    \begin{equation}
      \hat \varphi(\tau) =  \varphi(\tau_0) + \int_{\tau_0}^\tau \varepsilon^{-1} \omega(\hat v) + \omega_1(\hat v) d \tilde \tau, \qquad
      \hat v = (\hat h(\tilde \tau), \hat w(\tilde \tau)).
    \end{equation}
    The integral of $\omega_0$ is the same for $\varphi$ and $\psi$, by~\eqref{e:change-transversal} the difference between the integrals of $\omega_1$ is
    $\Delta \varphi(\hat v(\tau)) - \Delta \varphi(\hat v(\tau_0)) + E$, where
    \begin{equation}
      E = -\varepsilon \int_{\tau_0}^{\tau} \pdv{\Delta \varphi}{v} \hat f_{v, 2}(\hat v(\tau)) d\tau
    \end{equation}
    denotes the error term in~\eqref{e:change-transversal}.
    Thus
    \begin{equation} \label{e:E}
      \hat \psi(\tau) = \hat \varphi(\tau) + \Delta \varphi(\hat v(\tau)) + E.
    \end{equation}

    We can prove that $E = O(\varepsilon)$ using the estimates on $\pdv{\Delta \varphi}{v}$ from Lemma~\ref{l:d-delta-phi}. Let us use $h = \hat h(\tau)$ as independent variable instead of $\tau$. By~\eqref{e:d-h-tau} we have $d\tau \sim \ln h \; dh$.
    \begin{equation}
      \int_{\tau_1}^{\tau_2} \pdv{\Delta \varphi}{v} \hat f_{v, 2}(\hat v(\tau)) d\tau
      = \int_{\hat h(\tau_1)}^{\hat h(\tau_2)} \Big(O(h^{-1} \ln^{-2} h) + O(h^{-1} \ln^{-3} h)\Big)dh = O(1).
    \end{equation}
    Thus indeed $E=O(\varepsilon)$. This and~\eqref{e:E} imply that the evolution of the angle for order 2 averaging systems written using $\varphi$ and $\psi$ is the same with $O(\varepsilon)$ error.

    \item Transversal $\psi=0$ passes far from the saddle $C(z)$ for all $z$. The proof is the same in this case, we simply use Lemma~\ref{l:d-delta-phi-2} instead of Lemma~\ref{l:d-delta-phi} to estimate $E$.
  \end{enumerate}
\end{proof}
}

\section{Proof of the approximation lemma} \label{s:proof-approx}
In this section we prove Lemma~\ref{l:close-avg-init}.
As far from the separatrices solutions of the averaged system approximate solutions of the perturbed system in the averaged chart with accuracy $O(\varepsilon^2)$ for time intervals $O(\varepsilon^{-1})$, we may assume that $\overline h(0) > 0$ is small enough.
Then $\hat h(\tau)$ will decrease monotonically.
It will be convenient to use the notation
$\overline v(\tau), \overline h(\tau), \overline w(\tau), \overline \varphi(\tau) = \overline v(t), \overline h(t), \overline w(t), \overline \varphi(t)$
with $t = \varepsilon^{-1} \tau$.

Let us start with the estimates for $\overline{h}(\tau) - \hat{h}(\tau)$ and $\overline{w}(\tau) - \hat{w}(\tau)$.
We will first only consider what happens up to some moment $\tau_{fin}$ such that for all $\tau < \tau_{fin}$ we have
\begin{equation} \label{e:2-times}
  0.5 \hat h(\tau) < \overline h(\tau) \le 2 \hat h(\tau), \qquad
  \hat h(\tau) \ge C_2 \varepsilon.
\end{equation}
In order to receive a better estimate, let us switch from $h$ to the action $I$.
Denote $\overline I = I(\overline v)$, $\hat I = I(\hat v)$, $\overline r = (\overline I, \overline w)$, $\hat r = (\hat I, \hat w)$.
Denote $\hat f_{I, i} = \pdv{I}{h} \hat f_{h, i} + \pdv{I}{w} \hat f_{w, i}$
and $\hat f_{r, i} = (\hat f_{I, i}, \hat f_{w, i})$.
We need the following lemma.
\begin{lemma} \label{l:d-f-0-I}
  We have
  \begin{align}
  \begin{split}
    &\pdv{h}{I} = \omega, \qquad
    \norm{\pdv{r}{v}} = O(\ln h), \qquad
    \norm{\pdv{h}{r}} = O(\ln^{-1} h), \\
    &\norm{\pdv{\hat f_{r, 1}}{r}} = O(h^{-1} \ln^{-3} h), \qquad
    \norm{\pdv{\hat f_{r, 2}}{r}} = O(h^{-2} \ln^{-1} h), \qquad
    \norm{\hat f_{r, 3}} = O_*(h^{-2} \ln h) + O(h^{-2}).
\end{split}
\end{align}
\end{lemma}
\begin{proof}[Proof of Lemma~\ref{l:d-f-0-I}]
From the Hamiltonian equations we have $\pdv{h}{I} = \omega$.
By~\cite[Corollary~3.2]{neishtadt17} we have
\[
  \pdv{I}{w} = O(1), \; \norm{\pdv[2]{I}{w}{h}} = O(\ln h), \; \norm{\pdv[2]{I}{w}} = O(1).
\]
As $\pdv{I}{h} = \omega^{-1}$, the first estimate implies $\norm{\pdv{r}{v}} = O(\ln h)$.
We have
$\pdv{I}{h} (\pdv{h}{w})_{I=const} + \pdv{I}{w} = 0$, this gives $(\pdv{h}{w})_{I=const} = O(\ln^{-1} h)$ and $\norm{\pdv{h}{r}} = O(\ln^{-1} h)$.

We have
\begin{equation} \label{e:f-I-i}
  \hat f_{I, i} = \pdv{I}{h} \hat f_{h, i} + \pdv{I}{w} \hat f_{w, i}.
\end{equation}
For $i=1$ this rewrites as
\begin{equation} \label{e:local-f-I-0}
  \hat f_{I, 1} = (2\pi)^{-1} \oint_{H=h} f^0_h dt + \pdv{I}{w} \hat f_{w, 1}.
\end{equation}
By~\cite[Lemma~3.2]{neishtadt17} we have
\begin{equation}
  \oint_{H=h} f^0_h dt = O(1), \qquad
  \pdv{}{h} \oint_{H=h} f^0_h dt = O(\ln h), \qquad
  \pdv{}{w} \oint_{H=h} f^0_h dt = O(1).
\end{equation}
Plugging the first estimate in~\eqref{e:local-f-I-0} gives $f_{I, 1} = O(1)$.
Plugging these estimates and the estimates in Table~\ref{t:est} in the derivatives of~\eqref{e:local-f-I-0} gives
\begin{equation}
  \pdv{\hat f_{I, 1}}{h} = O(h^{-1} \ln^{-2} h), \qquad
  \norm{\pdv{\hat f_{I, 1}}{w}} = O(1).
\end{equation}
As $\norm{\pdv{h}{r}} = O(\ln^{-1} h)$, we have
$\norm{\pdv{\hat f_{I, 1}}{r}} = \norm{\pdv{\hat f_{I, 1}}{h} \pdv{h}{r} + \pdv{\hat f_{I, 1}}{w}\pdv{w}{r}} = O(h^{-1} \ln^{-3} h)$.
We can prove that $\norm{\pdv{\hat f_{w, 1}}{r}} = O(h^{-1} \ln^{-3} h)$ in the same way, so
$\norm{\pdv{\hat f_{r, 1}}{r}} = O(h^{-1} \ln^{-3} h)$.

By~\eqref{e:f-I-i} for $a=h, w_i$ we have
\begin{equation}
  \pdv{\hat f_{I, 2}}{a} = \pdv[2]{I}{h}{a} \hat f_{h, 2} + \pdv{I}{h} \pdv{\hat f_{h, 2}}{a}
  + \pdv[2]{I}{w}{a} \hat f_{w, 2} + \pdv{I}{w} \pdv{\hat f_{w, 2}}{a},
\end{equation}
this yields $\pdv{\hat f_{I, 2}}{h} = O(h^{-2})$, $\norm{\pdv{\hat f_{I, 2}}{w}} = O(h^{-1})$ by Table~\ref{t:est}.
As $\norm{\pdv{h}{r}} = O(\ln^{-1} h)$, we have $\norm{\pdv{\hat f_{r, 2}}{r}} = O(h^{-2} \ln^{-1} h)$.

Finally, the estimate on $\hat f_{I, 3}$ follows from~\eqref{e:f-I-i} and Table~\ref{t:est}, while the estimate on $\hat f_{w, 3}$ is given by Table~\ref{t:est}.
\end{proof}

As $\overline v, \overline \varphi$ is a solution of~\eqref{e:init-avg}, it is also a solution of~\eqref{e:init-avg-trimmed}.
Rewriting~\eqref{e:init-avg-trimmed} and~\eqref{e:avg} using $r$ instead of $v$ gives
\begin{align} \label{e:I-dot}
\begin{split}
  \dot{\overline r} &= \varepsilon \hat f_{r, 1}(\overline r) + \varepsilon^2 \hat f_{r, 2}(\overline r) + \varepsilon^3 \hat f_{r, 3}(\overline r, \overline \varphi, \varepsilon), \\
  \dot{\hat{r}} &= \varepsilon \hat f_{r, 1}(\hat r) + \varepsilon^2 \hat f_{r, 2}(\hat r).
\end{split}
\end{align}
Denote $\Delta(\tau) = \norm{\overline r(\tau) - \hat r(\tau)}$.
From~\eqref{e:I-dot} we have the following differential inequality for $\Delta$:
\begin{equation} \label{e:a-b}
  \dv{\Delta}{\tau} \le a(\tau) \Delta + \varepsilon^2 b(\tau),
\end{equation}
where $a(\tau) = \norm{\big(\pdv{\hat f_{r, 1}}{r} \big)_{int} + \varepsilon \big(\pdv{\hat f_{r, 2}}{r}\big)_{int}}$ and $b(\tau) = \norm{\hat f_{r, 3}(\overline r(\tau), \overline \varphi(\tau))}$. Here the notation $\big(\pdv{\hat f_{r, i}}{r} \big)_{int}$ means that each row of this matrix is taken at some intermediate point in $[\overline r(\tau), \hat r(\tau)]$.
By~\eqref{e:2-times},~\eqref{e:bound-h-eps} and Lemma~\ref{l:d-f-0-I}
we have $a(\tau) = O(h^{-1} \ln^{-3} h) + \varepsilon O(h^{-2} \ln^{-1} h)$. By Lemma~\ref{l:d-f-0-I} we have $b(\tau) = O_*(h^{-2} \ln h) + O(h^{-2})$.

As in~\cite{neishtadtVasiliev2005}, we use the following estimate for $\Delta$ obtained by solving~\eqref{e:a-b}:
\begin{equation}
  \Delta(\tau) \le \exp \bigg( \int_0^\tau a(\tau')d\tau' \bigg)
  \bigg( \Delta(0) + \varepsilon^2 \int_0^\tau b(\tau') d \tau' \bigg).
\end{equation}
Using~\eqref{e:d-tau-h} and the estimates for $a$ and $b$, we can make a change of variable and compute the integrals above as integrals $d \hat h$. We have
\[
  \int_0^\tau a(\tau')d\tau' = O(1) + \varepsilon O(h^{-1}) = O(1), \text{ where } h = \hat h(\tau).
\]
The integral of $b$ can be estimated in the same way.
As $\int O_*(h^{-2} \ln h) dt$ during each wind of the trajectory of the perturbed system around the figure eight is $O(h^{-2} \ln h)$ and this wind takes time $O(\ln h)$, we can replace this function with its average $O(h^{-2})$ if we also add the integral over the last incomplete wind:
\[
  \int_0^\tau O_*(h^{-2} \ln h) d \tau'
  \sim
  \int_0^\tau O(h^{-2}) d \tau' + O(\varepsilon h^{-2} \ln h),
  \text{ where in the last term } h = \hat h(\tau).
\]
Hence,
\[
  \int_0^\tau b(\tau') d \tau'
  = \int_0^\tau O(h^{-2}) d \tau' + O(\varepsilon h^{-2} \ln h)
  = O(h^{-1} \ln h).
\]
Note that $O(\varepsilon h^{-2} \ln h)$ is $O(h^{-1} \ln h)$ as $\varepsilon h^{-1} < C_2^{-1}$.

As $\Delta(0) = O(\varepsilon^2)$, this gives the estimate
$\Delta (\tau) = O(\varepsilon^2 h^{-1}\ln h)$.
As $\norm{\pdv{h}{r}} = O(\ln^{-1} h)$ (here $h \sim \hat h(\tau)$ by~\eqref{e:2-times}), we have
$|\overline h(\tau) - \hat h(\tau)| = O(\varepsilon^2 h^{-1})$.
\begin{equation} \label{e:est-Delta}
  |\overline h(\tau) - \hat h(\tau)| = O(\varepsilon^2 h^{-1}),
  \qquad
  |\overline w(\tau) - \hat w(\tau)| = O(\varepsilon^2 h^{-1}\ln h).
\end{equation}
From the estimate on $\overline h(\tau) - \hat h(\tau)$ we have just proved and~\eqref{e:bound-h-eps} we get that
$\overline{h}(\tau) - \hat{h}(\tau) = C_2^{-1} O(\varepsilon) < C_2^{-2} O(\hat{h}(\tau)) < 0.5 \hat h(\tau)$ for large enough $C_2$, so the condition~\eqref{e:2-times} actually holds for all $t$ considered in this lemma.

The estimate for the difference in $w$ in~\eqref{e:est-Delta} can be improved. By~\eqref{e:est-Delta} and Table~\ref{t:est} for $h = \hat h(\tau)$, $\hat v = \hat v(\tau)$ and $\overline v = \overline v(\tau)$ we have
\[
  \norm{\hat f_{w, 1}(\overline v) + \varepsilon \hat f_{w, 2}(\overline v)
  + \varepsilon^2 \hat f_{w, 3}(\overline v, \overline \varphi(\tau), \varepsilon)
  - \hat f_{w, 1}(\hat v) - \varepsilon \hat f_{w, 2}(\hat v)
  }
  = O(\varepsilon^2 h^{-2} \ln^{-1} h) + O_*(\varepsilon^2 h^{-2}).
\]
Arguing as above, we can estimate the integral of this expression $d\tau$:
\[
  |\overline w(\tau) - \hat w(\tau)| = O(\varepsilon^2 h^{-1}).
\]

Let us now prove the estimate for $\varphi$. Denote $\omega_{0, 1}(v) = \omega(v) + \varepsilon \omega_1(v)$. Then from~\eqref{e:init-avg-trimmed} we have
\[
\overline{\varphi}(t) - \overline \varphi(0)
= \varepsilon^{-1} \int_0^\tau \Big(
  \omega_{0, 1}(\overline v(\tau'))
  + \varepsilon^2 \hat f_{\varphi, 2}(\overline v(\tau'), \overline \varphi(\tau'))
\Big) d\tau'.
\]
From Table~\ref{t:est} and~\eqref{e:bound-h-eps} we have $\pdv{\omega_{0, 1}}{h} = O(h^{-1}\ln^{-2} h)$. We also have $\norm{\pdv{\omega_{0, 1}}{w}} = O(\ln^{-1} h)$. Thus from~\eqref{e:est-Delta} we have
$|\omega_{0, 1}(\overline v(\tau)) - \omega_{0, 1}(\hat v(\tau))| = O(\varepsilon^2 h^{-2} \ln^{-2} h)$.
From Table~\ref{t:est} we have
$\hat f_{\varphi, 2} = O(h^{-2}\ln^{-2} h) + O_*(h^{-2}\ln^{-1} h)$.
So
\[
  \overline{\varphi}(\tau) - \overline \varphi(0)
  = \varepsilon^{-1} \int_0^\tau \omega_{0, 1}(\hat v(\tau')) d\tau'
  + \varepsilon \int_0^\tau O(\hat h^{-2}(\tau') \ln^{-2} \hat h(\tau'))
  + O_*(\hat h^{-2}(\tau') \ln^{-1} \hat h(\tau')) d\tau'.
\]
The second integral can be estimated in the same way as $\int_{0}^{\tau} b(\tau') d \tau'$ above:
\[
  \varepsilon \int_0^\tau O(\hat h^{-2}(\tau') \ln^{-2} \hat h(\tau'))
  + O_*(\hat h^{-2}(\tau') \ln^{-1} \hat h(\tau')) d\tau
  = O(\varepsilon h^{-1} \ln^{-1} h).
\]
This proves the formula for $\varphi$.

\section{Cancellation lemma} \label{s:cancel}
In this section we prove the following lemma. It will be useful when we prove the formula for the pseudo-phase, because due to this lemma two terms will cancel out.
Denote $\omega_1(h, w) = \hat f_{\varphi, 1} = \overline f^0_{\varphi, 1}$ to match the notation in~\cite{neishtadtVasiliev2005}.
\begin{lemma} \label{l:cancel}
  Consider a solution $\hat h(\tau), \hat w(\tau)$ of the order~2 averaged system~\eqref{e:avg-form}. Take $\tau_1 < \tau_2 < \tau_*$ such that $\hat h(\tau_1)$ is small enough. Denote $h_1 = \hat h(\tau_1)$, $h_2 = \hat h(\tau_2)$, $w_* = \hat w(\tau_*)$,
  $\Theta_{3*} = \Theta_3(w_*)$. Then
  \begin{equation} \label{e:cancel}
    \int_{\tau_1}^{\tau_2} \omega_1(\hat h(\tau), \hat w(\tau)) d\tau
    = - \frac{2\pi}{\Theta_{3*}} \Bigl|_{\tau_1}^{\tau_2} u^0_{h, 1}(\hat h(\tau), \hat w(\tau), 0)
    + O(h_1^{1/2})
    + O(\varepsilon \ln^{-1} h_1).
  \end{equation}
\end{lemma}
\noindent Let us first estimate $\omega_1$. Denote $\mathcal I(h, w) = \int_0^{2\pi} t(\varphi) f^0_h(\varphi) d \varphi$.
\begin{lemma} \label{l:omega1-I1}
\begin{equation} \label{e:omega1-I1}
  \omega_1 = \frac{1}{T} \pdv{\mathcal I}{h} + O(h^{-1/2} \ln^{-1} h).
\end{equation}
\end{lemma}
\begin{proof}
  Integrating by parts, we can write
  \[
    2 \pi \omega_1 = \int_0^{2 \pi} f^0_\varphi d \varphi = \Bigl|_0^{2 \pi} \varphi f^0_\varphi - \int_0^{2 \pi} \varphi \pdv{f^0_{\varphi}}{\varphi} d\varphi.
  \]
  Using~\eqref{e:trace} and the equality $\frac 1 T  \pdv{}{h}\/(T f^0_h) = \pdv{f^0_h}{h} + \frac{1}{T} \frac{dT}{dh} f^0_h$, this rewrites as
  \[
    2 \pi \omega_1 = 2 \pi f^0_\varphi(h, w, 0)
    - \int_0^{2 \pi} \varphi \Div(f^0) d\varphi
    + \sum_{i=1}^k \int_0^{2 \pi} \varphi \bigg(
                  \pdv{f^0_{w_i}}{w_i} + \frac{1}{T} \pdv{T}{w_i} f^0_{w_i}
                  \bigg) d\varphi
    + \int_0^{2 \pi} \varphi \frac 1 T  \pdv{}{h}\/(T f^0_h) d\varphi.
  \]
  By Table~\ref{t:est} the first term is $O(h^{-1/2} \ln^{-1} h)$. The second term is $O(1)$ as $\Div(f^0)$ is bounded.
  The third term is $O(1)$ by Table~\ref{t:est}.
  As $\pdv{}{h}$ commutes with integrating by $\varphi$, we can rewrite the last term as
  $\frac 1 T \pdv{}{h} \int_0^{2 \pi} \varphi T f^0_h d\varphi = \frac {2\pi} {T} \pdv{\mathcal I}{h}$.
  We have obtained~\eqref{e:omega1-I1}.
\end{proof}

\begin{lemma} \label{l:d-I-w}
  \begin{equation}
    \pdv{}{w_i} \Big( \int_{0}^{2\pi} t(\varphi) f^0_h d\varphi \Big) = O(1), \qquad i = 1, \dots k.
  \end{equation}
\end{lemma}
\begin{proof}
As $t=(2\pi)^{-1}T\varphi$, we have
\begin{align}
\begin{split}
  &\pdv{}{w_i} \Big( \int_{0}^{2\pi} t(\varphi) f^0_h d\varphi \Big) =
  (2\pi)^{-1} \pdv{T}{w} \int_{0}^{2\pi} \varphi f^0_h d\varphi
  + (2\pi)^{-1} T \int_{0}^{2\pi} \varphi \pdv{f^0_h}{w} d\varphi = \\
  &= T^{-1} \pdv{T}{w} \int_{0}^{2\pi} \varphi O_*(1) dt
  + \int_{0}^{2\pi} \varphi O_*(1) dt = O(1).
\end{split}
\end{align}
\end{proof}

\begin{proof}[Proof of Lemma~\ref{l:cancel}]
For small enough $h_1$ the value of $\hat h(\tau)$ decreases, so we may use $h = \hat h(\tau)$ as a coordinate along the solution of the averaged system. We will also take $\dv{\tau}{h}$, $\dv{w}{h}$ and $\dv{\mathcal I}{h}$ along this solution.
For convenience let us recall~\eqref{e:d-tau-h} here:
\[
  \dv{\tau}{h} = - \frac{T}{\Theta_3(\hat w)} (1 + O(\hat h \ln \hat h) + O(\varepsilon)) = O(\ln h).
\]
By Lemma~\ref{l:d-I-w} we have $\norm{\pdv{\mathcal I}{w}} = O(1)$. We can write
\begin{equation}
  \norm{\dv{w}{h}} = \norm{\dv{\tau}{h}(\hat f_{w, 1} + \varepsilon \hat f_{w, 2})} = O(\ln h) + \varepsilon O(h^{-1} \ln^{-2} h),
\end{equation}
\begin{equation}
  \dv{\mathcal I}{h} = \pdv{\mathcal I}{h} + \pdv{\mathcal I}{w} \dv{w}{h}
  = \pdv{\mathcal I}{h} + O(\ln h) + \varepsilon O(h^{-1} \ln^{-2} h).
\end{equation}
As $\omega_1 = O(h^{-1} \ln^{-3} h)$ and so $\int_0^{h_1} \abs{T\omega_1} dh = O(\ln^{-1} h_1)$, we have
\begin{equation}
  \int_{\tau_1}^{\tau_2} \omega_1(\hat h(\tau), \hat w(\tau)) d\tau
  = - \int_{h_1}^{h_2} \frac{T \omega_1}{\Theta_3(\hat w)} dh
  + O(h_1) + O(\varepsilon \ln^{-1} h_1).
\end{equation}
Integrating the estimate for $\norm{\dv{w}{h}}$, we get
\begin{equation} \label{e:w-star}
  \norm{\hat w - w_*} = O(h_1 \ln h_1) + O(\varepsilon \ln^{-1} h_1).
\end{equation}
As by~\cite[Lemma~3.2]{neishtadt17}
\begin{equation} \label{e:d-theta-w}
  \pdv{\Theta_3}{w} = O(1),
\end{equation}
this means
\begin{equation} \label{e:theta-star}
  \abs{\Theta_3(\hat w) - \Theta_{3*}} = O(h_1 \ln h_1) + O(\varepsilon \ln^{-1} h_1)
\end{equation}
and
\begin{equation}
  \int_{\tau_1}^{\tau_2} \omega_1(\hat h(\tau), \hat w(\tau)) d\tau
  = - \frac 1 {\Theta_{3*}} \int_{h_1}^{h_2} T \omega_1 dh
  + O(h_1) + O(\varepsilon \ln^{-1} h_1).
\end{equation}
By Lemma~\ref{l:omega1-I1} this can be rewritten as
\begin{align} \label{e:int-omega1}
\begin{split}
  \int_{\tau_1}^{\tau_2} \omega_1 d\tau
  &= - \frac 1 {\Theta_{3*}} \int_{h_1}^{h_2} \pdv{\mathcal I}{h} dh
  + O(h_1^{1/2}) + O(\varepsilon \ln^{-1} h_1) = \\
  &= - \frac 1 {\Theta_{3*}} \int_{h_1}^{h_2} \dv{\mathcal I}{h} dh
  + O(h_1^{1/2}) + O(\varepsilon \ln^{-1} h_1) = \\
  &= - \frac 1 {\Theta_{3*}} \Bigl|_{\tau_1}^{\tau_2} \mathcal I(\hat h(\tau), \hat w(\tau))
  + O(h_1^{1/2}) + O(\varepsilon \ln^{-1} h_1).
\end{split}
\end{align}

As $dt = \frac{T d\varphi}{2 \pi}$, by~\eqref{e:u-int-t} we have
\[
  u^0_{h, 1}(h, w, 0)
  = \frac{1}{T} \int_0^T \Big(t - \frac T 2\Big) f^0_h(t) dt
  = \frac{1}{2\pi} \int_0^{2\pi} t f^0_h(t) d\varphi - \frac 1 2 \int_0^T f^0_h dt.
\]
By~\cite[Corollary 3.1]{neishtadt17} $\int_0^T f^0_h(t) dt = -\Theta_3(w) + O(h \ln h)$.
Hence,
\[
  \Bigl|_{\tau_1}^{\tau_2} u^0_{h, 1}(\hat h(\tau), \hat w(\tau), 0) = \frac 1 {2\pi} \Bigl|_{\tau_1}^{\tau_2} \mathcal I(\hat h(\tau), \hat w(\tau))
  + O(h_1 \ln h_1) + O(\varepsilon \ln^{-1} h_1).
\]
Comparing this with~\eqref{e:int-omega1}, we get~\eqref{e:cancel}.
\end{proof}

\section{Proof of the formula for the pseudo-phase} \label{s:formula-proof}
In this section we prove the formula~\eqref{e:phase} for the pseudo-phase.
We use the notation from Section~\ref{s:formula}.
First let us prove some auxiliary statements.
\begin{lemma} \label{l:lim-u}
  We have $\lim \limits_{\tau \to \tau_* - 0} u_{h, 1}^0(\hat v(\tau), 0) = \frac 1 4 (\Theta_1(w_*) - \Theta_2(w_*)) = u_*$.
\end{lemma}
\begin{proof}
  Recall that $\Theta_2$ corresponds to $0 < \varphi < \pi$ and $\Theta_1$ to $\pi < \varphi < 2\pi$. For $\tau \to \tau_* - 0$ we have $\hat h(\tau) \to 0$.
  Let us split the integral expression~\eqref{e:u-int-phi} (with $f$ replaced by $f^0$) for~$u^0_{h, 1}(\hat v(\tau), 0)$ into the integrals over the part of the trajectory near $l_1$ and near $l_2$. For the first part the value of $\varphi(t) - \pi$ is close to $\pi/2$ far away from the saddle $C$. But close to $C$ we have $f^0_h \approx 0$, so the integral near $l_1$ is close to $\Theta_1 / 4$. Similarly, the integral near $l_2$ is close to $-\Theta_2 / 4$.
\end{proof}

\begin{lemma}
  Take $\tau_1 < \tau_*$, denote $h_1 = \hat h(\tau_1)$.
  Then we have
  \begin{equation} \label{e:int-omega}
    \int_{\tau_1}^{\tau_*} \omega(\hat v(\tau)) d \tau
    = \frac{2\pi}{\Theta_{3*}} h_1
    + O(\varepsilon h_1) + O(h_1^2 \ln h_1).
  \end{equation}
\end{lemma}
\begin{proof}
  As $\omega T = 2\pi$,~\eqref{e:d-tau-h} and~\eqref{e:theta-star} implies that
  \[
  \int_{\tau_1}^{\tau_*} \omega(\hat h(\tau)) d \tau
  = - 2\pi \int_{h_1}^{0} \Theta_3^{-1} \Big( 1 + O(\hat h \ln \hat h) + O(\varepsilon) \Big) d\hat h
  = - \frac{2\pi}{\Theta_{3*}} \int_{h_1}^{0} \Big( 1 + O(\hat h \ln \hat h) + O(\varepsilon) \Big) d\hat h,
  \]
  which gives the required estimate.
\end{proof}

\begin{lemma} \label{l:delta-h-near-sep}
  Assume $\Theta_1, \Theta_2 > 0$. Then there exist $c_1, c_2 > 0$ such that for all small enough $\varepsilon$ the following holds. Take a point $(h_0, w_0, 0)$ on the transversal $\varphi = 0$ with $\varepsilon \Theta_3(w_0) + 2 c_1 \varepsilon^{3/2} \le h_0 < c_2$.
  Then the orbit of this point intersects the transversal $\varphi = 0$ once more with
  \begin{equation} \label{e:delta-h-near-sep}
    h = h_0 - \varepsilon \Theta_3(w_0) + O(\varepsilon h_0 \ln h_0) + O(\varepsilon^2 h_0^{-1/2}) > c_1 \varepsilon^{3/2}
  \end{equation}
  and the time passed between these two intersections is $O(\ln h_0)$.
\end{lemma}
{ \color{myblue}
\noindent This lemma is proved in Appendix~\ref{s:near-sep}, the proof is similar to the proof of~\cite[Proposition~5.1]{neishtadt17}.
}

\noindent Recall that $h_{-2}, h_{-3}, \dots$ denote the values of $h$ at the consecutive crossings of the transversal $\varphi = 0$ before $h_{-1}$.
\begin{lemma} \label{l:summation}
  For $\varepsilon^{0.9} <  h_{-n} < \varepsilon^{0.1}$ we have
  \begin{equation} \label{e:n-last-crossings}
    h_{-n} = h_{-1} + \varepsilon (n - 1) \Theta_3(w_{-n})  + O(h_{-n}^2 \ln h_{-n}) + \varepsilon O(h_{-n}^{1/2}).
  \end{equation}
\end{lemma}
\begin{proof}
First, let us note that as the considered $h_{-i}$ are in $[c_1 \varepsilon^{3/2}, \varepsilon^{0.1}]$, we have $\ln h \sim \ln \varepsilon \sim \ln h_{-n}$. So the time passed between two consecutive intersections is $O(\ln h_{-n})$. As $n \sim \varepsilon^{-1} h_{-n}$, the total time between the moments corresponding to $h_{-n}$ and $h_{-1}$ is $O(\varepsilon^{-1} h_{-n} \ln h_{-n})$. As $\dot w = O(\varepsilon)$, we have $\norm{w - w_{-n}} = O(h_{-n} \ln h_{-n})$ for all encountered values of $w$ and by~\eqref{e:d-theta-w} we have
$\Theta_3(w) - \Theta_3(w_{-n}) = O(h_{-n} \ln h_{-n})$.
Now the required estimate follows from Lemma~\ref{l:delta-h-near-sep} by summation.
\end{proof}

Let us return to the proof of the formula for pseudo-phase.
Denote by
$v(\tau), \varphi(\tau) = v(t), \varphi(t)$, where $t = \varepsilon^{-1} \tau$ and $v(\tau) = (h(\tau), w(\tau))$, the solution of the perturbed system written using the slow time $\tau$.
We denote by $\overline v(\tau), \overline \varphi(\tau)$, where $\overline v(\tau) = (\overline h(\tau), \overline w(\tau))$, the solution $v(\tau), \varphi(\tau)$ of the perturbed system, written in the averaged chart~\eqref{e:varchange}. Denote $\overline v_0 = (\overline h_0, \overline w_0) = \overline v(0)$, $\overline \varphi_0 = \overline \varphi(0)$. We have $\norm{\hat v_0 - \overline v_0} = O(\varepsilon^2)$, so we may use Lemma~\ref{l:close-avg-init}.
\begin{lemma}
  There exists $C_3 > 0$ such that for all $\tau$ with
  \begin{equation} \label{e:h-larger-eps}
    \hat h(\tau) > C_3 \varepsilon
  \end{equation}
  the solutions $v(\tau)$, $\overline v(\tau)$ and $\hat v(\tau)$ are close:
  \begin{equation} \label{e:h-close}
    \norm{\hat v - \overline v} = O(\varepsilon^2 \hat h^{-1}) < \hat h/4,
    \qquad
    \norm{\overline v - v} = O(\varepsilon) < \hat h/4.
  \end{equation}
\end{lemma}
\begin{proof}
  The first two estimates are given by Lemma~\ref{l:close-avg-init}.
  To obtain the last one, we just plug the estimates from Table~\ref{t:est} into the equation $v = \overline v + \varepsilon u_{v, 1} + \varepsilon^2 u_{v, 2}$ from~\eqref{e:varchange}.
\end{proof}

Now we are ready to prove~\eqref{e:phase}.
Consider a moment $\tau_1$ such that $\varphi(\tau_1) = 0$ and $\hat h(\tau_1)$ is as close as possible to $\varepsilon^{2/3} \ln^{-1/3} \varepsilon$. Note that we have~\eqref{e:h-larger-eps} for $\tau = \tau_1$. We may check that under the condition~\eqref{e:h-larger-eps} the difference between $\hat h(\tau)$ for consecutive times $\tau$ with $\varphi(\tau) = 0$ is $O(\varepsilon)$. Indeed, the time between consecutive fast times of crossing the transversal $\varphi = 0$ is $O(T)$ and $\dot{\hat{h}}$ is $O(T^{-1})$. Hence,
\begin{equation} \label{e:pick-h-k}
  \hat h(\tau_1) = (1 + o(1)) \varepsilon^{2/3} \ln^{-1/3} \varepsilon.
\end{equation}
Denote by
  $h_1, \varphi_1, \overline h_1, \overline \varphi_1, \hat h_1, \hat v_1$
the values of
$h, \varphi, \overline h, \overline \varphi, \hat h, \hat v$
at the slow time $\tau_1$.
As justified by~\eqref{e:h-close}, we may write $h_1$ instead of $\hat h_1$ and $\overline h_1$ in the error terms. For brevity let us even denote $h = h_1$ for the error terms and write simply $O(h)$.

\begin{lemma}
  For any $\tau \ge \tau_1$ until separatrix crossing (i.e. with $h(\tau) > 0$) we have
  \begin{equation} \label{e:theta-and-w}
    w(\tau) = w_* + O(h_1 \ln h_1),
    \qquad
    \Theta_3(w(\tau)) = \Theta_{3*} + O(h_1 \ln h_1).
  \end{equation}
\end{lemma}
\begin{proof}
  By~\eqref{e:w-star} we have $\hat w(\tau_1) = w_* + O(h_1 \ln h_1) + O(\varepsilon \ln^{-1} h_1)$. By~\eqref{e:h-close} we have $w(\tau_1) = \hat w(\tau_1) + O(\varepsilon^2 h_1^{-1})$. Arguing as in the proof of Lemma~\ref{l:summation} gives
  $w(\tau) = w(\tau_1) + O(h_1 \ln h_1)$. Combining these estimates gives the first statement (the term $O(h_1 \ln h_1)$ absorbs other terms for $h_1 \gtrsim \varepsilon$); the second statement follows from the first by~\eqref{e:d-theta-w}.
\end{proof}
{ \color{myblue}
\begin{remark} \label{r:est-h-m-1}
  We have $0 \le h_{-1} < \varepsilon \Theta_{i*} + O(\varepsilon^{3/2})$.
  Indeed, by~\cite[Proposition~5.1]{neishtadt17} if $h > \Theta_i + c \varepsilon^{3/2}$ for some $c$ at a crossing of transversal $\varphi=0$, there will be next crossing of this transversal. By~\eqref{e:theta-and-w} $\Theta_{i*}$ and $\Theta_{i}$ (with $h \sim \varepsilon$) are $O(h_1 \ln h_1)$-close, this is $\ll \varepsilon^{1/2}$ by~\eqref{e:pick-h-k}.
\end{remark} }

\noindent Let us split the integral in~\eqref{e:phase} into integrals from $0$ to $\tau_1$ and from $\tau_1$ to $\tau_*$. First, let us check that
\begin{equation} \label{e:tau-0-to-k}
  \varphi_0 + \frac{1}{\varepsilon} \int_{\tau = 0}^{\tau_1} \bigg(
    \omega(\hat v(\tau)) + \varepsilon \omega_1(\hat v(\tau))
  \bigg) d\tau
  = 2 \pi m + O(\varepsilon h^{-1} \ln^{-1} h),
\end{equation}
where $m \in \mathbb Z$.
By Lemma~\ref{l:close-avg-init} we have
\begin{equation}
  \frac{1}{\varepsilon} \int_{\tau = 0}^{\tau_1} \bigg(
    \omega(\hat v(\tau)) + \varepsilon \omega_1(\hat v(\tau))
  \bigg) d\tau = \overline \varphi_1 - \overline \varphi_0
  + O(\varepsilon h^{-1} \ln^{-1} h).
\end{equation}
We also have $\varphi = \overline \varphi + \varepsilon u_{\varphi, 1}$. By Table~\ref{t:est} $u_{\varphi, 1} = O(h^{-1} \ln^{-1} h)$, so
$\overline \varphi_1 - \overline \varphi_0
= \varphi_1 - \varphi_0 + O(\varepsilon h^{-1} \ln^{-1} h)$.
As $\varphi_1$ = $2\pi m$, this gives the required equality~\eqref{e:tau-0-to-k}.

Now let us use~\eqref{e:int-omega} and~\eqref{e:cancel} (in~\eqref{e:cancel} we pass to the limit for $\tau_2 \to \tau_* - 0$, by Lemma~\ref{l:lim-u} we have $u^0_{h, 1}(\hat v(\tau_2), 0) \to u_*$) to compute the remaining terms in~\eqref{e:phase}.
We have
\begin{align} \label{e:tau-k-tau-star}
  &\frac{1}{2\pi \varepsilon} \bigg(
    \int_{\tau=\tau_1}^{\tau_*} \big(
      \omega(\hat{v}(\tau)) + \varepsilon \omega_1(\hat{v}(\tau))
    \big) d\tau
 \bigg)
 + \frac{u_*}{\Theta_{3*}} = \\
 &\qquad= \frac{1}{\varepsilon \Theta_{3*}} \bigg(
    \hat h_1 + \varepsilon u_{h, 1}^0(\hat v_1, 0)
 \bigg)
 + O(h^{1/2}) + O(\varepsilon^{-1} h^2 \ln h). \nonumber
\end{align}
Note that the term $O(\varepsilon \ln^{-1} h)$ from~\eqref{e:cancel} is absorbed into $O(h^{1/2})$ by~\eqref{e:h-larger-eps}.
As $h_1 \gtrsim \varepsilon$, by Table~\ref{t:est} we have $\norm{\pdv{}{v}\/ \varepsilon u^0_{h, 1}} = O(\ln^{-1} h)$. Hence, by~\eqref{e:h-close} we have
\[
  \hat h_1 + \varepsilon u^0_{h, 1}(\hat v_1, 0)
  = \overline h_1 + \varepsilon u^0_{h, 1}(\overline v_1, 0)
  + O(\varepsilon^2 h^{-1})
  = \overline h_1 + \varepsilon u_{h, 1}(\overline v_1, 0, \varepsilon)
  + O(\varepsilon^2 h^{-1}).
\]
The last equality is justified by Lemma~\ref{l:est-f-hat}. The error term $O(\varepsilon^2)$ appears, but it is absorbed into $O(\varepsilon^2 h^{-1})$.
As $0 = \varphi_1 = \overline \varphi_1 + \varepsilon u_{\varphi, 1}$, by Table~\ref{t:est} we have
$\overline \varphi_1 = O(\varepsilon h^{-1} \ln^{-1} h).$
Hence, by the estimate
$\pdv{u_{h, 1}}{\varphi} = O(\ln h)$
from Table~\ref{t:est} we get
\[
  \varepsilon u_{h, 1}(\overline v_1, 0, \varepsilon)
  = \varepsilon u_{h, 1}(\overline v_1, \overline \varphi_1, \varepsilon)
  + O(\varepsilon^2 h^{-1})
\]
and
\[
  \hat h_1 + \varepsilon u^0_{h, 1}(\hat v_1, 0)
  = \overline h_1 + \varepsilon u_{h, 1}(\overline v_1, \overline \varphi_1, \varepsilon)
  + O(\varepsilon^2 h^{-1}).
\]
As by~\eqref{e:varchange}
$$
  h_1 = \overline h_1 + \varepsilon u_{h, 1}(\overline v_1, \overline \varphi_1, \varepsilon) + \varepsilon^2 u_{h, 2}(\overline v_1, \overline \varphi_1, \varepsilon),
$$
the estimate $\varepsilon^2 u_{h, 2} = O(\varepsilon^2 h^{-1})$ from Table~\ref{t:est} yields
\[
  \hat h_1 + \varepsilon u^0_{h, 1}(\hat v_1, 0)
  = h_1 + O(\varepsilon^2 h^{-1}).
\]
Combining this with~\eqref{e:tau-k-tau-star}, we get
\begin{equation} \label{e:tau-k-tau-star-final}
  \frac{1}{2\pi \varepsilon} \bigg(
    \int_{\tau=\tau_1}^{\tau_*} \big(
      \omega(\hat{v}(\tau)) + \varepsilon \omega_1(\hat{v}(\tau))
    \big) d\tau
  \bigg)
  + \frac{u_*}{\Theta_{3*}}
  = \frac{h_1}{\varepsilon \Theta_{3*}}
  - R(h_1)
\end{equation}
with the error term
\[
  R = O(h^{1/2}) + O(\varepsilon h^{-1}) + O(\varepsilon^{-1} h^2 \ln h).
\]
After taking a sum with~\eqref{e:tau-0-to-k}, we get
\begin{equation}
 \frac{h_1}{\varepsilon \Theta_{3*}}
 = \frac{1}{2\pi}
   \bigg(
     \varphi_0 + \frac 1 \varepsilon \int_{\tau=0}^{\tau_*} \big(
       \omega(\hat{v}(\tau)) + \varepsilon \omega_1(\hat{v}(\tau)) d\tau
     \big)
   \bigg)
   + \frac{u_*}{\Theta_{3*}}
   - m
   + R(h_1).
\end{equation}
Note that $R$ absorbs the error term in~\eqref{e:tau-0-to-k}. Let us now apply~\eqref{e:n-last-crossings} for $v_{-n} = v_1$. We have $n \sim \varepsilon^{-1} h_1$; by~\eqref{e:theta-and-w} we have
$(n-1)\Theta_3(w_1) = (n-1)\Theta_{3*} + O(\varepsilon^{-1} h_1^2 \ln h_1)$,
so this yields the required formula~\eqref{e:phase}, but with the error term $R(h_1)$ depending on $h_1$. Note that the error term above and the error term in~\eqref{e:n-last-crossings} divided by $\varepsilon$ are not greater than $R$. Then we just plug in the expression~\eqref{e:pick-h-k} for $h_1$ and obtain $R = O(\varepsilon^{1/3} \ln^{1/3} \varepsilon)$. One may check that~\eqref{e:pick-h-k} minimizes the error term. Indeed, first we check that up to some power of $\ln \varepsilon$ the value of $R$ is minimal for $h \approx \varepsilon^{2/3}$. Then $\ln h \approx (2/3) \ln \varepsilon$, and from this we see that $R$ is minimal for $h$ given by~\eqref{e:pick-h-k}.
This completes the proof of formula~\eqref{e:phase}. \qed

{ \color{myblue}
So far we have proved formula~\eqref{e:phase} for a specific choice of transversal $\varphi=0$ tangent to the bisector of the angle between separatrices.
\begin{remark} \label{r:phase-transversal}
  Formula~\eqref{e:phase} holds for any transversal $\varphi=0$ tangent to the bisector of the angle between separatrices.
\end{remark}
\begin{proof}
  Suppose we have other transversal, denote by $\psi$ the angle variable such that $\psi=0$ on this other transversal.
  Let us show that the right-hand side of~\eqref{e:phase} is the same for $\varphi$ and $\psi$ with $O(\varepsilon)$ accuracy. By Lemma~\ref{l:change-phi} the evolution of slow variables $\hat v(\tau)$ for order 2 averaged system is the same for $\varphi$ and $\psi$, so $w_*$ and $\int \omega d \tau$ (for brevity, we omit the limits of integration) are also the same. Only the terms $\varphi_0$ and $\int \omega_1 d \tau$
  are different. Let us show that their sum is almost the same (up to $O(\varepsilon)$) using~\eqref{e:change-transversal} with $\tau_1 = 0$ and $\tau_2 = \tau_*$.
  The term $\Delta \varphi$~\eqref{e:change-transversal} at $\tau_1 = 0$ is the difference between $\varphi_0$ for $\varphi$ and $\psi$.
  Close to separatrices the difference in $\varphi$ between our two transversals tends to zero by~\eqref{e:delta-t-tangent}, so at $\tau_2=\tau_*$ we have $\Delta \varphi = 0$.
  The error term in~\eqref{e:change-transversal} is $O(\varepsilon)$, as verified in the proof of Remark~\ref{r:any-transversal}. This means that $\varphi_0 + \int \omega_1 d\tau$ is the same for $\varphi$ and $\psi$ with $O(\varepsilon)$ accuracy.

  Denote by $h_{-1}^\psi$ and $\xi^\psi$ the values of $h^{-1}$ and $\xi$, respectively, computed using $\psi$ instead of $\varphi$. By~\eqref{e:delta-t-tangent} the time for unperturbed system between transversals $\varphi=0$ and $\psi=0$ is $O(\sqrt{h})$. One can check that this estimate also holds for perturbed system in the same way, as the estimate $h_{-1} \gtrsim \varepsilon^{3/2}$ implies that the amplitude of the vector field of the unperturbed system near the bisector is $\gtrsim \sqrt{h_{-1}} \gtrsim \varepsilon^{3/4} \gg \varepsilon$, i.e., much greater than the amplitude of the perturbation.
  The trajectory between the two transversals lies $O(\sqrt{h})$-close to the saddle, thus we have $f_h = O(\sqrt{h_{-1}})$.
  As $\dot h = \varepsilon f_h = O(\varepsilon \sqrt{h_{-1}})$, this estimate for the time implies that $h_{-1} - h_{-1}^\psi = O(\varepsilon h_{-1}) = O(\varepsilon^2)$ and $\xi - \xi^\psi = O(\varepsilon)$.

  Thus~\eqref{e:phase} holds for any transversal $\varphi=0$ tangent to the bisector.
\end{proof}
}

\section{Probabilities (proofs)} \label{s:prob-proofs}
A trajectory starting in $G_3$ may be captured into $G_1$ or $G_2$ after separatrix crossing with the outcome determined by the pseudo-phase as stated in~\cite[Proposition~5.1]{neishtadt17}. Let us state a corollary of this proposition here.
\begin{corollary} \label{c:cature-pseudo-phase}
  Any solution of the perturbed system
  with the pseudo-phase in
  $[O(\varepsilon^{1/2}), \frac{\Theta_2}{\Theta_3} - O(\varepsilon^{1/2})]$ is captured in $G_2$
  and with the pseudo-phase in
  $[\frac{\Theta_2}{\Theta_3} + O(\varepsilon^{1/2}), 1 - O(\varepsilon^{1/2})]$ is captured in $G_1$.
\end{corollary}
{ \color{myblue}
\noindent The general reasoning in~\cite[Proposition~5.1]{neishtadt17} is as follows. The last wind before separatrix crossing starts with $h = \Theta_3 \psi$, where $\psi$ is the pseudo-phase. After solution passes near $l_2$ (cf. Figure~\ref{f:up}), $h$ decreases by $\approx \Theta_2$. If after that $h<0$ (i.e, $\psi < \Theta_2/\Theta_3$), we have capture in $G_2$; otherwise, in $G_1$.
}

\begin{proof}[Proof of Proposition~\ref{p:prob-Anosov}]
  As $\mu(U_1(\varepsilon_0)) + \mu(U_2(\varepsilon_0)) \le \varepsilon_0$,
  it is enough to show that
  \begin{equation} \label{e:est-prob-capture-one-side}
    \frac{\mu(U_j(\varepsilon_0))}{\varepsilon_0}
    >
    \frac{\Theta_j}{\Theta_3} + O(\varepsilon_0^{1/3} \ln^{1/3} \varepsilon_0).
  \end{equation}
  Denote by $\psi(\varepsilon)$ the right hand side of~\eqref{e:phase} without the fractional part and the error term. Note that the integrals in~\eqref{e:phase} are computed along the solution of the averaged system of order 2 and so they depend on $\varepsilon$.
  Denote $a = (2\pi)^{-1} \int_{\tau=0}^{\tau_*} \omega(\hat v_1(\tau)) d\tau$, where $\hat v_1(\tau)$ is the solution of the first order averaged system~\eqref{e:avg-1-order} with the initial condition $\hat v_1(0) = \hat v(0)$.
  We may check that
  \begin{equation}
    \psi = \varepsilon^{-1}a + O(1), \qquad
    \dv{\psi}{\varepsilon} = - \varepsilon^{-2} a + O(\varepsilon^{-1}).
  \end{equation}
  Hence, for small enough $\varepsilon_0$ we have $\psi \to \infty$ monotonically when $\varepsilon$ decreases from $\varepsilon_0$ to $0$ and
  \begin{equation} \label{e:d-eps-psi}
    \dv{\varepsilon}{\psi} = -\psi^{-2}a + O(\psi^{-3}).
  \end{equation}
  Without loss of generality we can take $j = 1$ in~\eqref{e:est-prob-capture-one-side}. By~Corollary~\ref{c:cature-pseudo-phase} and~\eqref{e:phase} for all $\varepsilon$ such that $\psi \in [n + \frac{\Theta_2}{\Theta_3} + O(\varepsilon^{1/3} \ln^{1/3} \varepsilon), \; n + 1 - O(\varepsilon^{1/3} \ln^{1/3} \varepsilon)]$, $n \in \mathbb N$, the trajectory with the initial condition $v_0, \varphi_0$ is captured in $G_1$.  Using~\eqref{e:d-eps-psi}, we can estimate the length of the union of the preimages of such segments for the map $\varepsilon \mapsto \psi$ as follows ($n_0 \sim \varepsilon_0^{-1}$ in the formula below):
  \begin{equation}
    \mu(U_1(\varepsilon_0)) \ge a \sum_{n \ge n_0} n^{-2} \Big(\frac{\Theta_1}{\Theta_3} + O(n^{-1/3} \ln^{1/3} n) \Big)
    + O(\varepsilon_0^2).
  \end{equation}
  On the other hand, we have
  \begin{equation}
    \varepsilon_0 + O(\varepsilon_0^2) = a \sum_{n \ge n_0} (n^{-2} + O(n^{-3})),
  \end{equation}
  as the union of the preimages of the segments $[n, n+1]$, $n \ge n_0$ is $(0, \varepsilon_0 + O(\varepsilon_0^2)]$. These two formulas imply~\eqref{e:est-prob-capture-one-side}.
\end{proof}

\begin{remark} \label{r:prob-Arnold}
  Let us sketch how formula~\eqref{e:phase} implies the formula~\eqref{e:prob-capture} for the probability of capture for Definition~\ref{d:prob-Arnold}.
  Given initial data $v_0, \varphi_0$, let us fix some $\varphi$ close to $\varphi_0$ and vary $v$ near $v_0$. Denote by $\xi(v)$ the pseudo-phase of the solution of the perturbed system with the initial condition $v, \varphi$.  From~\eqref{e:phase} we have $\norm{\dv{\xi}{v}} \sim \varepsilon^{-1}$. Using~\eqref{e:phase} and Corollary~\ref{c:cature-pseudo-phase}, it is possible to show that most points in the neighborhood of $v_0$ are covered by disjoint interchanging stripes of width $O(\varepsilon)$ formed by values of $v$ such that the trajectory is captured in $G_1$ and $G_2$, and the widths of the stripes captured into $G_j$ are proportional to $\Theta_j$. Then, naturally, the relative measure of the values of $v$ captured into $G_j$ is $\frac{\Theta_j}{\Theta_3}$. Integrating this by $\varphi$, we get the formula for the probability of capture.
\end{remark}

\section{Formulas for the averaging chart} \label{s:varchange-formulas}
In this section we present formulas for $\overline f_{\varphi, 2}$ and $\overline f_{h, 3}$ from Lemma~\ref{l:varchange-short} and prove this lemma. We use the notation introduced in Section~\ref{s:varchange}. We will also need the following notation.
\begin{itemize}
  \item Denote by $x$ the column vector $(h, w, \varphi)$ and by
    $\overline x$ the column vector $(\overline h, \overline w, \overline \varphi)$. Let
    $f_v = (f_h, f_w), \;
    \overline f_{x, i} = (\overline f_{h, i}, \overline f_{w, i}, \overline f_{\varphi, i}), \;
    u_{x, i} = (u_{h, i}, u_{w, i}, u_{\varphi, i})$.

  \item Given $k = x, v, h, \varphi$, let us denote $u_{k, 1, 2} = u_{k, 1} + \varepsilon u_{k, 2},
    \; \overline f_{k, 1, 2} = \overline f_{k, 1} + \varepsilon \overline f_{k, 2},
    \; \overline f_{k, 2, 3} = \overline f_{k, 2} + \varepsilon \overline f_{k, 3},
    \; \overline f_{k, 1, 2, 3} = \overline f_{k, 1} + \varepsilon \overline f_{k, 2}
    + \varepsilon^2 \overline f_{k, 3}$. For $k=x$ the terms $u_{\varphi, 2}, \; \overline f_{\varphi, 3}$ appear, we set $u_{\varphi, 2} = \overline f_{\varphi, 3} = 0$.

  \item Given a vector-function $g(x) = (g_1, \dots, g_l)$, denote
    $\big( \pdv{g}{x} \big)_{int} = (\pdv{g_1}{x}\/(\xi_1), \dots, \pdv{g_l}{x}\/(\xi_l))$, \;
    $\big( \pdv[2]{g}{x} \big)_{int} = (\pdv[2]{g_1}{x}\/(\eta_1), \dots, \pdv[2]{g_l}{x}\/(\eta_l))$,
    where $\xi_i, \eta_i$ are some intermediate points on the segment $[x, \overline x]$.
\end{itemize}
\begin{lemma} \label{l:varchange}
  \noindent We have the following system of linear equations determining $\overline f_{\varphi, 2}$ and $\overline f_{v, 3} = (\overline f_{h, 3}, \overline f_{w, 3})$:
  \begin{align}
  \begin{split} \label{e:f3}
     &(1 + \varepsilon \pdv{u_{\varphi, 1}}{\varphi})\overline f_{\varphi, 2}
     + \varepsilon^2 \pdv{u_{\varphi, 1}}{v} \overline f_{v, 3} = \\
     &\quad = \pdv{\omega}{v} u_{v, 2}
     + \frac 1 2 u_{v, 1, 2}^T \Big(\pdv[2]{\omega}{v}\Big)_{int} u_{v, 1, 2}
     + \Big(\pdv{f_\varphi}{x}\Big)_{int} u_{x, 1, 2}
     - \pdv{u_{\varphi, 1}}{v} \overline f_{v, 1, 2}
     - \pdv{u_{\varphi, 1}}{\varphi} \overline f_{\varphi, 1}, \\
     &(1 + \varepsilon \pdv{u_{v, 1, 2}}{v}) \overline f_{v, 3}
     + \pdv{u_{v, 1, 2}}{\varphi} \overline f_{\varphi, 2} = \\
     &\quad = \pdv{f_v}{v} u_{v, 2}
     + \frac 1 2 u_{x, 1, 2}^T \Big(\pdv[2]{f_v}{x}\Big)_{int} u_{x, 1, 2}
     - \pdv{u_{v, 1}}{v} \overline f_{v, 2}
     - \pdv{u_{v, 2}}{v} \overline f_{v, 1, 2}
     - \pdv{u_{v, 2}}{\varphi} \overline f_{\varphi, 1}. \\
  \end{split}
  \end{align}
\end{lemma}

\begin{proof}[Proof of lemmas~\ref{l:varchange-short} and~\ref{l:varchange}]
  We shall differentiate the coordinate change~\eqref{e:varchange} with respect to the time and rewrite all emerging terms as functions of $\overline x$. For brevity the equations on $h$ and $w$ will be grouped together as an equation on $v$. The derivatives of the left hand sides of~\eqref{e:varchange} are given by~\eqref{e:init}. They are functions of $x$, let us write Taylor's expansions at the point $\overline x$. We group together the terms of order at least $3$ for the coordinate change in $v$ and $2$ for the change in $\varphi$
  \begin{align*}
    \dot{v} &= \varepsilon f_v(x) = \varepsilon f_v(\overline x + \varepsilon u_{x, 1} + \varepsilon^2 u_{x, 2}) = \\
    &= \varepsilon f_v(\overline x)
    + \varepsilon^2 \pdv{f_v}{x} u_{x, 1}
    + \varepsilon^3 \bigg(
      \pdv{f_v}{v} u_{v, 2}
      + \frac 1 2 u_{x, 1, 2}^T \Big(\pdv[2]{f_v}{x}\Big)_{int} u_{x, 1, 2}
    \bigg), \\
    \dot{\varphi} &= \omega(v) + \varepsilon f_{\varphi}(x)
    = \omega(\overline v + \varepsilon u_{v, 1} + \varepsilon^2 u_{v, 2})
    + \varepsilon f_{\varphi}(\overline x + \varepsilon u_{x, 1} + \varepsilon^2 u_{x, 2}) = \\
    &= \omega(\overline v)
    + \varepsilon \bigg(
      \pdv{\omega}{v} u_{v, 1}
      + f_{\varphi}(\overline x)
    \bigg) + \\
    & \hspace{32pt} + \varepsilon^2 \bigg(
      \pdv{\omega}{v} u_{v, 2}
      + \frac 1 2 u_{v, 1, 2}^T \Big(\pdv[2]{\omega}{v}\Big)_{int} u_{v, 1, 2}
      + \Big(\pdv{f_\varphi}{x}\Big)_{int} u_{x, 1, 2}
    \bigg).
  \end{align*}
  Now we write the terms containing the derivatives of $u_{k, i}$.
  \begin{align*}
    &\varepsilon \dot{u}_{v, 1}(\overline x) + \varepsilon^2 \dot{u}_{v, 2}(\overline x) = \\
    &= \varepsilon \pdv{u_{v, 1}}{\varphi} \omega
    + \varepsilon^2 \bigg(
      \pdv{u_{v, 2}}{\varphi} \omega
      + \pdv{u_{v, 1}}{x} \overline f_{x, 1}
    \bigg)
    + \varepsilon^3 \bigg(
      \pdv{u_{v, 1}}{x} \overline f_{x, 2, 3}
      + \pdv{u_{v, 2}}{x} \overline f_{x, 1, 2, 3}
    \bigg), \\
    &\varepsilon \dot{u}_{\varphi, 1}(\overline x)
    = \varepsilon \pdv{u_{\varphi, 1}}{\varphi} \omega
    + \varepsilon^2 \pdv{u_{\varphi, 1}}{x} \overline f_{x, 1, 2, 3}.
  \end{align*}
  Let us plug these expressions together with~\eqref{e:init-avg} into the time derivative of~\eqref{e:varchange}.
  Equating the terms of the same order in~$\varepsilon$ (grouping together the terms with order at least $3$ for the equation on $v$ and $2$ for the equation on $\varphi$), we get~\eqref{e:homological} and~\eqref{e:homological-rhs}, as well as the following equations:
\begin{align*}
   \overline f_{\varphi, 2} &= \pdv{\omega}{v} u_{v, 2}
   + \frac 1 2 u_{v, 1, 2}^T \Big(\pdv[2]{\omega}{v}\Big)_{int} u_{v, 1, 2}
   + \Big(\pdv{f_\varphi}{x}\Big)_{int} u_{x, 1, 2}
   - \pdv{u_{\varphi, 1}}{x} \overline f_{x, 1, 2, 3}, \\
   \overline f_{v, 3} &= \pdv{f_v}{v} u_{v, 2}
   + \frac 1 2 u_{x, 1, 2}^T \Big(\pdv[2]{f_v}{x}\Big)_{int} u_{x, 1, 2}
   - \pdv{u_{v, 1}}{x} \overline f_{x, 2, 3}
   - \pdv{u_{v, 2}}{x} \overline f_{x, 1, 2, 3}, \\
\end{align*}
which are equivalent to~\eqref{e:f3}, we just expand some terms like $\overline f_{x, 1, 2, 3}$ in order to move the terms containing $\overline f_{\varphi, 2}$ and $\overline f_{v, 3}$ to the left hand side.
\end{proof}

\begin{proof}[Proof of Lemma~\ref{l:u-int-t}]
The function $u_{a, 1}$ is uniquely determined by two properties.
The first one is that $\pdv{u_{a, 1}}{t} = f_a(t) - \langle f_a \rangle_t$
(this follows from~\eqref{e:homological},~\eqref{e:homological-rhs}).
Denote by $U$ the expression on the right hand side of~\eqref{e:u-int-t}. We have
\[
  \pdv{U}{t_0} = \frac{1}{T} \int_0^T \Big(t - \frac T 2\Big) \pdv{f_a}{t}\/(t+t_0) dt.
\]
Integrating by parts, this can be rewritten as
\[
  \pdv{U}{t_0}
  = \frac{1}{T} \Bigl|_{t=0}^T f_a(t+t_0) \Big(t - \frac T 2\Big) - \frac{1}{T} \int_0^T f_a(t+t_0) dt
  = f_a(t_0) - \langle f_a \rangle_t.
\]
Hence the first property of $u_{a, 1}$ holds for $U$.

The second property is that $\langle u_{a, 1} \rangle_t = 0$.
This also holds for $U$, it is checked by writing $\int U(t_0)dt_0$ as a double integral
and changing the order of integration.
\end{proof}

\section{Estimates related to the energy-angle variables} \label{s:est-action-angle}
\subsection{The coordinates \texorpdfstring{$\tilde h, \; \tilde w, \; \tilde t_i$}{tilde h, w, ti}} \label{s:Moser}
Our goal in this section is to estimate how $q, p$ (or, more generally, a smooth function $\psi(q, p, z)$) depend on  $h, w, \varphi$ for $h \to 0$. To do so, we introduce new coordinates $\tilde h, \tilde w, \tilde t_i$. The subscript $i$ is here because there will be different coordinate systems in different parts of the phase space. Then we will estimate how $q, p$ depend on $\tilde h, \tilde w, \tilde t_i$ and how $\tilde h, \tilde w, \tilde t_i$ depend on $h, w, \varphi$. Combining these estimates, we will get the required estimates of the dependence of $q, p$ on $h, w, \varphi$.

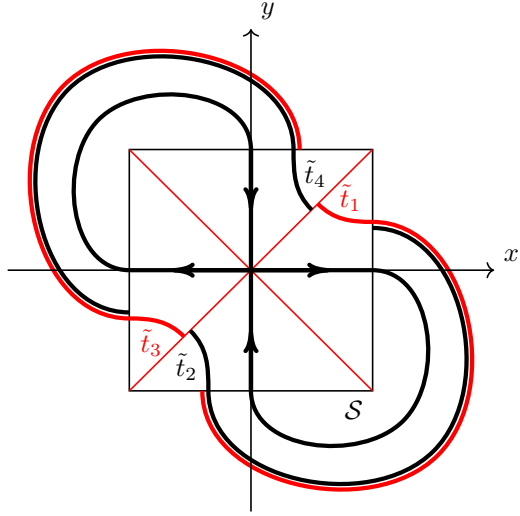
\begin{figure}[H]
  \centering
  \begin{tikzpicture}[scale=0.8]
    \draw[semithick] (-2,2) -- (-2,-2) -- (2, -2) node[below left]{$\mathcal S$} --
      (2, 2) -- cycle;
    \draw[semithick, red] (2, 2) -- (-2, -2);
    \draw[semithick, red] (-2, 2) -- (2, -2);
    \draw[ultra thick] (0, 2) -- (0, -2) to[out=270, in=225] (2.5, -2.5) %
      to [out=45, in=0] (2, 0) -- (-2, 0) to[out=180, in=225] (-2.5, 2.5)
      to [out=45, in=90] (0, 2);
    \draw[red, ultra thick] (1.1, 1.1) to[out=-45, in=180] (2, 0.8) node[above left]{$\tilde t_1$}
      to [out=0, in=45] (3.07, -3.07) to[out=-135, in=-90] (-0.8, -2);
    \draw[ultra thick] (2, 0.7) to [out=0, in=45] (3, -3) %
      to[out=-135, in=-90] (-0.7, -2) node[above left]{$\tilde t_2$} to [out=90, in=-45] (-1, -1);
    \draw[red, ultra thick] (-1.1, -1.1) to[out=135, in=0] (-2, -0.8) node[below right]{$\tilde t_3$}
      to[out=180, in=-135] (-3.07, 3.07) to [out=45, in=90] (0.8, 2);
    \draw[ultra thick] (-2, -0.7) to[out=180, in=-135] (-3, 3) %
      to [out=45, in=90] (0.7, 2) node[below right]{$\tilde t_4$} to[out=-90, in=135] (1, 1);
    \draw[arrows=->,semithick](0, -4)--(0,4) node[above right]{$y$};
    \draw[arrows=->,semithick](-4, 0)--(4, 0) node[above right]{$x$};
    \draw[-{To[length=3mm,width=2mm]}, ultra thick](0, 2)--(0, 1);
    \draw[-{To[length=3mm,width=2mm]}, ultra thick](0, -2)--(0, -1);
    \draw[-{To[length=3mm,width=2mm]}, ultra thick](0, 0)--(1.3, 0);
    \draw[-{To[length=3mm,width=2mm]}, ultra thick](0, 0)--(-1.3, 0);
  \end{tikzpicture}
  \caption{Domains where $\tilde t_i$ are defined.}
  \label{f:ti}
\end{figure}

For simplicity we will assume that the Hamiltonian $H$ is analytic. Then by~\cite{moser1956analytic}\footnote{
The result of~\cite{moser1956analytic} is for the case when $H$ periodically depends on the time, but one may check that when $H$ does not depend on the time the coordinate change constructed in~\cite{moser1956analytic} also does not depend on the time. The dependence on the parameter is also absent in~\cite{moser1956analytic}, but the proof may be easily adapted for the parametric case.}
one can find a new coordinate system $x, y$ in the neighborhood of the saddle $C$ such that this coordinate change is analytic and volume preserving, and the unperturbed system in the new coordinates is determined by a Hamiltonian $H_{x, y} = H_{x, y}(xy, z)$
with $H_{x, y}(C) = 0$ for all $z$ (we may subtract $H_{x, y}(C)$ from $H_{x, y}$ if this does not hold).
Let $\tilde{h} = xy$, $\tilde w = w = z$, denote $a(\tilde h, \tilde w) = \dv{H_{x, y}}{\tilde{h}}$ (we have $a \ne 0$).
Then in the new chart the unperturbed system rewrites as
\begin{align}
\begin{split} \label{e:ode-xy}
  \dot x = a(\tilde h, \tilde w) x, \;
  \dot y = -a(\tilde h, \tilde w) y.
\end{split}
\end{align}
Note that $\tilde h$ is a first integral of this system.
Also note that $\tilde h$ is a smooth function of $h, z$, as one can find $\tilde h$ from the equality $H_{x, y}(\tilde h, z) = h$. This also means that $\tilde h$ is defined on the whole phase space, even far from $C$.
We also have for any fixed value of $z$
\begin{equation} \label{e:lim-tilde-h-h}
  \lim_{h \to 0} \frac {h} {\tilde h(h, z)} = a(0, z).
\end{equation}
We will assume that the coordinates $x, y$ are as drawn in Figure~\ref{f:ti}, else we can rotate this coordinate system by $\pi m / 2$. Then, as $h>0$ for $\tilde h = xy > 0$, we have $a > 0$.
Rescaling $p$, $q$, $x$ and $y$ if needed, we may assume that
the neighborhood of $C$ where the new coordinates are defined contains the square
$\mathcal S = \{ x, y : -1 \le x, y \le 1 \}$ for all $z$.
{\color{myblue}
We also need the image of the bisector of the angle between separatrices in the $p, q$ coordinate tangent to the transversal $\varphi=0$ (cf. Section~\ref{s:assumptions}) under the map $p, q \mapsto x, y$ to be tangent to the line $x=y$ for all $z$. This can be achieved by a $z$-dependent area-preserving rescaling of $x$ and $y$.
}
\noindent The diagonals $x = \pm y$ split $\mathcal S$ into four triangles adjacent to each of its sides. In each such triangle let us introduce the time $\tilde t_i$ (it can be positive or negative) that passes after the trajectory of the unperturbed system intersects the adjacent side of $\mathcal S$.
The time $\tilde t_i$ can also be continued outside the square to the neighborhood of the separatrix crossing the transversal $\tilde t_i = 0$ (it is a side of $\mathcal S$). Domains where each $\tilde t_i$ is defined are drawn in figure~\ref{f:ti}.
Note that the coordinate systems $\tilde h, \tilde w, \tilde t_i$ cover the whole phase space (we only consider non-negative values of $h$ close to zero here).

We will assume that $\varphi=0$ corresponds to the transversal $\Gamma$ given by $x = y \ge 0$. Note that here we consider the angle coordinate in the domain $G_3$, for the domains $G_1$ and $G_2$ the transversal $\Gamma$ would be given by $\{ x = \pm y \} \cap G_i$.

\subsection{Estimates on how $q, p$ depend on \texorpdfstring{$\tilde h, \; \tilde w, \; \tilde t_i$}{tilde h, w, ti}} \label{s:est-tilde-t}
Outside of $\mathcal S$ each point of the phase space is covered by two coordinate systems $\tilde h, \tilde w, \tilde t_i$. For both of them the coordinate change $p, q, z \leftrightarrow \tilde h, \tilde w, \tilde t_i$ is defined and is smooth without singularities. So we only need to consider what happens inside $\mathcal S$.
For definiteness, let us restrict ourselves to the triangle $\{ 1 \ge x \ge y \ge 0 \}$.
For brevity we will write just $\tilde t$ for the coordinate $\tilde t_i$ defined in this triangle. This means that $\tilde t$ is the time after the trajectory intersects the line $x=1$.
Note that $\tilde t \le 0$ inside our triangle.
We have
\begin{align} \label{e:hx}
\begin{split}
  &x = e^{a(\tilde h, \tilde w) \tilde t}, \; y = \tilde h e^{-a(\tilde h, \tilde w) \tilde t}, \; z =\tilde w; \\
  &\tilde h = xy, \; \tilde w = z, \; \tilde t = \frac{\ln x}{a(xy, z)};
\end{split}
\end{align}
\begin{align} \label{e:dxh}
\begin{split}
  \pdv{x}{\tilde h} &= \pdv{a}{\tilde h}\/(\tilde h, \tilde w) \tilde t x, \;
  \pdv{x}{\tilde w} = \pdv{a}{\tilde w}\/(\tilde h, \tilde w) \tilde t x, \;
  \pdv{x}{\tilde t} = a(\tilde h, \tilde w) x; \\
  \pdv{y}{\tilde h} &= -\pdv{a}{\tilde h}\/(\tilde h, \tilde w) \tilde t y + \frac{1}{x}, \;
  \pdv{y}{\tilde w} = -\pdv{a}{\tilde w}\/(\tilde h, \tilde w) \tilde t y, \;
  \pdv{y}{\tilde t} = -a(\tilde h, \tilde w) y; \\
  \pdv{z}{\tilde w} &= 1, \; \pdv{z}{\tilde h}, \pdv{z}{\tilde t} = 0.
\end{split}
\end{align}
Note that $\tilde t x = \tilde t e^{a \tilde t} = O(1)$, as $a \tilde t < 0$. We also have $x \ge \tilde h^{1/2}$, as $x \ge y$.
It follows that
\begin{equation} \label{e:dxh-est}
  \pdv{y}{\tilde h} = O(h^{-1/2}); \;
  \pdv{x}{\tilde h}, \; \pdv{y}{\tilde w}, \; \pdv{x}{\tilde w}
  = O(|\tilde t|+1)e^{-a \abs{\tilde t}}; \;
  \pdv{y}{\tilde t}, \; \pdv{x}{\tilde t} = O(1)e^{-a \abs{\tilde t}}; \;
  \pdv{z}{\tilde w} = 1; \;
  \pdv{z}{\tilde h}, \pdv{z}{\tilde t} = 0.
\end{equation}
Note that by~\eqref{e:lim-tilde-h-h} we may write   $O(h^k)$ instead of $O(\tilde h^k)$.
It also follows from~\eqref{e:dxh} that
\begin{align} \label{e:d2xh}
\begin{split}
  &\pdv[2]{y}{\tilde h} = - \frac{2 \tilde t}{x}\pdv{a}{\tilde h} + \dots = O(|\tilde t| + 1)h^{-1/2}; \;
  \pdv[2]{y}{\tilde w}{\tilde h} = O(|\tilde t| + 1)h^{-1/2}; \\
  &\pdv[2]{y}{\tilde t}{\tilde h} = O(h^{-1/2}); \;
  \pdv[2]{y}{\tilde w}{\tilde t}, \pdv[2]{x}{\tilde w}{\tilde t}, \pdv[2]{x}{\tilde h}{\tilde t}
  = O(\abs{\tilde t} + 1) e^{-a \abs{\tilde t}};
  \pdv[2]{y}{\tilde t}, \pdv[2]{x}{\tilde t}
  = O(e^{-a \abs{\tilde t}}); \\
  &\pdv[2]{y}{\tilde w},
  \pdv[2]{x}{\tilde w}{\tilde h}, \pdv[2]{x}{\tilde h}, \pdv[2]{x}{\tilde w}
  = O((|\tilde t|+1)^2)e^{-a \abs{\tilde t}}; \;
  \pdv{z}{\text{\textasteriskcentered \textasteriskcentered}} = 0.
\end{split}
\end{align}
Now let us return from $(x, y)$ to $(q, p)$. Let us consider a smooth function $\psi(x, y, z)$ without singularities, e.g. $\psi = q$ or $\psi = p$.
We will use the following formula ($a_i, \; b_i$ are some coordinate systems and $c$ is some function)
\begin{equation} \label{e:d2-composition}
  \pdv[2]{c}{a_i}{a_j} = \sum_{l} \pdv[2]{b_l}{a_i}{a_j} \pdv{c}{b_l}
  + \sum_{k, l} \pdv{b_l}{a_j} \pdv{b_k}{a_i} \pdv[2]{c}{b_k}{b_l}.
\end{equation}
We can estimate the derivatives of $\psi$, using the chain rule for the first derivatives and~\eqref{e:d2-composition} for the second derivatives, and~\eqref{e:dxh-est}, \eqref{e:d2xh}.
This gives us
\begin{align} \label{e:dqh}
\begin{split}
  &\pdv{\psi}{\tilde h} = O(h^{-1/2}); \;
  \pdv{\psi}{\tilde w} = O(1); \; \pdv{\psi}{\tilde t_i} = O(e^{-a \abs{\tilde t_i}}); \\
  &\pdv[2]{\psi}{\tilde h} = O(h^{-1}); \;
  \pdv[2]{\psi}{\tilde w}{\tilde h} = O(\abs{\tilde t} + 1) h^{-1/2}; \;
  \pdv[2]{\psi}{\tilde t}{\tilde h} = O(h^{-1/2}); \\
  &\pdv[2]{\psi}{w} = O(1); \;
  \pdv{\psi}{tw} = O(\abs{\tilde t} + 1) e^{-a \abs{\tilde t}}; \;
  \pdv[2]{\psi}{t} = O(e^{-a \abs{\tilde t}}).
\end{split}
\end{align}
Outside of $\mathcal S$ we can take as $\tilde t$ any of the two coordinates $\tilde t_i$ defined near each separatrix, we have $\abs{\tilde t} + 1 \sim 1$. These estimates are valid everywhere: we obtained them in a part of $\mathcal S$, in other parts of $\mathcal S$ they can be obtained similarly, and outside of $S$ we even have $O(1)$ on all right hand sides as the considered coordinate change is smooth.

Let us also consider a function $\psi_0$ with $\psi_0(C) = 0$ (e.g. $\psi_0 = f_h$). As $C$ corresponds to $x=y=0$, the functions $\psi_0, \pdv{\psi_0}{z}, \pdv[2]{\psi_0}{z}$ (here the derivatives are taken for fixed $x, y$) all vanish at $C$ and so are $O(e^{-a\abs{\tilde t}})$. Some of the estimates above turn out to be better for $\psi_0$:
\begin{align} \label{e:est-psi0-tilde}
\begin{split}
  \psi_0 = O(e^{-a\abs{\tilde t}}); \;
  \pdv{\psi_0}{\tilde w} = O(e^{-a|\tilde t|})(\abs{\tilde t} + 1); \;
  \pdv[2]{\psi_0}{\tilde w} = O(e^{-a|\tilde t|})(\abs{\tilde t} + 1)^2.
\end{split}
\end{align}

{ \color{myblue}
\begin{remark} \label{r:O-star}
  Now we can precisely define the notation $O_*$ from Table~\ref{t:est}.
  We write $g = O_*(h^\alpha \ln^\beta h)$ if $g = O(h^\alpha \ln^\beta h) e^{-a\abs{\tilde t}} (|\tilde t| + 1)^\gamma$ for some $\gamma$, where $\tilde t$ is one of the coordinates $\tilde t_i$. At each point one or two coordinates $\tilde t_i$ are defined. If there are two, they are both $O(1)$, so we may choose any of them as $\tilde t$.
\end{remark}
}

\subsection{Estimates on how \texorpdfstring{$\tilde h, \; \tilde t_i$}{tilde h, ti} depend on \texorpdfstring{$h, \; w, \; \varphi$}{h, w, phi}} \label{s:est-T}
First, recall that $\tilde h$ is an analytic function of $h, w$. As $\tilde h(0, w) = 0$, all summands in the series for $\tilde h$ contain $h$.
Hence, we can write $\tilde h = h \tilde h_0(h, w)$, where $h_0(h, w)$ is also analytic.
From this we have
\begin{equation} \label{e:d-tilde-h-w}
  \pdv{\tilde h}{w}, \pdv[2]{\tilde h}{w} = O(h).
\end{equation}

Denote by $S(h, w)$ the time that the solution of the unperturbed system with given $h, w$ takes to get from the diagonal of the square $\mathcal S$ to its side. Then the total time spent inside $\mathcal S$ during each period is $4 S$. From~\eqref{e:hx} we have $S = - \frac{\ln \tilde h}{2 a(\tilde h, \tilde w)}$. Hence, by~\eqref{e:d-tilde-h-w} we have
\begin{align} \label{e:est-S}
\begin{split}
  S &= O(\ln h), \; \pdv{S}{h} = O(h^{-1}), \; \pdv{S}{w} = O(\ln h), \\
  \pdv[2]{S}{h} &= O(h^{-2}), \; \pdv[2]{S}{h}{w} = O(h^{-1}), \; \pdv[2]{S}{w} = O(\ln h).
\end{split}
\end{align}
Denote by $T_{reg, 1}(h, w)$ and $T_{reg, 2}(h, w)$ the times that the solution of the unperturbed system spends outside $\mathcal S$ near each of the separatrix loops during each period. These are smooth functions of $h, w$.
Then
\begin{equation} \label{e:T-S}
  T = 4S + T_{reg, 1} + T_{reg, 2}.
\end{equation}
From~\eqref{e:est-S} we get the estimates on $T$, $\omega$ from Table~\ref{t:est}.

Let us recall that for the unperturbed system we denote by $t = T \varphi / (2\pi)$ the time passed after crossing the transversal $\varphi = 0$ given by $x=y>0$. For each $\tilde t_i$ we have $\tilde t_i = t - t_{0, i}$, where $t_{0, i}$ is the value of $t$ corresponding to $\tilde t_i = 0$. We have (see Figure~\ref{f:ti})
\begin{equation} \label{e:t-0-i}
  t_{0, i} = k S + k_1 T_{reg, 1} + k_2 T_{reg, 2} \text{ with } k \in \{ 1, 3 \}; \; k_1, k_2 \in \{ 0, 1 \}.
\end{equation}
Hence, we have
\begin{equation}
  \tilde t_i = \Big(
    4S + T_{reg, 1} + T_{reg, 2}
  \Big)
  \frac{\varphi}{2\pi}
  - k S - k_1 T_{reg, 1} - k_2 T_{reg, 2}.
\end{equation}
This may also be rewritten as
\begin{equation}
  \tilde t_i = S(h, w)
  \Big(\frac{2\varphi}{\pi} - k \Big)
  + T_{reg}(h, w, \varphi),
\end{equation}
where $T_{reg}$ has no singularities and $\frac{2\varphi}{\pi} - k = O(\ln^{-1} h) (\abs{\tilde t_i} + 1)$.

From these formulas,~\eqref{e:est-S}, smooth dependence of $\tilde h$ on $h, w$ and~\eqref{e:d-tilde-h-w} we get
\begin{align} \label{e:dtphi}
\begin{split}
  &
  \pdv{\tilde t_i}{h} = O(h^{-1} \ln^{-1} h)  (\abs{\tilde t_i} + 1); \;
  \pdv{\tilde t_i}{w} = O(\abs{\tilde t_i} + 1); \; \pdv{\tilde t_i}{\varphi} = O(\ln h); \\
  &\pdv{\tilde h}{h} = O(1); \;
  \pdv{\tilde w}{w} = 1; \;
  \pdv{\tilde h}{w} = O(h); \;
  \pdv{\tilde h}{\varphi}, \pdv{\tilde w}{\varphi}, \pdv{\tilde w}{h} = 0; \\
  &\pdv[2]{\tilde t_i}{h} = O(h^{-2} \ln^{-1} h)  (\abs{\tilde t_i} + 1); \;
  \pdv[2]{\tilde t_i}{h}{w} = O(h^{-1} \ln^{-1} h)  (\abs{\tilde t_i} + 1); \; \pdv[2]{\tilde t_i}{h}{\varphi} = O(h^{-1}); \\
  &\pdv[2]{\tilde t_i}{w} = O(|\tilde t_i| + 1); \; \pdv[2]{\tilde t_i}{w}{\varphi} = O(\ln h); \;
  \pdv[2]{\tilde t_i}{\varphi} = 0; \\
  &\pdv[2]{\tilde h}{h}, \pdv[2]{\tilde h}{h}{w} = O(1); \;
  \pdv[2]{\tilde h}{w} = O(h); \\
  &\pdv{\tilde h}{\text{\textasteriskcentered}}{\varphi} = 0; \;
  \pdv[2]{\tilde w}{\text{\textasteriskcentered}}{\text{\textasteriskcentered}} = 0
  \text{\; for $\text{\textasteriskcentered} = h, w, \varphi$}.
\end{split}
\end{align}

\subsection{Estimates on how $q, p$ depend on  \texorpdfstring{$h, w, \varphi$}{h, w, phi}}
As above, let $\psi(x, y, z)$ be a smooth function without singularities, e.g. $\psi = q$ or $\psi = p$. Applying to~\eqref{e:dqh} and~\eqref{e:dtphi} the chain rule for first derivatives and formula~\eqref{e:d2-composition} for second derivatives, we get the following estimates (here $\tilde t$ is one of the coordinates $\tilde t_i$ as in Section~\ref{s:est-tilde-t}):
\begin{align} \label{e:est-psi}
\begin{split}
  &\pdv{\psi}{h} = O(h^{-1}\ln^{-1}h) e^{-a|\tilde t|}(\abs{\tilde t} + 1); \;
  \pdv{\psi}{w} = O(1); \; \pdv{\psi}{\varphi} = O(\ln h) e^{-a|\tilde t|}; \\
  &\pdv[2]{\psi}{h} = O(h^{-2}\ln^{-1}h) e^{-a|\tilde t|}(\abs{\tilde t} + 1); \;
  \pdv[2]{\psi}{h}{w} = O(h^{-1}\ln^{-1}h) e^{-a|\tilde t|}(\abs{\tilde t} + 1)^2; \\
  &\pdv[2]{\psi}{h}{\varphi} = O(h^{-1}) e^{-a|\tilde t|}(\abs{\tilde t} + 1); \;
  \pdv[2]{\psi}{w} = O(1); \\
  &\pdv[2]{\psi}{w}{\varphi} = O(\ln h) e^{-a|\tilde t|}(\abs{\tilde t} + 1); \; \pdv[2]{\psi}{\varphi} = O(\ln^2 h) e^{-a|\tilde t|}.
\end{split}
\end{align}
Let us also note that for a function $\psi_0$ with $\psi_0(C) = 0$ we can use~\eqref{e:est-psi0-tilde} and some of the estimates above turn out to be better:
\begin{align} \label{e:est-psi0}
\begin{split}
  \psi_0 = O(e^{-a\abs{\tilde t}}); \;
  \pdv{\psi_0}{w} = O(e^{-a|\tilde t|})(\abs{\tilde t} + 1); \;
  \pdv[2]{\psi_0}{w} = O(e^{-a|\tilde t|})(\abs{\tilde t} + 1)^2.
\end{split}
\end{align}
Finally, as
$\pdv{x}{w}, \pdv{y}{w} = O(|\tilde t|+1)e^{-a|\tilde t|}$
by~\eqref{e:est-psi0}, we have
\begin{equation} \label{e:est-dw-dz}
  \pdv{\psi}{w} = \pdv{\psi}{z} + O(|\tilde t|+1)e^{-a|\tilde t|}.
\end{equation}

\subsection{Estimates on $f$} \label{s:est-f}
Here we obtain the estimates on $f_h$, $f_{w_i}$ and $f_\varphi$ from Table~\ref{t:est}. The estimates on $f_{w_i}$ together with its derivatives follow from~\eqref{e:est-psi} as $f_{w_i} = f_{z_i}$ is smooth without singularities.
The estimates on $f_h$ follow from~\eqref{e:est-psi} and ~\eqref{e:est-psi0}, as $f_h = f_q \pdv{h}{q} + f_p \pdv{h}{p} + f_z \pdv{h}{z}$ is smooth without singularities and $f_h(C) = 0$ (as by~\cite[Lemma 2.1]{neishtadt17} we have $\pdv{h}{p}\/(C), \pdv{h}{q}\/(C), \pdv{h}{z}\/(C) = 0$).

Let us estimate $f_\varphi(h, w, 0)$. Recall that $t$ is the time passed after the solution of the unperturbed system crosses the transversal $x = y > 0$.
For $x, y > 0$ we have $t = \frac{1}{2 a(h, w)} (\ln x - \ln y)$, this is obtained by solving~\eqref{e:ode-xy} with initial conditions $x = y = \tilde h^{1/2}$ for $t=0$.For $\varphi=0$ (and therefore $t=0$, $x=y=\tilde h^{1/2}$) we have
\[
  \pdv{\varphi}{x} = \pdv{}{x}\/(\omega t) = \omega \pdv{t}{x}
  = \frac{\omega}{2a(h, w)x}, \qquad
  \pdv{\varphi}{y}
  = - \frac{\omega}{2a(h, w)y}, \qquad
  \pdv{\varphi}{z} = 0, \qquad
\]
\[
  f_\varphi(h, w, 0) = f_x \pdv{\varphi}{x} + f_y \pdv{\varphi}{y}
  = \frac{\omega}{2a(h, w)}(x^{-1} f_x - y^{-1} f_y)
  = \frac{\omega \tilde h^{-1/2}}{2a(h, w)}(f_x - f_y).
\]
Here $f_x, f_y$ are the components of the vector field $f$ written in the $x, y$ chart, they are $O(1)$. Hence, $f_\varphi(h, w, 0) = O(h^{-1/2} \ln^{-1} h)$. We can apply~\eqref{e:est-psi} to $\psi = \frac{f_x - f_y}{2 a}$, together with~\eqref{e:d-tilde-h-w} this gives
\begin{equation}
  \pdv{f_\varphi(h, w, 0)}{h} = O(h^{-3/2}\ln^{-1} h), \qquad
  \norm{\pdv{f_\varphi(h, w, 0)}{w}} = O(h^{-1/2}\ln^{-1} h).
\end{equation}

Denote $g(h, w) = f_\varphi(h, w, \varphi_g)$, where $\varphi_g \approx \pi$ corresponds to $x = y < 0$.
As $\varphi_g$ corresponds to $t = 2S + T_{reg, 1}$, we have $\varphi_g = \pi + 0.5 \omega (T_{reg, 1} - T_{reg, 2})$.
We can write $g = g_0 + g_1$, where $g_0$ is computed as if $x = y < 0$ corresponds to $\varphi = 0$ and
$g_1 = \pdv{\varphi_g}{h} f_h + \pdv{\varphi_g}{w} f_w$.
As for $x=y$ we have $e^{-a\abs{\tilde t}} = O(h^{1/2})$ and $f_h = O(h^{1/2}), \pdv{f_h}{w} = O(h^{1/2} \ln h), \pdv{f_h}{h} = O(h^{-1/2})$ by~\eqref{e:est-psi} and~\eqref{e:est-psi0}, we have $g_1 = O(h^{-1/2} \ln^{-2} h)$
and $\pdv{g_1}{h} = O(h^{-3/2} \ln^{-2} h)$, $\pdv{g_1}{w} = O(h^{-1/2} \ln^{-1} h)$.
For $g_0$ and its derivatives we can use the estimates for $f_\varphi(h, w, 0)$ proved above.
Hence, the estimates for $f_\varphi(h, w, 0)$ proved above also hold for $g$, $\pdv{g}{h}$ and $\pdv{g}{w}$.

\begin{lemma} \label{l:int-star}
Suppose that for a function $\alpha(h, w, \varphi)$ we have the estimate $\alpha(h, w, \varphi) = O_*(\tilde \alpha)$ with $\tilde \alpha = \tilde \alpha(h)$.
Then for any $\varphi_0, \varphi_1 \in [0, 2\pi]$ we have
$\int_{\varphi_0}^{\varphi_1} \alpha d\varphi = O(\tilde \alpha \ln^{-1} h)$.
\end{lemma}
\begin{proof}
  It is enough to show that $\int_0^T \abs{\alpha} dt = O(\tilde \alpha)$. This integral $dt$ can be splitted into four integrals $d \tilde t_i$, and each of them is $O(1)$ as the estimate on $\alpha$ contains a term that decays exponentially with the growth of $\tilde t_i$.
\end{proof}

From~\eqref{e:trace} and the estimates on $f_h$ and $f_w$ from Table~\ref{t:est} we have $\pdv{f_\varphi}{\varphi} = O_*(h^{-1}\ln^{-1} h)$.
Note that the estimates for $\Div(f)$ are given by~\eqref{e:est-psi}, as this function is smooth.
For given $h, w$ denote by $\varphi_*(\varphi)$ the angle corresponding to the "nearest" intersection of the solution with given $h$ and the line $x = y$. We have $\varphi_* = 0$ for $\varphi < \pi/2$ or $\varphi > 3\pi/2$ and $\varphi_*(\varphi) \approx \pi$ for $\pi/2 \le \varphi \le 3\pi/2$.
Arguing as in the proof of Lemma~\ref{l:int-star}, from
$\pdv{f_\varphi}{\varphi} = O_*(h^{-1}\ln^{-1} h)$
we can obtain
$f_\varphi(\varphi) - f_\varphi(\varphi_*) = O_*(h^{-1}\ln^{-2} h)$.
As for $x=y$ we have $e^{-a\abs{\tilde t}} \sim h^{1/2}$ and
$f_\varphi(\varphi_*) = O(h^{-1/2} \ln^{-1} h) = O(h^{-1} \ln^{-2} h) e^{-a\abs{\tilde t}} (|\tilde t| + 1)$,
we have $f_\varphi = O_*(h^{-1}\ln^{-2} h)$.

Let us apply $\pdv{}{h}$ to \eqref{e:trace}, this gives $\pdv[2]{f_\varphi}{\varphi}{h} = O_*(h^{-2}\ln^{-1} h)$. Arguing as above, we get $\pdv{f_\varphi}{h} = O_*(h^{-2}\ln^{-2} h)$. The estimate $\norm{\pdv{f_\varphi}{w}} = O_*(h^{-1}\ln^{-2} h)$ is obtained in the same way.

{ \color{myblue}
\subsection{Estimates on the choice of transversal}

\begin{lemma} \label{l:d-delta-phi}
  Take a transversal $\Gamma'$ to the union of separatrices passing through $C$ for all $z$.
  Denote by $\Delta \varphi(h, w)$ the value of the angle variable $\varphi$ (counted from $\Gamma$) on $\Gamma'$.
  Then
\begin{equation} \label{e:d-delta-phi}
  \pdv{\Delta \varphi}{h} = O(h^{-1} \ln^{-2} h), \qquad
  \pdv{\Delta \varphi}{w} = O(\ln^{-1} h).
\end{equation}
\end{lemma}
\begin{proof}
The transversal $\varphi = 0$ is given by $x=y$ in the $x, y$ chart. In the $x, y, z$ chart $\Gamma'$ can be defined by a function $y = G(x, z)$ with $G(0, z) = 0$. Denote $\pdv{G}{x}\/(0, z) = b(z)$,
\begin{equation}
  G(x, z) = b(z) x + O(x^2).
\end{equation}

Denote by $x_{\Gamma'}(\tilde h, \tilde w)$ the value of $x$ on the transversal $\Gamma'$ determined by given $\tilde h = xy$ and $\tilde w = z$. We have
\begin{equation}
  x_{\Gamma'} G(x_{\Gamma'}, \tilde w) = \tilde h.
\end{equation}
This gives
\begin{equation}
  x_{\Gamma'} = \sqrt{\tilde h / b(\tilde w)} + O(\tilde h), \qquad
  \pdv{x_{\Gamma'}}{\tilde h} = \frac{1}{2 \sqrt{b(\tilde w) \tilde h}} + O(1), \qquad
  \pdv{x_{\Gamma'}}{\tilde w} = O(\sqrt{\tilde h}).
\end{equation}
Denote by $\Delta t(\tilde h, \tilde w)$ the difference in time (for unperturbed system) between $\Gamma'$ and $\Gamma$, $\Delta \varphi = \omega \Delta t$.
When $\varphi=0$, we have $x=y=\sqrt{\tilde h}$.
By the formula~\eqref{e:hx} for $\tilde t$ we have
\begin{equation}
  \Delta t = a(\tilde h, \tilde w)^{-1} (\ln x_{\Gamma'} - 0.5 \ln \tilde h)
  = a(\tilde h, \tilde w)^{-1} \ln \frac{x_{\Gamma'}}{\sqrt{\tilde h}}.
\end{equation}
This gives
\begin{equation}
  \Delta t = O(1), \qquad
  \pdv{\Delta t}{\tilde h} = O(\tilde h^{-1/2}), \qquad
  \pdv{\Delta t}{\tilde w} = O(1).
\end{equation}
As $\pdv{\tilde h}{w} = O(h)$ and $h \sim \tilde h$, we have
\begin{equation}
  \Delta t = O(1), \qquad
  \pdv{\Delta t}{h} = O(h^{-1/2}), \qquad
  \pdv{\Delta t}{w} = O(1).
\end{equation}
Using $\Delta \varphi = \omega \Delta t$, we get~\eqref{e:d-delta-phi}.

\begin{remark} Denote by $\Delta t(h, w)$ the difference in time (for unperturbed system) between $\Gamma'$ and $\Gamma$. Suppose $\Gamma'$ is tangent to the bisector of the angle between separatrices. Then, arguing as in the proof of Lemma~\ref{l:d-delta-phi} (with $b=1$), we get
  \begin{equation} \label{e:delta-t-tangent}
    \Delta t = O(\sqrt{h}).
  \end{equation}
\end{remark}

\end{proof}
\begin{lemma} \label{l:d-delta-phi-2}
  Take a transversal $\Gamma'$ to one of the separatices that is far from $C$ for all $z$.
  Denote by $\Delta \varphi(h, w)$ the value of the angle variable $\varphi$ (counted from $\Gamma$) on $\Gamma'$.
  Then
\begin{equation} \label{e:d-delta-phi-2}
  \pdv{\Delta \varphi}{h} = O(h^{-1} \ln^{-2} h), \qquad
  \pdv{\Delta \varphi}{w} = O(\ln^{-1} h).
\end{equation}
\end{lemma}
\begin{proof}
  Suppose $\Gamma'$ is a transversal to the separatrix $l_2$, the proof is similar for $l_1$.
  We can write $\Delta \varphi = 2\pi \frac{S + b(h, w)}{4 S + d(h, w)}$, where $S$ is as in Section~\ref{s:est-T} and $b(h, w), d(h, w)$ are smooth functions; $4S + d = T$.
  We can rewrite $\Delta \varphi = \pi/2 + 2\pi \frac{b - d/4}{4 S + d} = \pi/2 + \omega (b - d/4)$.
  Product rule gives the required estimates.
\end{proof}

}

\section{Estimates for $\overline f_{h, 2}$ and $\overline f_{w, 2}$} \label{s:est-f-2}
For brevity, in this section we will write $u_h$ instead of $u_{h, 1}$, $\overline f_h$ instead of $\overline f_{h, 1}$, and so on.
\subsection{Expressions for $\overline f_{h, 2}$ and $\overline f_{w, 2}$}
\begin{lemma}
\begin{align}
\begin{split} \label{e:f-h-w-2}
  2\pi \overline f_{h, 2} &= \int_0^{2\pi} (\Div f - \sum_{w_i} \pdv{f_{w_i}}{w_i}) u_h + \sum_{w_i} \pdv{f_h}{w_i} u_{w_i} \; d\varphi
  \:-\:
  2\omega^{-1} \sum_{w_i} \pdv{\omega}{w_i} \int_0^{2\pi} f_h u_{w_i} d\varphi. \\
  2\pi \overline f_{a, 2}
  &=
  \omega \pdv{}{h} \int_0^T f_a u_h dt
  - 2 \omega^{-1} \sum_{w_i} \pdv{\omega}{w_i} \int_0^{2\pi} f_a u_{w_i} d\varphi \; + \\
  &+
  \int_0^{2\pi}
  \Div(f)u_a + \omega^{-1} \sum_{w_i} (\pdv{f_a}{w_i} u_{w_i} - \pdv{f_{w_i}}{w_i} u_a) d\varphi
  \text{ for } a = w_1, \dots, w_k.
\end{split}
\end{align}
\end{lemma}
\begin{proof}
Fix $a \in \{h, w_1, \dots, w_k\}$.
By~\eqref{e:f-h-2-init} we have
\begin{equation}
  2\pi \overline f_{a, 2} = \int_0^{2\pi} \pdv{f_a}{h} u_h + \pdv{f_a}{\varphi} u_\varphi + \sum_{w_i} \pdv{f_a}{w_i} u_{w_i} \; d\varphi.
\end{equation}
By Lemma~\ref{l:varchange-short} we have
$\pdv{u_b}{\varphi} = \frac{1}{\omega} (f_b - \overline f_b), \; b = h, w_1, \dots w_k$;
$\pdv{u_\varphi}{\varphi} = \frac{1}{\omega} (f_\varphi - \overline f_\varphi + \pdv{\omega}{h}u_h + \sum_{w_i} \pdv{\omega}{w_i}u_{w_i})$.
We also have $\int_0^{2\pi} u_b d\varphi = 0$, $b=h, \varphi, w_1, \dots, w_k$.
Hence, $\int_0^{2\pi} \overline f_b u_c d\varphi = 0$ for $b, c = h, \varphi, w_1, \dots, w_k$.
Integrating by parts, we have
\begin{align}
  &\int_0^{2\pi} \pdv{f_a}{\varphi} u_{\varphi} d\varphi
  = - \int_0^{2\pi} f_a \pdv{u_\varphi}{\varphi} d\varphi = \\
  &\qquad= - \omega^{-1} \int_0^{2\pi} f_a f_\varphi d\varphi
  + \frac{2\pi}{\omega} \overline f_\varphi \overline f_a
  -  \omega^{-1} \sum_{b=h, w_1, \dots w_k} \pdv{\omega}{b} \int_0^{2\pi} f_a u_b d\varphi.
\end{align}
Similarly, we have
\begin{equation}
  \int_0^{2\pi} \pdv{f_\varphi}{\varphi} u_a d\varphi
  = - \int_0^{2\pi} f_\varphi \pdv{u_a}{\varphi} d\varphi
  = - \omega^{-1} \int_0^{2\pi} f_a f_\varphi d\varphi
  + \frac{2\pi}{\omega} \overline f_\varphi \overline f_a.
\end{equation}
Hence, we have
\begin{equation}
  2\pi \overline f_{a, 2} = \int_0^{2\pi} \pdv{f_a}{h} u_h + \pdv{f_\varphi}{\varphi} u_a + \sum_{w_i} \pdv{f_a}{w_i} u_{w_i} \; d\varphi
  \: - \: \omega^{-1} \sum_{b=h, w_1, \dots w_k} \pdv{\omega}{b} \int_0^{2\pi} f_a u_b d\varphi.
\end{equation}
Expressing $\pdv{f_\varphi}{\varphi}$ through~\eqref{e:trace} gives
\begin{equation}
\int_0^{2\pi} \pdv{f_\varphi}{\varphi} u_a d\varphi =
\int_0^{2\pi} \bigg(\Div f - \pdv{f_h}{h} - \sum_{w_i} \pdv{f_{w_i}}{w_i}
- T^{-1}\Big(\sum_{b=h, w_1, \dots, w_k} \pdv{T}{b} f_b \Big)
\bigg)u_ad\varphi.
\end{equation}
Integrating by parts gives
\begin{equation} \label{e:fa-ub}
  \int_0^{2\pi} f_a u_b d\varphi = -\int_0^{2\pi} f_b u_a d\varphi, \qquad a=h, w_1, \dots, w_k, b=h, w_1, \dots, w_k.
\end{equation}
We also have
$T^{-1}\pdv{T}{b} + \omega^{-1} \pdv{\omega}{b} = \pdv{}{b}\/(\ln T + \ln \omega) = 0$ for $b=h, w_1, \dots, w_k$.
Thus,
\begin{equation}
\int_0^{2\pi} \pdv{f_\varphi}{\varphi} u_a d\varphi=
\int_0^{2\pi} \bigg(\Div f - \pdv{f_h}{h} - \sum_{w_i} \pdv{f_{w_i}}{w_i}
\bigg)u_a d\varphi
- \omega^{-1}\Big(\sum_{b=h, w_1, \dots, w_k} \pdv{\omega}{b} \int_0^{2\pi} f_a u_b d\varphi \Big)
\end{equation}
So, we have
\begin{align} \label{e:f-a-2}
\begin{split}
  2\pi \overline f_{a, 2} = &\int_0^{2\pi} \pdv{f_a}{h} u_h + \sum_{w_i} \pdv{f_a}{w_i} u_{w_i}
  + \big(\Div f - \pdv{f_h}{h} - \sum_{w_i} \pdv{f_{w_i}}{w_i}
  \big)u_a\; d\varphi \:- \\
  &-\: 2\omega^{-1} \sum_{b=h, w_1, \dots w_k} \pdv{\omega}{b} \int_0^{2\pi} f_a u_b d\varphi.
\end{split}
\end{align}

Note that for $b = h, w_1, \dots w_k$ we have
\begin{equation} \label{e:int-ua-fa}
  \int_0^{2\pi} f_b u_b d\varphi = \omega \int_0^{2\pi} u_b \pdv{u_b}{\varphi} d\varphi = 0.
\end{equation}
Therefore, for $a = h$ we have
\begin{equation}
  2\pi \overline f_{h, 2} = \int_0^{2\pi} \sum_{w_i} \pdv{f_h}{w_i} u_{w_i}
  + \big(\Div f - \sum_{w_i} \pdv{f_{w_i}}{w_i}
  \big)u_h\; d\varphi
  \:-\:
  2\omega^{-1} \sum_{w_i} \pdv{\omega}{w_i} \int_0^{2\pi} f_h u_{w_i} d\varphi,
\end{equation}
which is the first formula in~\eqref{e:f-h-w-2}.

Now we assume that $a = w_1, \dots, w_k$. We can compute
$
  \omega \pdv{}{h} \int_0^T f_a u_h dt =
  \omega \pdv{}{h} \bigg( \omega^{-1} \int_0^{2\pi} f_a u_h d\varphi \bigg) =
  - \omega^{-1} \pdv{\omega}{h} \int_0^{2\pi} f_a u_h d\varphi + \pdv{}{h} \int_0^{2\pi} f_a u_h d\varphi.
$
As
$\pdv{}{h} \int_0^{2\pi} f_a u_h d\varphi = \int_0^{2\pi} \pdv{f_a}{h} u_h d\varphi + \int_0^{2\pi} f_a \pdv{u_h}{h} d\varphi$,
we get
\begin{equation}
  \int_0^{2\pi} \pdv{f_a}{h} u_h d\varphi
  =
  \omega \pdv{}{h} \int_0^T f_a u_h dt
  -
  \int_0^{2\pi} f_a \pdv{u_h}{h}  d\varphi
  +
  \omega^{-1} \pdv{\omega}{h} \int_0^{2\pi} f_a u_h d\varphi.
\end{equation}
Plugging this into~\eqref{e:f-a-2} yields
\begin{align}
\begin{split}
  2\pi \overline f_{a, 2}
  =&
  \omega \pdv{}{h} \int_0^T f_a u_h dt
  +
  \int_0^{2\pi}
  \sum_{w_i} \pdv{f_a}{w_i} u_{w_i}
  + \big(\Div f - \sum_{w_i} \pdv{f_{w_i}}{w_i}
  \big)u_a
  - \big(\pdv{f_h}{h}u_a + \pdv{u_h}{h}f_a \big)
  \; d\varphi \:- \\
  &-\: 2\omega^{-1} \sum_{w_i} \pdv{\omega}{w_i} \int_0^{2\pi} f_a u_{w_i} d\varphi
  -
  \omega^{-1} \pdv{\omega}{h} \int_0^{2\pi} f_a u_h d\varphi.
\end{split}
\end{align}
As $\langle \pdv{u_h}{h} \rangle_\varphi = 0$, we have
\begin{align}
\begin{split}
  &\int_0^{2\pi} f_a \pdv{u_h}{h} d\varphi = \omega \int_0^{2\pi} \pdv{u_h}{h} \pdv{u_a}{\varphi} d\varphi
  = - \omega \int_0^{2\pi} u_a \pdv{}{h}\/\bigg( \omega^{-1} (f_h - \overline f_h) \bigg) d\varphi
  = \\
  &\qquad = - \int_0^{2\pi} u_a \pdv{f_h}{h} d\varphi
  + \omega^{-1} \pdv{\omega}{h} \int_0^{2\pi} u_a f_h d\varphi.
\end{split}
\end{align}
Hence,
\begin{align}
\begin{split}
  2\pi \overline f_{a, 2}
  =&
  \omega \pdv{}{h} \int_0^T f_a u_h dt
  + \int_0^{2\pi} \sum_{w_i} \pdv{f_a}{w_i} u_{w_i}
  + \big(\Div f - \sum_{w_i} \pdv{f_{w_i}}{w_i}
  \big)u_a
  \; d\varphi \:- \\
  &-\: 2\omega^{-1} \sum_{w_i} \pdv{\omega}{w_i} \int_0^{2\pi} f_a u_{w_i} d\varphi
  - \omega^{-1} \pdv{\omega}{h} \int_0^{2\pi} f_a u_h d\varphi
  - \omega^{-1} \pdv{\omega}{h} \int_0^{2\pi} u_a f_h  d\varphi.
\end{split}
\end{align}
By~\eqref{e:fa-ub} the last two terms cancel out, and we get the second formula in~\eqref{e:f-h-w-2}.
\end{proof}

\subsection{Estimate for $\overline f_{h, 2}$} \label{s:est-f-h-2}
\begin{lemma} \label{l:est-int-f-u}
  Let $\psi$ be either a smooth function $\psi(p, q, z)$ or the function
  $\pdv{f_b}{w_i}, \; b=h, w_1, \dots, w_k, \; i=1, \dots, k$
  and $u(h, w, \varphi)$ be a function with $u = O(1)$, $\langle u \rangle_\varphi = 0$ (e.g. $u_{b, 1}$ for $b=h, w_1, \dots, w_k$).
  Then
  \begin{equation}
    \int_0^{2\pi} \psi \cdot u \; d\varphi = O(\ln^{-1} h).
  \end{equation}
\end{lemma}
\begin{proof}
  First, let us prove that
  $\int_0^T |\psi - a| dt = O(1)$, where $a=\psi(C)$ for smooth $\psi$ and $a = \pdv{f_b}{z_i}\/(C)$ for $\psi=\pdv{f_b}{w_i}$.
  For smooth $\psi$ this follows from~\cite[Lemma 3.2]{neishtadt17}.
  For $\psi=\pdv{f_b}{w_i}$ by \eqref{e:est-dw-dz} we have
  $\psi - \pdv{f_b}{z_i} = O(|\tilde t_i|+1)e^{-a|\tilde t_i|}$ (here $\pdv{f_b}{z_i}$ is smooth). From this we have $\int |\psi - \pdv{f_b}{z_i}| dt = O(1)$, so the required statement follows from the smooth case with $\psi = \pdv{f_b}{z_i}$.

  As $u = O(1)$, the estimate above implies $\int_0^T (\psi - a)u \; dt = O(1)$.
  As $\langle u \rangle_\varphi = 0$, we have $\int_0^T a \cdot u \; dt = 0$ and $\int_0^T \psi \cdot u \; dt = O(1)$.
  Changing the variable, we obtain the required estimate.
\end{proof}

\begin{lemma}
  \begin{equation}
    \overline f_{h, 2} = O(\ln^{-1} h)
  \end{equation}
\end{lemma}
\begin{proof}
  By Table~\ref{t:est} we have $\omega^{-1} \pdv{\omega}{w_i} = O(1)$. Plugging this and the estimate of Lemma~\ref{l:est-int-f-u} in~\eqref{e:f-h-w-2} yields the required estimate.
\end{proof}

\subsection{Estimate for $\overline f_{w, 2}$} \label{s:est-f-w-2}
\begin{lemma} \label{l:d-int-h}
  For any smooth function $\psi(p, q, z)$ we have
  \begin{equation}
    \pdv{}{h}\/\Big(\int_{t=0}^{T} \psi u_h dt\Big) = O(h^{-1} \ln^{-2} h)
  \end{equation}
\end{lemma}
\begin{proof}
  We will assume $\psi(C) = 0$, as we can replace $\psi$ by $\psi - \psi(C)$ due to $\langle u_h \rangle_\varphi = 0$.
  We will use the integral expression for $u$ given by~\eqref{e:u-int-phi}:
  \begin{equation} \label{e:int-for-u}
    u_h(t_0) = \int_{0}^{T} \Big(\frac{t}{T} - \frac 1 2 \Big) f_h(t+t_0) dt.
  \end{equation}
  We have (in the formula below $t_1=t+t_0+kT$, where $k \in \mathbb Z$ is such that $t_1 \in [0, T)$; $\{ x \}$ denotes the fractional part of $x$, i.e. such number $y$ in $[0, 1)$ that $x - y \in \mathbb Z$)
  \begin{equation}
    \int_{t_0=0}^{T} \psi u_h dt_0
    = \int_{0}^{T}\int_{0}^{T} \Big(\frac{t}{T} - \frac 1 2 \Big) f_h(t+t_0) \psi(t_0) dt dt_0
    = \int_{0}^{T}\int_{0}^{T} \Big(\Big\{ \frac{t_1 - t_0}{T} \Big\} - \frac 1 2 \Big) f_h(t_1) \psi(t_0) dt_1 dt_0.
  \end{equation}
    We will use the following notation from sections~\ref{s:Moser}-\ref{s:est-T}: $x, y, \tilde t_i, \tilde h, S$.
  The phase space can be splitted by the lines $x = y$, $\tilde t_1 = - \tilde t_2$ and $\tilde t_3 = - \tilde t_4$ (see Figure~\ref{f:ti}) into four parts, such that each part is covered by one of the coordinates $\tilde t_1, \dots, \tilde t_4$. Let all possible values of $\tilde t_i$ in its part span the segment $[a_i, b_i]$. Note that for the second two lines the values of the coordinates $\tilde t_i$ defined there (and also the values of the corresponding $a_i$ or $b_i$) are smooth functions of $h, w$ without singularities. For example, for $\tilde t_1 = - \tilde t_2$ we have $\tilde t_1 = b_1 = T_{reg, 1}/2$ and $\tilde t_2 = a_2 = -T_{reg, 1}/2$, where $T_{reg, 1}(h, w)$ is the time between the points with $\tilde t_1 = 0$ and $\tilde t_2 = 0$.

   Denote $\lambda_{i, j}(t_0, t_1) = \Big\{ \frac{t(\tilde t_j = t_1) - t(\tilde t_i = t_0)}{T} \Big\} - \frac 1 2$. The integral above can be split into a sum of the following $16$ integrals for $i, j = 1, \dots, 4$:
  \begin{equation} \label{e:16-integrals}
    \int_{t_0 = a_i}^{b_i} \int_{t_1 = a_j}^{b_j}
    \lambda_{i, j}(t_0, t_1) f_h(\tilde t_j = t_1) \psi(\tilde t_i=t_0) dt_1 dt_0.
  \end{equation}
  As $\pdv{}{h} = \pdv{\tilde h}{h}\pdv{}{\tilde h} = O(1)\pdv{}{\tilde h}$, we will estimate the $\tilde h$-derivative of~\eqref{e:16-integrals} instead of its $h$-derivative.
  Let us first note that the discontinuity of $\lambda_{i, j}(t_0, t_1)$ corresponds to $\tilde t_i = t_0$ and $\tilde t_j = t_1$ giving the same point, so $i=j$ and $t_0=t_1$, and this discontinuity does not create additional terms in the $\tilde h$-derivative of~\eqref{e:16-integrals}.
  By~\eqref{e:t-0-i} we have
  $\lambda_{i, j}(t_0, t_1) + 0.5 = \Big\{ \frac{t_1 - t_0 + k_{i, j} S + T_{reg, i, j}}{T} \Big\}$ (here $k_{i,j} \in \mathbb Z$ and $T_{reg, i, j}$ is a smooth function of $h, w$),
  so $\pdv{\lambda_{i, j}(t_0, t_1)}{\tilde h} = O(h^{-1}\ln^{-2} h)(\abs{t_1 - t_0} + 1)$. As $f_h(C) = \psi(C) = 0$, by~\eqref{e:est-psi0} we have $f_h(\tilde t_s), \psi(\tilde t_s) = O(e^{-a|\tilde t_s|})$ for $s=1, \dots, 4$, this means
  $\int \int \abs{t_0 - t_1} f_h(t_1) \psi(t_0) dt_0 dt_1 = O(1)$
  and so
  $\int \int \pdv{\lambda_{i, j}}{\tilde h} f_h(t_1) \psi(t_0) dt_0 dt_1 = O(h^{-1} \ln^{-2} h)$.
  By~\eqref{e:dqh} we have
  $\int \int \lambda_{i, j} \pdv{f_h}{\tilde h} \psi dt_0 dt_1 = O(h^{-0.5} \ln^2 h)$
  and
  $\int \int \lambda_{i, j} f_h \pdv{\psi}{\tilde h} dt_0 dt_1 = O(h^{-0.5} \ln^2 h)$.
  The $\tilde h$-derivative of~\eqref{e:16-integrals} also has terms associated with the change of the domain of integration. There are four similar terms, let us consider just one of them:
  \begin{equation}
    \pdv{a_i}{\tilde h} \int_{t_1 = a_j}^{b_j}
    \lambda_{i, j}(a_i, t_1) f_h(\tilde t_j = t_1) \psi(\tilde t_i = a_i) dt_1.
  \end{equation}
  There are two cases. First, $a_i$ may correspond to $x=y=O(\sqrt{h})$, then $\pdv{a_i}{\tilde h} = O(h^{-1})$, $\psi = O(\sqrt{h})$,  and our term is $O(h^{-0.5} \ln h)$. Otherwise, we have $\pdv{a_i}{\tilde h} = O(1)$ and our term is $O(\ln h)$.

  Combining these estimates, we see that the $\tilde h$-derivative (and so also the $h$-derivative) of~\eqref{e:16-integrals} is $O(h^{-1} \ln^{-2} h)$. This proves the lemma.
\end{proof}

\begin{lemma}
  \begin{equation}
    \overline f_{w_i, 2} = O(h^{-1} \ln^{-3} h), \qquad i = 1, \dots, k.
  \end{equation}
\end{lemma}
\begin{proof}
  By Table~\ref{t:est} we have $\omega^{-1} \pdv{\omega}{w_i} = O(1)$. This and the estimate of Lemma~\ref{l:est-int-f-u} gives the estimate $O(\ln^{-1} h)$ for all terms of the expression~\eqref{e:f-h-w-2} for $\overline f_{w_i, 2}$ except the first one. The first term is estimated by Lemma~\ref{l:d-int-h}.
\end{proof}

\section{Estimates related to the averaging chart} \label{s:est-u}
In this section we prove the estimates from Table~\ref{t:est} for the functions $u_{k, i}$ and $\overline f_{k, i}$. The following lemma allows to mass-produce such estimates.
However, these estimates are not always good, so we will estimate some of these functions in a different way.
\begin{lemma} \label{l:est-u}
Given a function $Y(h, w, \varphi)$, let
\begin{equation} \label{e:Y-f}
  \overline f = \langle Y \rangle_\varphi
\end{equation}
and let the function $u$ be determined by the equation
$\omega \pdv{u}{\varphi} = Y - \overline f$
and the condition $\langle u \rangle_\varphi = 0$.
Denote $Y_T = T \cdot Y$.
Let $v = h, w, hh, hw, ww$ and $\pdv{}{v}$ denote the corresponding first or second derivative.
Then we can estimate the functions $\overline f$ and $u$ and their derivatives using estimates for $Y$ and $Y_T$ and their derivatives (these estimates are denoted by $\tilde Y, \tilde Y_T, \tilde Y_v, \tilde Y_{T, v}$ below, they depend only on $h$) in the following way:
\begin{enumerate}
  \item $\overline f = O(\tilde Y)$ for $Y = O(\tilde Y)$;
  $\overline f = O(\tilde Y \ln^{-1} h)$ for $Y = O_*(\tilde Y)$.

  \item $\pdv{\overline f}{v} = O(\tilde Y_v)$ for $\pdv{Y}{v} = O(\tilde Y_v)$;
  $\pdv{\overline f}{v} = O(\tilde Y_v \ln^{-1} h)$ for $\pdv{Y}{v} = O_*(\tilde Y_v)$.

  \item $\pdv{u}{\varphi} = O(\tilde Y_T)$ for $Y_T = O(\tilde Y_T)$ or $Y_T = O_*(\tilde Y_T)$.

  \item $\pdv[2]{u}{v}{\varphi} = O(\tilde Y_{T, v})$ for $\pdv{Y_T}{v} = O(\tilde Y_{T, v})$ or $\pdv{Y_T}{v} = O_*(\tilde Y_{T, v})$.

  \item $u = O(\tilde Y_T)$ for $Y_T = O(\tilde Y_T)$;
  $u = O(\tilde Y_T \ln^{-1} h)$ for $Y_T = O_*(\tilde Y_T)$.

  \item $\pdv{u}{v} = O(\tilde Y_{T, v})$ for $\pdv{Y_T}{v} = O(\tilde Y_{T, v})$;
  $\pdv{u}{v} = O(\tilde Y_{T, v} \ln^{-1} h)$ for $\pdv{Y_T}{v} = O_*(\tilde Y_{T, v})$.
\end{enumerate}
\end{lemma}
\begin{remark}
  As the maps $Y \mapsto u$ and $Y \mapsto \overline f$ are linear, for $Y = Y_1 + Y_2 = O(\tilde Y_1) + O_*(\tilde Y_2)$ we can estimate $u(Y)$ as $u(Y_1) + u(Y_2)$ and $\overline f(Y)$ as $\overline f(Y_1) + \overline f(Y_2)$.
\end{remark}
\begin{proof}
  Item~1 follows from~\eqref{e:Y-f} and Lemma~\ref{l:int-star}. Item~2 is proved in the same way, as $\pdv{}{v}$ commutes with averaging with respect to $\varphi$.

  We have
  \[
    2\pi \pdv{u}{\varphi} = Y_T - \langle Y_T \rangle_\varphi.
  \]
  This equation implies item~3. As $\langle u \rangle_\varphi = 0$, integrating this estimate for $\pdv{u}{\varphi}$ gives the first part of item~5. Together with Lemma~\ref{l:int-star} the equation above implies that for $Y_T = O_*(\tilde Y_T)$ we have $\langle Y_T \rangle_\varphi = O(\tilde Y_T \ln^{-1} h)$ and $u(\varphi_1) - u(\varphi_0) = O(\tilde Y_T \ln^{-1} h)$ for any $\varphi_0, \varphi_1$. This proves the second part of item~5.
  Items~4 and~6 are proved like items~3 and~5, we just need to take $\pdv{}{v}$ of the equation above.
\end{proof}
The functions $u_{*, *}$ and $\overline f_{*, *}$ are given by Lemma~\ref{l:varchange-short} and Lemma~\ref{l:varchange}.
Lemma~\ref{l:est-u} allows to obtain the estimates for $u_{h, 1}$, $u_{w, 1}$ and the derivatives of these functions. Let us note that $u_{w_i, 1}$ is determined by $Y = f_{w_i}$. However, for $Y = f_{w_i} - f_{w_i}(C)$ we get the same value of $u_{w_i, 1}$, but better estimates, as we may use~\eqref{e:est-psi0}.

The estimate for the functions $\overline f_{*, 1}$ and their derivatives are also obtained by the lemma above. Note that $\overline f_{\varphi, 1} = \langle f_\varphi \rangle_\varphi$, as $\langle u_{v, 1} \rangle_\varphi = 0$.
Using the estimates above, we can estimate $u_{\varphi, 1}$, $u_{h, 2}$, $u_{w, 2}$, $\overline f_{h, 2}$, $\overline f_{w_i, 2}$ and their derivatives by Lemma~\ref{l:est-u}.
However, for the functions  $\overline f_{h, 2}$ and $\overline f_{w, 2}$ themselves better estimates are obtained in sections~\ref{s:est-f-h-2} and~\ref{s:est-f-w-2}.

To estimate the functions $\overline f_{h, 3}$ and $\overline f_{\varphi, 2}$, we need to assume that
\begin{equation} \label{cond:f-h-3}
  h > C_h \varepsilon.
\end{equation}
The large enough constant $C_h > 0$ will be chosen below. It will be greater than the constant $C_{inv}$ from Lemma~\ref{l:invertible}.
By~\eqref{e:f3} we have the following system of equations (here $v = (h, w)$ and $A_\varphi$, $A_v$ denote the right hand sides of~\eqref{e:f3}):
\begin{align}
\begin{split}
   &(1 + \varepsilon \pdv{u_{\varphi, 1}}{\varphi}) \overline f_{\varphi, 2}
   + \varepsilon^2 \pdv{u_{\varphi, 1}}{v} \overline{f}_{v, 3} = A_\varphi, \\
   &(1 + \varepsilon \pdv{u_{v, 1, 2}}{v}) \overline{f}_{v, 3}
   + \pdv{u_{v, 1, 2}}{\varphi} \overline{f}_{\varphi, 2} = A_v.
\end{split}
\end{align}
From~\eqref{cond:f-h-3}, \eqref{e:varchange} and the estimates on $u_{h, 1}$ and $u_{h, 2}$ for large enough $C_h$ we have
\begin{equation}
  h \in [0.5 \overline h, 2 \overline h].
\end{equation}
This allows us to estimate the intermediate values from~\eqref{e:f3} as if they were at the point $\overline h$.
Using Table~\ref{t:est} and~\eqref{cond:f-h-3}, we have $A_\varphi = O(h^{-2} \ln^{-2} h), \; A_v = O_*(h^{-2}) + O(h^{-2} \ln^{-1} h)$.
We can substitute the expression for $\overline{f}_{v, 3}$ from the second equation into the first one.
This yields
\begin{align}
\begin{split}
  &\overline{f}_{\varphi, 2} \bigg(
    1 + \varepsilon \pdv{u_{\varphi, 1}}{\varphi}
    - \varepsilon^2 \pdv{u_{\varphi, 1}}{v}
    (1 + \varepsilon \pdv{u_{v, 1, 2}}{v})^{-1} \pdv{u_{v, 1, 2}}{\varphi}
  \bigg)
  = \\
  & \qquad =
  A_\varphi - \varepsilon^2 \pdv{u_{\varphi, 1}}{v}\/(1 + \varepsilon \pdv{u_{v, 1, 2}}{v})^{-1} A_v.
\end{split}
\end{align}
From~\eqref{cond:f-h-3} and Table~\ref{t:est} we see that for large enough $C_h$
\[
  \norm{\pdv{u_{v, 1, 2}}{\varphi}} = O(\ln h); \qquad
  \norm{\varepsilon \pdv{u_{v, 1, 2}}{v}}, \; \abs{\varepsilon \pdv{u_{\varphi, 1}}{\varphi}} < 0.1, \qquad
  \norm{\varepsilon^2 \pdv{u_{\varphi, 1}}{v}} = O(C_h^{-2} \ln^{-1} h).
\]
For large enough $C_h$ we have
$\abs{\varepsilon \pdv{u_{\varphi, 1}}{\varphi}
- \varepsilon^2 \pdv{u_{\varphi, 1}}{v}
(1 + \varepsilon \pdv{u_{v, 1, 2}}{v})^{-1} \pdv{u_{v, 1, 2}}{\varphi}} < 0.5$.
Hence, we have
$\overline{f}_{\varphi, 2} = O_*(h^{-2} \ln^{-1} h) + O(h^{-2} \ln^{-2} h)$
(let us note that for $h > C_h \varepsilon \abs{\ln \varepsilon}^{0.5}$ we have
$\norm{\varepsilon^2 \pdv{u_{\varphi, 1}}{v}} = O(\ln^{-2} h)$
and this yields slightly better estimate $\overline{f}_{\varphi, 2} = O(h^{-2} \ln^{-2} h)$). Then from the second equation we obtain $\overline f_{v, 3} = O_*(h^{-2}) + O(h^{-2} \ln^{-1} h)$.

\begin{lemma} \label{l:est-f-hat}
  The estimates for the functions
  $\overline f_{a, i}$ (for $a = h, w$, $i = 1, 2, 3$ and $a = \varphi$, $i = 1, 2$) and their derivatives stated in Table~\ref{t:est} also hold for the corresponding functions
  $\hat f_{a, i}$
   and their derivatives.
  Moreover, we have $|u_{a, 1}(h, w, \varphi, \varepsilon) - u^0_{a, 1}(h, w, \varphi)| = O(\varepsilon)$ for $a=h, w_1, \dots, w_k$.
\end{lemma}
\begin{proof}
Recall that the expressions $\overline f^0_{*, *}$ are computed by the same formulas as $\overline f_{*, *}$, with the perturbation $f$ replaced by $f^0$. This means that the estimates we have for $\overline f_{*, *}$ (they are valid for any smooth perturbation $f$) also hold for $\overline f^0_{*, *}$. By~\eqref{e:def-f-hat} we have $\hat f_{h, 1} = \overline f^0_{h, 1}$, $\hat f_{w, 1} = \overline f^0_{w, 1}$ and $\hat f_{\varphi, 1} = \overline f^0_{\varphi, 1}$, so for these expressions and their derivatives the lemma holds.

By~\eqref{e:def-f-hat} we also have $\hat f_{h, 2} = \overline f^0_{h, 2} + \langle f_h^1(h, w, \varphi) \rangle_\varphi$. Denote $\alpha = \langle f_h^1(h, w, \varphi) \rangle_\varphi$. Similarly to the estimate on $\overline f_{h, 1}$ above ($\alpha$ is computed exactly as $\overline f_{h, 1}$ if we start with $f^1$ instead of $f$), we have $\alpha = O(\ln^{-1} h)$, $\pdv{\alpha}{h} = O(h^{-1} \ln^{-2} h)$ and $\pdv{\alpha}{w} = O(1)$. Therefore, the estimates for $\overline f_{h, 2}$ and $\pdv{\overline f_{h, 2}}{h}$ from Table~\ref{t:est} also hold for $\hat f_{h, 2}$.
The estimates for $\hat f_{w, 2}$ are obtained in the same way.

We have $\hat f_{\varphi, 2} = \overline f_{\varphi, 2} + \varepsilon^{-1}(\overline f_{\varphi, 1} - \overline f^0_{\varphi, 1})$. Using~\eqref{e:f-series}, we get $\varepsilon^{-1}(\overline f_{\varphi, 1} - \overline f^0_{\varphi, 1}) = \langle f^1_\varphi + \varepsilon f^2_\varphi \rangle_\varphi$, where $f^i_\varphi$ is the $\varphi$-component of $f^i$ written in the energy-angle coordinates. As the estimate for $\overline f_{\varphi, 1} = \langle f \rangle_\varphi$ holds for any smooth $f$, we can plug in $f^1 + \varepsilon f^2$ instead of $f$ and get the estimate
$\langle f^1_\varphi + \varepsilon f^2_\varphi \rangle_\varphi = O(h^{-1}\ln^{-3} h)$.
As $f^2_p$ and $f^2_q$ are uniformly bounded by a constant independent of $\varepsilon$, one may check that this estimate is uniform in $\varepsilon$.
Therefore, the estimate for $\overline f_{\varphi, 2}$ also holds for $\hat f_{\varphi, 2}$.

Before estimating $\hat f_{h, 3}$, let us prove the second statement of the lemma.
For $b = h, w, \varphi$ the map $\mathcal U: f \to u_{b, 1}$ is linear by~\eqref{e:homological} and~\eqref{e:homological-rhs}. Hence, for
$\psi=f^1 + \varepsilon f^2$ and $u^\psi_{b, 1} = \mathcal U(\psi)$
we have
\[
  u_{b, 1}(h, w, \varphi, \varepsilon) = u^0_{b, 1}(h, w, \varphi) + \varepsilon u^\psi_{b, 1}(h, w, \varphi, \varepsilon).
\]
 As $\psi$ is smooth with respect to $p, q, z$ and uniformly bounded with respect to $\varepsilon$, the estimate $u_{a, 1} = O(1)$ for $a = h, w_1, \dots, w_k$ also holds for $u^\psi_{a, 1}$.

We have $\hat f_{h, 3} = \overline f_{h, 3} + \langle f^2_h \rangle_\varphi + \varepsilon^{-1}(\overline f_{h, 2} - \overline f^0_{h, 2})$. Clearly, $\langle f^2_h \rangle_\varphi = O(1)$. As we have $f_{h, 2} = \sum_a \langle \pdv{f_h}{a} u_a \rangle_\varphi, \; a =h, \varphi, w_1, \dots, w_k$ and the functions $u_{1, *}$ linearly depend on $f$, for $\psi=f^1 + \varepsilon f^2$ we can write
$\overline f_{h, 2} = \sum_a \langle (\pdv{f^0_h}{a} + \varepsilon \pdv{f^\psi_h}{a}) (u^0_a + \varepsilon u^\psi_a)  \rangle_\varphi$
and
\[
  \varepsilon^{-1}(\overline f_{h, 2} - \overline f^0_{h, 2})
  =
  \sum_a \Big\langle \pdv{f^0_h}{a} u^\psi_a +
  \pdv{f^\psi_h}{a} u^0_a +
  \varepsilon \pdv{f^\psi_h}{a} u^\psi_a \Big\rangle_\varphi, \qquad
  a = h, \varphi, w_1, \dots, w_k.
\]
Here the upper index $\psi$ means that the function is obtained using $\psi$ instead of $f$. The estimates on $\pdv{f_h}{a}, \; u_a$ from Table~\ref{t:est} are also valid for $\pdv{f^\psi_h}{a}, \; u^\psi_a$. Using these estimates, we obtain $\varepsilon^{-1}(\overline f_{h, 2} - \overline f^0_{h, 2}) = O(h^{-1})$, thus proving the estimate for $\hat f_{h, 3}$.
In a similar way we obtain $\varepsilon^{-1}(\overline f_{w_i, 2} - \overline f^0_{w_i, 2}) = O(h^{-1})$, thus proving the estimate for $\hat f_{w, 3}$.
\end{proof}

\begin{appendices}
{ \color{myblue}
\section{Estimates near the separatrices} \label{s:near-sep}
\begin{proof}[Proof of Lemma~\ref{l:delta-h-near-sep}]
By~\cite[Lemma~3.5, Corollary~3.4]{neishtadt17} there are $c_2 > 0$ and $c_3 > 0$ such that for $h_0 \in [c_3 \varepsilon, c_2]$ the orbit crosses the transversal again after time $O(\ln h_0)$ passes and we have $h = h_0 - \varepsilon \oint_{h=h_0} f_h dt + O(\varepsilon^2 h_0^{-1/2})$ (where integral is taken with $w=w_0$) at the transversal crossing.
As $\oint_{h=h_0} f_h dt = -\Theta_3 + O(h_0 \ln h_0)$ by~\cite[Corollary 3.1]{neishtadt17}, this gives the required estimate.

Hence, we consider only $h_0 < c_3 \varepsilon$.
By~\cite[Proposition~5.1]{neishtadt17} for $h_0 < c_3 \varepsilon$ the orbit of our point intersects the transversal $\varphi = 0$ once more
(the condition $\Theta_1, \Theta_2 > 0$ is used here).
Denote by $h_1$ the value of $h$ for this new intersection.
From~\cite[Proposition~5.1]{neishtadt17} during the whole wind between $h_0$ and $h_1$ we have $h = h_0 + O(\varepsilon)$.
Hence, we also have\footnote{
Trajectory might cross separatrices of unperturbed systems and leave $G_3$ into $G_2$ or $G_1$, this is why we write $|h|$ here.
} $|h| < d_1 \varepsilon$ for some $d_1$.

As $h_0 \sim \varepsilon$, \eqref{e:delta-h-near-sep} follows from the estimate
\begin{equation} \label{e:h1-h0}
  h_1 - h_0 = \varepsilon \Theta_3(w_0) + O(\varepsilon^{3/2}).
\end{equation}
We prove this estimate, arguing close to the proof of~\cite[Proposition~5.1]{neishtadt17}.
The rough idea is to use coordinate chart $x, y$ (cf. Section~\ref{s:Moser}) in the phase space $p, q$ of unperturbed system for fixed $z$ such that coordinate axes coincide with separatrices of unperturbed system. Then one can use one of the variables $x, y$ as an independent variable and write the change of $h$ as a certain integral that we then estimate.

We will assume
\begin{equation} \label{e:w-close-w0}
  \norm{w - w_0} < \varepsilon \ln^2 \varepsilon
\end{equation}
on the whole considered part of trajectory. If this condition fails at some point, we only consider the part of trajectory before this moment (but we will see that this condition actually holds until the next crossing of the transversal $\varphi = 0$).

Recall that $-\Theta_3$ is the time integral of $f_h$ along the separatrices.
Let us split the trajectory into four segments near the saddle (separated by the boundary of $\mathcal S$ and the line $x = y$, we use the notation from Section~\ref{s:Moser}) and two segments near each separatix.
To prove the lemma, we compare the change of $h$ for each segment to the integral of $f_h$ along the nearby part of one of the separatrices.
Denote by $\mathcal B \subset \mathcal S$ the domain given by $0 < x < y$.
We will only prove this estimate for the part of trajectory near the saddle that lies in $\mathcal B$ (i.e. $\varphi$ starting from $3\pi/2 + o(1)$ and approaching $2\pi$), as for other parts near the saddle the proof is similar and for the two parts far from the saddle the proof is much easier. Denote the value of $h$ when entering $\mathcal B$ (i.e. when $y=1$) by $h_2$. We will prove
\begin{equation} \label{e:h2-h1}
  h_2 - h_1 = \varepsilon \Theta_3' + O(\varepsilon^{3/2}),
\end{equation}
where
$\Theta'_3 = -\int_{x=0, y \in [0, 1]} f_h(x, y, z=w_0) dt$
is the integral of $f_h$ along the part of separatices of unperturbed system in $\mathcal B$. Taking sum of these estimates over all six segments of the trajectory will give~\eqref{e:h1-h0}.

We will estimate the change of $\tilde h = xy$ instead of the change of $h$.
As $|h| < d_1 \varepsilon$, we have $\pdv{h}{\tilde h} = a(\tilde h, w) = a(0, w) + O(\varepsilon)$ (the notation $a$ was introduced in~\eqref{e:ode-xy}). We have (here and thereafter $\pdv{}{\alpha} \big|_\beta$ denotes that partial derivative is taken for fixed $\beta$) $\pdv{h}{w}\big|_{\tilde h} = O(h) = O(\varepsilon)$, as follows from
\begin{equation}
  \pdv{\tilde h}{w}\Big|_{\tilde h} = 0 =
  \pdv{\tilde h}{w}\Big|_{h} + \pdv{\tilde h}{h}\Big|_{w} \pdv{h}{w}\Big|_{\tilde h},
\end{equation}
\eqref{e:d-tilde-h-w} and $\pdv{\tilde h}{h} \sim 1$. Thus
\begin{equation} \label{e:dh-dh-tilde}
  h_2 - h_1 = (a(0, w_0) + O(\varepsilon))(\tilde h_2 - \tilde h_1) + O(\varepsilon^2 \ln^2 \varepsilon).
\end{equation}
Set $f_{\tilde h} = x f_y + y f_x$, then $f_h = a(\tilde h, w) f_{\tilde h} + \pdv{h}{w} f_z$ and on the separatrices $f_h = a(0, w) f_{\tilde h}$.
Set $\Theta' = -\int_{x=0, y \in [0, 1]} f_{\tilde h}(x, y, z=w_0) dt$, we have $\Theta_3' = a(0, w_0) \Theta'$.
We will prove that
\begin{equation} \label{e:delta-h-last-wind}
  \tilde h_2 - \tilde h_1 = \varepsilon \Theta' + O(\varepsilon^{3/2}),
\end{equation}
it will imply~\eqref{e:h2-h1} by~\eqref{e:dh-dh-tilde}.

The perturbed system is written as follows in the coordinates $x, y, z$:
\begin{align}
\begin{split}
  \dot x = a(\tilde h, w) x + \varepsilon f_x(x, y, z), \qquad
  \dot y = - a(\tilde h, w) y + \varepsilon f_y(x, y, z), \qquad
  \dot z = \varepsilon f_z(x, y, z).
\end{split}
\end{align}
Recall that $a(0, w) > 0$.
We have $\dot{\tilde h} = \varepsilon f_{\tilde h}$.
We take the solution with initial data $y_2 = 1, x_2 = \tilde h_2 > \Theta' + d_2 \varepsilon^{3/2}$ and $z_2$ with $\norm{z_2 - w_0} < \varepsilon \ln^2 \varepsilon$, where $d_2$ is a large enough constant.
We already know (\cite[Proposition~5.1]{neishtadt17}) that this solution crosses the transversal $y = x$, we denote by $\tilde h_1$ the corresponding value of $\tilde h$.

Let us prove~\eqref{e:delta-h-last-wind}.
We may use coordinates $\tilde h, w, y$ to cover considered part of the phase space. Denote
\begin{equation} \label{e:def-psi}
  \psi(\tilde h, w, y) = \frac{f_{\tilde h}}{ - a(\tilde h, w)y + \varepsilon f_y}, \qquad \psi_0(w, y) = \frac{f_{\tilde h}(x{=}0, y, z{=}w)}{-a(0, w)y} = - \frac{f_x(x{=}0, y, z{=}w)}{a(0, w)}.
\end{equation}
For some $d_3, d_4 > 0$ we have
\begin{equation} \label{e:est-y-dot}
  \dot y = - a(\tilde h, w) y + \varepsilon f_y < -d_3 \varepsilon, \qquad
  \dot y = - a(\tilde h, w) y + \varepsilon f_y < - \frac{a(\tilde h, w)}{2} y
\end{equation}
for
\begin{equation} \label{e:y-greater-eps}
  y > d_4 \varepsilon.
\end{equation}
Under the condition~\eqref{e:y-greater-eps} we can express the value of $\tilde h$ along our solution as a function of $y$ satisfying the equation
\begin{equation}
  \dv{\tilde h}{y} = \varepsilon \psi(\tilde h, w, y).
\end{equation}
We also have $\Theta' = \int_{1}^{0} \psi_0(w_0, y) dy$.

Take $d_4 \varepsilon < y_4 < y_3 \le 1$. We have
\begin{equation}
  \tilde h(y_4) - \tilde h(y_3) = \varepsilon \int_{y_4}^{y_3} \psi(\tilde h, w, y) dy.
\end{equation}
As $f_{\tilde h}$ is a smooth function of $x, y$, we have
$f_{\tilde h}(\tilde h, w, y) - f_{\tilde h}(0, w, y) = O(x) = O(\tilde h/y)$.
From this we get that under the condition~\eqref{e:y-greater-eps}
\begin{equation} \label{e:diff-psi}
  \psi(\tilde h, w, y) - \psi_0(w, y) = O(\tilde h/y^2) + O(\varepsilon / y).
\end{equation}
As the right hand side of~\eqref{e:diff-psi} is $O(\varepsilon/y^2)$,
\begin{equation}
  \varepsilon \int_{y_3}^{y_4} \Big( \psi(\tilde h, w, y) - \psi_0(w, y) \Big) dy = O(\varepsilon^2 / y).
\end{equation}
We have $\norm{z - z_0} = O(\varepsilon \ln^2 \varepsilon)$ by~\eqref{e:w-close-w0} for all considered $z$, thus (by~\eqref{e:def-psi} $\psi_0(w)$ is smooth) $\psi_0(w, y) = \psi_0(w_0, y) + O(\varepsilon \ln^2 \varepsilon)$.
Take $y_3 = 1$ and $y_4 = \varepsilon^{1/2}$.
Then
\begin{equation} \label{e:diff-tilde-h-to-y2}
  \tilde h(y_2) = \tilde h_0 + \varepsilon \int_{1}^{\varepsilon^{1/2}} \psi_0(w, y) dy + O(\varepsilon^{3/2})
  = \tilde h_0 + \varepsilon \int_{1}^{\varepsilon^{1/2}} \psi_0(w_0, y) dy + O(\varepsilon^{3/2}).
\end{equation}
As $\Theta' =  \int_{1}^{0} \psi_0(w_0, y) dy$ and
\begin{equation} \label{e:int-psi-0-near-saddle}
  \int_{\varepsilon^{1/2}}^{0} |\psi_0(w, y)| dy = O(\varepsilon^{1/2}),
\end{equation}
this means
$\tilde h(y_4) > (d_2 - O(1)) \varepsilon^{3/2}$
and $x(y_4) > (d_2 - O(1)) \varepsilon$.
As $\dot x = a(\tilde h, w)x + \varepsilon f_x$ with $a > 0$, this means that for large $d_2$ we will have $x > (d_2 - O(1)) \varepsilon$ and $x$ is increasing from the moment corresponding to $y_4$ and until crossing the transversal.
Hence, for the moment of transversal crossing we have $x_1 = y_1 \ge (d_2 + O(1)) \varepsilon$.
This yields that~\eqref{e:y-greater-eps} holds until our trajectory crosses the transversal.

As $x \le y$, we have $\tilde h/y^2 = x/y \le 1$.
Thus the right hand side of~\eqref{e:diff-psi} is $O(1)$ and so
\begin{equation}
  \varepsilon \int_{\varepsilon^{1/2}}^{y_1} \Big( \psi(\tilde h, w, y) - \psi_0(w, y) \Big) dy = O(\varepsilon^{3/2}).
\end{equation}
By~\eqref{e:int-psi-0-near-saddle} this means
\begin{equation}
  \tilde h(y_1) - \tilde h(y_2) = \varepsilon \int_{y_2}^{0} \psi_0(w, y) dy + O(\varepsilon^{3/2}) = \varepsilon \int_{y_2}^{0} \psi_0(w_0, y) dy + O(\varepsilon^{3/2}).
\end{equation}
Together with~\eqref{e:diff-tilde-h-to-y2} this provides the required estimate~\eqref{e:delta-h-last-wind}.

Finally, let us estimate the time passed. As~\eqref{e:y-greater-eps} holds, we can use~\eqref{e:est-y-dot}: $\dot y < -(a/2) y$. Then $\dv{\ln y}{t} < -a/2$. As $\ln y$ changes from $O(1)$ to $O(\ln \varepsilon)$, this means that the time spent in $\mathcal B$ is $O(\ln \varepsilon)$. Similar estimate for other parts of phase space gives that the total time between two transversal crossings is also $O(\ln \varepsilon)$.
As $\dot w = O(\varepsilon)$, this also gives
$w - w_0 = O(\varepsilon \ln \varepsilon)$, yielding~\eqref{e:w-close-w0}.
\end{proof}

}
\end{appendices}

\printbibliography

\newpage

\vskip 5mm

\noindent Anatoly Neishtadt,

\noindent {\small Department of Mathematical Sciences,}

\noindent {\small Loughborough University, Loughborough LE11 3TU, United Kingdom;}

\noindent {\small Space Research Institute, Moscow 117997, Russia}

\noindent {\footnotesize{E-mail : a.neishtadt@lboro.ac.uk}}

\vskip 5mm

\noindent Alexey Okunev,

\noindent {\small Department of Mathematical Sciences,}

\noindent {\small Loughborough University, Loughborough LE11 3TU, United Kingdom}

\noindent {\footnotesize{E-mail : a.okunev@lboro.ac.uk}}

\end{document}